\documentclass[reqno,10pt]{amsart}
\usepackage {amsmath,amssymb}
\usepackage {color,graphicx}
\usepackage [latin1]{inputenc}
\usepackage {times}
\usepackage {float}
\usepackage{hyperref}
\usepackage{amscd}

\usepackage {amsmath,amssymb}
\usepackage {color,graphicx}
\usepackage [latin1]{inputenc}
\usepackage[T1]{fontenc}
\usepackage[ngerman,english]{babel}
\usepackage {times}
\usepackage {float}
\usepackage{amscd}
\usepackage{hyperref}
\usepackage{todonotes}

\DeclareMathOperator {\ev}{ev}
\DeclareMathOperator {\mult}{mult}
\DeclareMathOperator {\rc}{rc}
\DeclareMathOperator {\Star}{Star}
\DeclareMathOperator {\ind}{ind}
\DeclareMathOperator {\codim}{codim}
\DeclareMathOperator {\ft}{ft}
\DeclareMathOperator {\divisor}{div}
\DeclareMathOperator {\val}{val}

\DeclareMathOperator {\id}{id}

\newcommand {\PP}{{\mathbb P}}
\newcommand {\QQ}{{\mathbb Q}}
\newcommand {\RR}{{\mathbb R}}
\newcommand {\ZZ}{{\mathbb Z}}

\newcommand {\FF}{{\mathbb F}}
\newcommand {\Bl}{{\mathcal Bl}}
\newcommand {\calM}{{\mathcal M}}

\newcommand {\Zci}{Z^{\text{c.i.}}}

\newcommand {\M}[3]{\calM_{#1}^{\text{\upshape lab}}(\RR^{#2},#3)}

\newcommand {\dcup}{\hspace{-0.7ex}
                      \begin{array}{c}
                        \\[-3.12ex]
                        \cdot\\[-2.3ex]
                        \cup
                      \end{array}
                      \hspace{-0.7ex}
                   }

\newcommand{\Mbar}{{\overline{M}}}

\newcommand{\Fan}{{\Theta}} 

\newcommand {\bX}{{\mathbf X}}

\newcommand{\vertex}{\nu}
\newcommand{\ray}{\varrho}

\newcommand{\cell}{\tau}

\newcommand {\arxiv}[3]
            {
              #1, 
              \emph{#2},
              preprint \href{http://arxiv.org/abs/#3}{arxiv:#3}.
            }
            
\newcommand {\arxivjournal}[4]
            {
              #1, 
              \emph{#2}, 
              #3;
              also at \href{http://arxiv.org/abs/#4}{arxiv:#4}.
            }

\newtheorem {theorem}{Theorem}[section]
\newtheorem {proposition}[theorem]{Proposition}
\newtheorem {lemma}[theorem]{Lemma}

\newtheorem {corollary}[theorem]{Corollary}
\theoremstyle {definition}
\newtheorem {definition}[theorem]{Definition}

\newtheorem {notation}[theorem]{Notation}
\theoremstyle {remark}
\newtheorem {remark}[theorem]{Remark}

\begin {document}

\title {Intersections on tropical moduli spaces}
\author {Johannes Rau}
\address {\foreignlanguage{ngerman}{Johannes Rau, Universität des Saarlandes, Fachrichtung Mathematik, 
          Postfach 151150, 66041 Saarbrücken}, Germany}
\email {johannes.rau@math.uni-sb.de}


\begin{abstract}
  This article explores to which extent the
  algebro-geometric theory of rational descendant Gromov-Witten
	invariants can be carried over to the tropical world. 
	Despite the fact that the tropical moduli-spaces we work with are non-compact,
	the answer is surprisingly positive. 
  We discuss the string, divisor and dilaton equations, 
	we prove a splitting lemma describing the intersection with
  a ``boundary'' divisor and we prove general tropical versions of the 
	WDVV resp.\ topological recursion equations (under some assumptions).
	As a direct application, we prove that 
	the toric varieties $\PP^1$,  $\PP^2$, $\PP^1 \times \PP^1$ 
	and with Psi-conditions only in combination with point conditions, 
	the tropical and classical descendant Gromov-Witten invariants coincide
	(which extends the result for $\PP^2$ in \cite{MR08}). 
	Our approach uses tropical intersection theory and 
	can unify and simplify some parts of 
  the existing tropical enumerative geometry (for rational curves).
\end{abstract}

\maketitle

\section* {Introduction} \label {introduction}


Over the last few years, the list of results in tropical enumerative geometry 
became quite long. However, lacking an appropriate tropical intersection theory,
most existing results are obtained by 
\begin{itemize}
  \item
    relating the tropical numbers
    directly to the classical ones (cf.\ \cite{Mi03})
    and then using the algebro-geometric theory, or
  \item
    involved ad hoc computations (eg.\ \cite{GM05}, \cite{KM06}, \cite{FM},
    \cite{MR08}, \cite{CJM08}), which moreover have to be repeated 
    for each new class of enumerative problem.
\end{itemize}
On the other hand, based on \cite{Mi06}, the basic constructions of tropical intersection theory
are now developed in \cite{AR07}
(see also \cite{Ka2, AR08}). 
Furthermore, in \cite{GKM07} the authors show that the moduli spaces
of rational tropical curves are tropical varieties (i.e.\ satisfy the balancing condition).
Hence we can apply intersection theory to them. In \cite{Mi07}
G.~Mikhalkin proposes the definition of tropical Psi-divisors in tropical moduli spaces 
of abstract curves, and they were first studied in \cite{KM07}.
In summary, all the tools needed to develop a tropical analogue
of classical Gromov-Witten theory for rational curves are at our disposal,
and the present article tries to carry out this program 
(as mentioned before, first steps are contained e.g.\ in
\cite{GKM07}, \cite{KM07} and \cite{MR08}).

The ``ready for use'' main theorems \ref{CompOfNumbersPlane} and \ref{CompOfNumbersR=1}
state that for $\PP^1$, $\PP^2$ and $\PP^1 \times \PP^1$,
and with Psi-conditions only
in combination with point conditions, the tropical and conventional
descendant Gromov-Witten invariants coincide. For the case
of $\PP^2$, this equality was already proven in the previous paper
\cite{MR08} joint with Hannah Markwig. One should emphasize
that both results are obtained by checking that the involved numbers satisfy the same recursive formulas, 
and \emph{not} by proving some sort of correspondence theorem. 

This article is a continuation of \cite{MR08}, and some statements can be found
in older versions there. The focus here is to consequently
replace older ad hoc computations by more appropriate tools from tropical intersection theory.
As a consequence, we typically obtain more general statements 
(e.g.\ working in any dimension).

We work with non-compact tropical moduli spaces, i.e.\ fans in $\RR^N$,
due to the fact that compactifications and their intersection theory
have not yet been constructed satisfactorily. However, the non-compact
approach has its limitations. This will become visible e.g.\ from the
assumptions we need in our general WDVV and topological recursion statements. 

Let us also mention that in subsection \ref{fandisplacementrule} we show
that the fan displacement rule for Minkowski weights
describing toric intersection theory (cf.\ \cite{FS94}) coincides with the intersection
product of tropical cycles introduced in \cite[section 9]{AR07} (see also \cite{Ka2}).



\smallskip
The article contains the following parts. 
Section \ref{intersectiontheory} repeats the basics
of tropical intersection theory from \cite{AR07}
and adds some results which will be important later.
In section \ref{abstractcurves}, we study the intersection
ring of $\calM_n$, the space of abstract rational tropical curves.
In section \ref{parametrizedcurves} we extend this to
$\M{n}{r}{\Delta}$, the space of parametrized curves in $\RR^r$.
In particular, we prove general versions of
the string, dilaton and divisor equations.
Section \ref{splittingcurves} deals with the intersection of a 
one-dimensional family of curves with a boundary divisor. 
By analogy with the classical case, we prove a "`splitting lemma"'
which allows to compute this intersection as

I would like to thank Andreas Gathmann, Eric Katz, Michael Kerber,
Hannah Markwig and Grigory Mikhalkin 
for many helpful discussions and/or for proofreading various versions of this work.

\section {Intersection theory} \label {intersectiontheory}


In this section we will establish the parts of tropical intersection theory 
that we will need to attack
the problems of tropical Gromov-Witten theory in a
satisfactory way. 
Subsection \ref{basics} 
gives a quick overview on the definitions and results from \cite{AR07} and \cite{AR08}
that we will need
(however, note that our notations will sometimes
slightly differ from the original ones). 
Sections \ref{local} -- \ref{general} contain some new material. 
In particular, subsection \ref{fandisplacementrule} contains a proof of the fact that
toric intersection theory (as described by the fan displacement rule in \cite{FS94})
and tropical intersection theory for fans are identical.
(An alternative proof can be found in \cite[Theorem 4.4]{Ka2}).

\subsection{The Basics} \label{basics}
\label{cycles}


A \emph{cycle $X$} is a balanced (weighted, pure-dimensional, rational and 
polyhedral) complex in a finite-dimensional vector space 
$V = \Lambda \otimes \RR$
with underlying lattice $\Lambda$ (the most common case
is $V=\RR^r$, whose underlying lattice, if not specified otherwise,
is $\ZZ^r$). The top-dimensional polyhedra
in $X$ are called \emph{facets}, the codimension one polyhedra are
called \emph{ridges}. \emph{Balanced} means that for each ridge $\tau \in X$ the following 
\emph{balancing condition at $\tau$} is satisfied: The weighted
sum of the primitive vectors of the 
facets $\sigma$ around $\tau$
$$
  \sum_{\substack{\sigma \in X^{(\dim(X))} \\ \tau < \sigma}} 
    \omega(\sigma) v_{\sigma / \tau}
$$
vanishes ``modulo $\tau$'', or, precisely, lies in the linear vector space
spanned by $\tau$, denoted by $V_\tau$.
Here, $\omega(\sigma)$ denotes the weight of a facet $\sigma$ and
a \emph{primitive vector 
$v_{\sigma / \tau}$ of $\sigma$ modulo
$\tau$} is a vector in $\Lambda$ that points from $\tau$ towards
$\sigma$ and fulfils the primitive condition: The lattice
$\ZZ v_{\sigma / \tau} + (V_\tau \cap \Lambda)$ must be equal to the 
lattice $V_\sigma \cap \Lambda$. Slightly differently, in \cite{AR07}
the class of $v_{\sigma / \tau}$ modulo $V_\tau$ is called primitive
vector and $v_{\sigma / \tau}$ is just a representative of it. 
We will abbreviate the lattice $V_\sigma \cap \Lambda$ by $\Lambda_\sigma$. 
\\
The \emph{support of $X$}, denoted by $X$, is the union of all facets
in $X$ with non-zero weight. We call $X$ \emph{irreducible} if for any
cycle $Y$ of the same dimension with $|Y| \subseteq |X|$ there
exists an integer $\mu \in \ZZ$ such that $Y = \mu \cdot X$.
The \emph{positive part of $X$}, denoted by $X^+$, is the set of
all faces contained in a facet with positive weight. 
A \emph{general element $x$ of $X$} is an element $x \in |X|$ that lies in the interior of a facet.
If the underlying polyhedral complex is a fan (i.e.\ if all polyhedra
are actually cones with vertex in $0$), we call $X$ a \emph{fan cycle}
(or sometimes just \emph{fan}). 
\\
In fact, given a cycle $X$ we do not really want to fix its structure
as a polyhedral complex but only its support and its weights. Therefore,
by abuse of notation, a \emph{cycle} $X$ also denotes the 
class of balanced polyhedral
complexes with the same support and agreeing weights (on the common
refinement).

%
%


\smallskip
\label{principaldivisor}


A \emph{(non-zero) rational function on
$X$} is a function $\varphi : |X| \rightarrow \RR$ 
that is integer affine on each polyhedron.
Here, \emph{integer linear} means that it maps lattice elements to integers
and \emph{integer affine} means that it is a sum of an integer linear function
(called the \emph{linear part}) and a real constant. 
If $X$ is a fan, we also assume $\varphi(0) = 0$.
The \emph{divisor of $\varphi$}, denoted by
$\divisor(\varphi) = \varphi \cdot X$, is the
balanced subcomplex of $X$ constructed in \cite[3.3]{AR07}, 
namely
the codimension one skeleton $X \setminus X^{(\dim X)}$ together with
the weights $\omega_{\varphi \cdot X}(\tau)$ for each ridge
$\tau \in X$. 
These weights are given by the formula
$$
  \omega_{\varphi \cdot X}(\tau) =
    \sum_{\substack{\sigma \in X^{(\dim X)} \\ \tau < \sigma}} 
    \omega(\sigma) \varphi_\sigma(v_{\sigma / \tau}) 
    - \varphi_\tau\Big( \sum_{\substack{\sigma \in X^{(\dim X)} \\ \tau < \sigma}} 
    \omega(\sigma) v_{\sigma / \tau}\Big),
$$
where $\varphi_\sigma : V_\sigma \rightarrow \RR$ denotes the 
linear part of the affine function
$\varphi|_\sigma$. Note that the balancing condition of $X$ around
$\tau$ ensures that the argument of $\varphi_\tau$ is indeed an element
of $V_\tau$. 
Essentially, this weight measures the change of slope of $\varphi$ when traversing the
ridge $\tau$, as illustrated in the following picture.
\begin{center}
  \input{pics/DivisorConstr.pstex_t}
\end{center}
To be more precise, let $\Gamma_\varphi$ be the graph of $\varphi$ in $X \times \RR$.
It is  a polyhedral complex
whose polyhedra are in one-to-one correspondence with those of $X$, but
in general $\Gamma_\varphi$ is not balanced.
However, it can be completed to a cycle by 
adding facets in $(0,-1)$-direction at each
ridge of $\Gamma_\varphi$, equipped with the above weights. Now, if we (imaginary)
intersect this tropically completed graph of $\varphi$ with
$X \times \{-\infty\}$ (i.e. compute the tropical zero locus of $\varphi$),
we obtain the cycle $\divisor(\varphi) = \varphi \cdot X$ of our
definition. 
\\
If $\varphi$ is globally affine (resp. linear), all weights
are zero, which we denote by $\varphi \cdot X = 0$. Let the \emph{support
of $\varphi$}, denoted by  $|\varphi|$, be the subcomplex of $X$ containing 
the points $x \in |X|$ where $\varphi$ is not locally affine.
Then we have $|\varphi \cdot X| \subseteq |\varphi|$.
Furthermore, the
intersection product is bilinear (see \cite[3.6]{AR07}). 
As the restriction of a rational function to a subcycle is again a 
rational function, we can also form multiple intersection products
$\varphi_1 \cdot \ldots \cdot \varphi_l \cdot X$. In this case we will 
sometimes omit ``$\cdot X$'' to keep formulas shorter. Note that multiple
intersection products are commutative (see \cite[3.7]{AR07}).

\smallskip
\label{projection}


A \emph{morphism of cycles $X \subseteq V = \Lambda \otimes \RR$ 
and $Y \subseteq V' = \Lambda' \otimes \RR$}
is a map $f : |X| \rightarrow |Y|$ that is induced by a linear
map b $\Lambda$ to $\Lambda'$ and that maps each polyhedron 
of $X$ into a polyhedron of $Y$. We call $f$ an \emph{isomorphism} and 
write $X \cong Y$, if there exists
an inverse morphism and if for all facets $\sigma \in X$ we have
$\omega_X(\sigma) = \omega_Y(f(\sigma))$.
\\
Such a morphism \emph{pulls back rational functions} $\varphi$ on $Y$ 
to rational functions $f^*(\varphi) = \varphi \circ f$ on $X$. Note that
the second condition of a morphism makes sure that we do not have to
refine $X$ further. $f^*(\varphi)$ is already affine 
on each cone. The inclusion $|f^*(\varphi)| \subseteq
f^{-1}(|\varphi|)$ holds, as the composition of an affine and a
linear function is again affine. 
\\
Furthermore, we can \emph{push forward subcycles} 
$Z$ of $X$ to subcycles $f_*(Z)$ of $Y$ of same dimension. This is due
\cite[2.24 and 2.25]{GKM07} in the case of fans and can be generalized
to complexes (see \cite[7.3]{AR07}). We can omit further refinements
here if we assume
that $f(\sigma) \in Y$ for all $\sigma \in X$. Then $f_*(Z)$ is
defined by assigning the following weights to the
$\dim(Z)$-dimensional polyhedra $\sigma' \in Y$:
$$
  \omega_{f_*(Z)}(\sigma') =
    \sum_{\substack{\sigma \in X \\ f(\sigma) = \sigma'}}
    |\Lambda_{\sigma'} / f(\Lambda_\sigma)| \cdot \omega_Z(\sigma)
$$
By definition we have
$|f_*(Z)| \subseteq f(|Z|)$. 
The \emph{projection formula} (see \cite[4.8]{AR07}) connects all the above constructions via
$$
  f_*(f^*(\varphi) \cdot X) = \varphi \cdot f_*(X).
$$

\smallskip
\label{productofcycles}

By \cite[definition 9.3]{AR07} it is also possible to form
the intersection product of
two cycles $X,Y$ in $V = \Lambda \otimes \RR$: We choose 
coordinates $x_1, \ldots, x_r$ on $\Lambda$ (and denote the same
coordinates on the second factor of $V \times V$ by 
$y_1, \ldots, y_r$). Then the diagonal $\Delta$ in $V \times V$ is
given by $\Delta = \max\{x_1,y_1\} \cdots \max\{x_r,y_r\} 
\cdot (V \times V)$. Furthermore we consider the function
$\pi : \Delta \rightarrow V, (x,x) \mapsto x$. Then the 
intersection product of $X$ and $Y$ in $V$ is given by
$$
  X \cdot Y := \pi_*\big(\max\{x_1,y_1\} \cdots \max\{x_r,y_r\} 
    \cdot (X \times Y)\big).
$$
This intersection product is independent of the chosen coordinates,
commutative, associative, bilinear, admits the identity element 
$V$ and satisfies $(\varphi \cdot X) \cdot Y = 
\varphi \cdot (X \cdot Y)$, where $\varphi$ is a rational function
on $X$.

\smallskip
\label{rationalequivalence}


Let us now turn to the concept of rational equivalence (we summarize \cite{AR08}).
Let $X$ be a zero-dimensional cycle. Then \emph{degree $\deg(X)$
of $X$} denotes the sum of 
the weights of all points in $X$. 
Now let $X$ be an arbitrary cycle and let $\varphi,
\widetilde{\varphi}$ be two rational functions
on $X$. We call them \emph{(rationally) equivalent} if 
$\varphi-\widetilde{\varphi}$ is the sum of a bounded and a globally
linear function. Obviously, this property is preserved when
pulled back. Furthermore, if $Y$ is an one-dimensional
subcycle of $X$, then $\deg(\varphi \cdot Y) =
\deg(\widetilde{\varphi} \cdot Y)$ holds
(see \cite[lemma 8.3]{AR07}). 
\\
Let $X$ be a cycle and let $Y$ be a subcycle. We call
$Y$ \emph{rationally equivalent to zero}, denoted by $Y \sim 0$,
if there exists a morphism $f : X' \rightarrow X$ and a
bounded rational function $\phi$ on $X'$ such that
$$
  f_*(\phi \cdot X') = Y.
$$
This property commutes with
taking Cartesian products, intersection products
(of functions as well as of cycles) and with
pushing forward. Moreover, if $Y$ is zero-dimensional,
then $Y \sim 0$ implies $\deg(Y) = 0$. 
\\
Let $\widetilde{Y}$ be another subcycle of $X$. Then
we call
$Y$ and $\widetilde{Y}$ \emph{rationally equivalent} if
$Y-\widetilde{Y}$ is rationally
equivalent to zero. 
The easiest example of rationally equivalent cycles are translations.
Let $X$ be a cycle in $V = \Lambda \otimes \RR$ and let
us denote by $X + v$ denote the translation of $X$
by an arbitrary vector $v \in V$. Then
$$
  X \sim X + v
$$
holds (see also \cite[Lemma 2.1]{MR08}). 
\\
If $X,Y$ live in $V = \Lambda \otimes \RR$, we call them
\emph{numerically equivalent} if for any cycle $Z$ in $V$
of complementary dimension the equation
$$
  \deg(X \cdot Z) = \deg(Y \cdot Z)
$$
holds. 
\\
Let $X$ be a cycle in $V = \Lambda \otimes \RR$.
We define the \emph{degree} or \emph{recession fan} 
of $X$, denoted by $\delta(X)$, as follows:
$\delta(X)$ is supported on the 
purely $\dim(X)$-dimensional part of the polyhedral set 
$$
  \bigcup_{\sigma \in X} \rc(\sigma).
$$
Here, the \emph{recession cone $\rc(\sigma)$} of a polyhedron $\sigma$ 
is defined to be the cone containing all vectors $v \in V$ such that, starting 
at an arbitrary point $x \in \sigma$, the ray $x + \RR v$ is contained in $\sigma$.
Now, for a fine enough fan structure on
this polyhedral set, the weights are given by
$$
  \omega_{\delta(X)}(\sigma') :=
    \sum_{\substack{\sigma \in X \\ \sigma' \subseteq \rc(\sigma)}}
    \omega_X(\sigma).
$$
In particular, if $X$ is a curve, then $\delta(X)$ is
just the union of all unbounded rays in $X$ and 
the weights are the sums of the weights
of the rays in $X$ of given direction. Geometrically,
we simply shrink all bounded parts of $X$ to a point and move
the final single vertex to the origin. 
\\
The main result of \cite{AR08} is that for cycles 
$X$ in $V = \Lambda \otimes \RR$, rational equivalence,
numerical equivalence and ``having the same degree'' coincides.
To prove this, an important substep is to show that $X$ is always
rationally equivalent to its degree,
$$
  X \sim \delta(X).
$$

\smallskip
\subsection{Local computation of intersection products}
\label{local}


Let $X$ be a cycle and let $\tau \in X$ be a polyhedron in $X$.
We define the \emph{star of $X$ at $\tau$} to be the
fan
$$
  \Star_X(\tau) := \{\bar{\sigma} | \tau < \sigma \in X\},
$$ 
where $\bar{\sigma}$ denotes the cone in $V/V_\tau$
spanned by the image
of $\sigma - \tau$ under the quotient map 
$q : V \rightarrow V/V_\tau$. We make it into
a cycle by defining 
$\omega_{\Star_X(\tau)}(\bar{\sigma}) =
\omega_X(\sigma)$ for all facets $\bar{\sigma}$
of $\Star_X(\tau)$ (note that $q$ preserves the 
codimension of the polyhedra).
This fan
contains all the local information of $X$ around $\tau$
and can be considered as the tropical version of a small
neighbourhood of an interior point of $\tau$ (divided by
the linearity space $V_\tau$).
Its dimension equals the codimension of $\tau$
in $X$. 
\\
Let furthermore $\varphi$ be a rational function
on $X$. Choose an arbitrary affine function $\psi$
with $\varphi|_\tau = \psi|_\tau$. Then 
$\varphi - \psi$ induces a rational function on
$\Star_X(\tau)$ which we denote by $\varphi^\tau$
(and call it a \emph{germ of $\varphi$ at $\tau$}).
This function is only unique up to adding a linear function,
which is enough for us as it does not change its divisor.
The following proposition shows that our intersection products are
local constructions (i.e.\ can be expressed in terms of stars and germs).

\begin{proposition} \label{LocalIntersectionOfFunctions}
  Let $X$ be a cycle with polyhedra $\tau < \sigma \in X$.
  Let $\varphi, \varphi_1, \ldots \varphi_l$ be rational functions on $X$. Then the following
  statements are true.
  \begin{enumerate}
    \item
      $\Star_{\Star_X(\tau)}(\bar{\sigma}) = \Star_X(\sigma)$
    \item
      $(\varphi^\tau)^\sigma = \varphi^\sigma$ on $\Star_X(\sigma)$ 
      (up to adding a linear function)
    \item
      $\Star_{\varphi \cdot X}(\tau) = \varphi^\tau \cdot \Star_X(\tau)$
    \item
      $\Star_{\varphi_1 \cdot \ldots \cdot \varphi_l \cdot X}(\tau)
       = \varphi_1^\tau \cdot \ldots \cdot \varphi_l^\tau \cdot \Star_X(\tau)$
    \item
      If $l = \dim(X) - \dim(\tau)$, then 
      $\omega_{\varphi_1 \cdot \ldots \cdot \varphi_l \cdot X}(\tau)
      = \omega_{\varphi_1^\tau \cdot \ldots \cdot 
      \varphi_l^\tau \cdot \Star_X(\tau)}(\{0\})$, i.e. we can compute
      the weight of $\tau$ in $\varphi_1 \cdot \ldots \cdot \varphi_l \cdot X$
      ``locally'' in $\Star_X(\tau)$.      
  \end{enumerate}
\end{proposition}

\begin{proof}
  (a) and (b) are immediate consequences of the definitions. (d) follows from
  (c) by induction and (e) is just a special case of (d), namely
  when $\varphi_1^\tau \cdot \ldots \cdot \varphi_l^\tau \cdot \Star_X(\tau)$ 
  is zero-dimensional. Hence we are left with (c). \\
  Let $r := \dim(X) - \dim(\tau)$ be the codimension of $\tau$ in $X$. The
  statement is trivial when $r=0$: Both sides are $0$. Assume $r=1$. In this 
  case, we only have to check
  $$
    \omega_{\varphi \cdot X}(\tau) = \omega_{\varphi^\tau \cdot \Star_X(\tau)}(\{0\}).
  $$
  By adding an affine function we can assume that $\varphi|_\tau = 0$ without
  changing the intersection product and in particular the weight of $\tau$ in 
  $\varphi \cdot X$. But then we can replace both weights according to their
  definition and observe that
  $$
    \omega_{\varphi \cdot X}(\tau) =
      \sum_{\substack{\sigma \in X^{(\dim(X))} \\ \tau < \sigma}} 
      \omega(\sigma) \varphi_\sigma(v_{\sigma / \tau}) =
      \sum_{\bar{\sigma} \in \Star_X(\tau)^{(1)}} 
        \omega(\bar{\sigma}) \varphi^\tau(v_{\bar{\sigma} / \{0\}}) =
    \omega_{\varphi^\tau \cdot \Star_X(\tau)}(\{0\}) 
  $$
  holds true, as 
  $[v_{\sigma / \tau}] = v_{\bar{\sigma} / \{0\}}
  \in V/V_\tau$. \\
  Now let us assume $r > 1$ and let $\tau'$ be a ridge in $X$. Then
  we can use the previous case as well as (a) and (b) and obtain
  $$
    \omega_{\varphi \cdot X}(\tau') \stackrel{r=1}{=}
      \omega_{\varphi^{\tau'} \cdot \Star_X(\tau')}(\{0\})
        \stackrel{\text{(a), (b)}}{=}
      \omega_{(\varphi^\tau)^{\tau'} \cdot 
        \Star_{\Star_X(\tau)}(\tau')}(\{0\}) \stackrel{r=1}{=}
    \omega_{\varphi^\tau \cdot \Star_X(\tau)}(\bar{\tau}'),   
  $$
  which proves the claim.
\end{proof}

We can extend this to the case of the intersection product of two
cycles.

\begin{proposition} \label{LocalIntersectionOfCycles}
  Let $X, Y$ be two cycles in $V = \RR \otimes \Lambda$.
  Then the equation
  $$
    \Star_{X \cdot Y}(\tau) = \Star_X(\tau) \cdot \Star_Y(\tau).
  $$
  holds for all polyhedra $\tau \in X \cdot Y$.
\end{proposition}

\begin{proof}
  First, we fix some notation. Let $x_1, \ldots, x_r$
  be a lattice basis of $\Lambda^\vee$ such that the first
  $d := \codim_V(\tau)$ elements generate $V_\tau^\perp$.
  When we consider the product $\Lambda \times \Lambda$, the 
  same coordinates on the second factor will be denoted
  by $y_1, \ldots, y_r$.
  Furthermore, let $\Delta : V \rightarrow V \times V,
  x \mapsto (x,x)$ denote the diagonal map.
  By definition of the intersection product of cycles and
  using \ref{LocalIntersectionOfFunctions} (d) we have to compute
  $$
    \Star_{\max\{x_1,y_1\} \cdots \max\{x_r,y_r\} 
          \cdot (X \times Y)}(\Delta(\tau))
    = \max\{x_1,y_1\} \cdots \max\{x_r,y_r\} \cdot
      \Star_{X \times Y}(\Delta(\tau)) 
  $$
  and 
  $$
    \max\{x_1,y_1\} \cdots \max\{x_d,y_d\} \cdot
      (\Star_X(\tau) \times \Star_Y(\tau))
  $$
  respectively.
  Thus the statement follows from the fact that
  \begin{eqnarray*}
    \max\{x_{d+1},y_{d+1}\} \cdots \max\{x_r,y_r\} \cdot 
        (V \times V / \Delta(V_\tau))
      & \rightarrow
      & V / V_\tau \times V / V_\tau, \\
    (x, y)
      & \mapsto
      & (x, y)
  \end{eqnarray*}
  is an isomorphism and can be restricted
  to an isomorphism of 
  $\max\{x_{d+1},y_{d+1}\} \cdots \max\{x_r,y_r\} \cdot
  \Star_{X \times Y}(\Delta(\tau))$
  and $\Star_X(\tau) \times \Star_Y(\tau)$.
\end{proof}

\subsection{Transversal Intersections}
\label{transversal}


Let us now consider ``generic'' intersections.

\begin{definition}
  Let $X, Y$ be two cycles in $V = \Lambda \otimes \RR$ of 
  codimension $c$ resp. $d$. We say
  \emph{$X$ and $Y$ intersect transversally} if 
  $X \cap Y$ is of pure codimension $c + d$ and if
  for each facet $\tau$ in $X \cap Y$ the corresponding
  neighbourhoods $\Star_X(\tau)$ and $\Star_Y(\tau)$ are 
  (transversal) affine subspaces of $V$.
\end{definition}

In this case, by locality of the intersection product, the 
computation of $X \cdot Y$ can be reduced to the intersection
of vector spaces. This motivates the following study of intersections of linear
functions and spaces.

\begin{lemma} \label{IntersectingLinearFunctions}
  Let $h_1, \ldots, h_l$ be integer linear functions on $V$ ($l \leq \dim(V) =: r$) 
  and define the rational functions $\varphi_i := \max\{h_i,0\}$
  on $V$. Let $H : V \rightarrow \RR^l$ be the linear function
  with $H(x) = (h_1(x), \ldots, h_l(x))$ and let us assume
  that $H$ has full rank. Then 
  $\varphi_1 \cdot \ldots \cdot \varphi_l \cdot V$ 
  is equal to the 
  subspace $\ker(H)$ with weight $\ind(H) := |\ZZ^l/H(\Lambda)|$.
  Here we give $V$ the fan structure consisting of all 
  cones where each of the $h_i$ is either positive or zero
  or negative, with all weights being $1$.
\end{lemma} 

\begin{proof}
  Let us assume $l=1$ first (i.e. $H=h_1$)
  In this case we have
  to compute the weight of the only ridge in $V$
  which is $h_1^\bot = \ker(H)$. This ridge is contained
  in the two facets corresponding to $h_i \geq 0$ and
  $h_i \leq 0$. Let $v_\geq = - v_\leq$ be  
  corresponding primitive vectors. This implies
  that for example $v_\geq$ generates the
  one-dimensional lattice $\Lambda/h_1^\bot \cong
  h_1(\Lambda)$ and therefore 
  $|\ZZ / h_1(\Lambda)| = h_1(v_\geq)$.
  On the other hand we can compute the weight of
  $h_1^\bot$ in $h_1 \cdot V$ to be
  $$
    \omega_{h_1 \cdot V}(h_1^\bot)
      = \varphi_1(v_\geq) + \varphi_1(v_\leq)
      = h_1(v_\geq) + 0 = |\ZZ / h_1(\Lambda)|.
  $$
  Now we make induction for $l > 1$. 
  The induction hypothesis says that
  $\varphi_2 \cdot \ldots \varphi_l \cdot V$
  is equal to the subspace $\ker(H')$ with weight
  $\ind(H')$, where $H' = h_2 \times \ldots \times h_l$.
  By applying the case $l=1$ to the vector space
  $\ker(H') = (\ker(H') \cap \ZZ^r) \otimes \RR$, we
  obtain that $\varphi_1 \cdot \ldots \varphi_l \cdot V$
  is equal to the subspace $h_1^\bot \cap \ker(H') = \ker(H)$ 
  with weight $\ind(h_1|_{\ker(H')}) \cdot \ind(H')$. We have to
  show that this weight coincides with $\ind(H)$. This follows
  from the exact sequence
  $$
    \begin{array}{ccccccccc}
      0 & \rightarrow & h_1(\ker(H') \cap \Lambda)
        & \rightarrow & H(\Lambda) & \rightarrow
        & H'(\Lambda) & \rightarrow & 0 \\
      & & h_1(x) & \mapsto & H(x) = (h_1(x), 0) \\
      & & & & H(x) & \mapsto & H'(x) 
    \end{array}
  $$
  and its induced quotient sequence
  $$
    \begin{array}{ccccccccc}
      0 & \rightarrow & \ZZ^{l-1}/H'(\Lambda) 
        & \rightarrow & \ZZ^l / H(\Lambda) & \rightarrow
        & \ZZ / h_1(\ker(H') \cap \Lambda) & \rightarrow & 0 
    \end{array}.
  $$
\end{proof}

\begin{remark}
  In the special case $l = r$ the weight of 
  $\{0\}$ in the intersection product
  $\varphi_1 \cdot \ldots \cdot \varphi_r \cdot V$
  is $|\ZZ^r / H(\Lambda)|$, which equals $|\det(M)|$
  where $M$ is a matrix representation of $H$
  with respect to a lattice basis of $\Lambda$ and
  the standard basis of $\ZZ^r$. This special case
  of the statement is proven in \cite[Lemma 5.1]{MR08}. 
	Note that in this case, if $\det(M)$ is zero,
	the intersection product is zero as well.
	Hence this version can be extended to
  the case where $H$ has not full rank.
\end{remark}

Now we use this lemma to compute the intersection of 
two linear subspaces.

\begin{lemma} \label{IntersectingLinearSpaces}
  Let $U,W$ be two subspaces of $V = \RR \otimes \Lambda$ (with
  rational slope) such that $U + W = V$. If we consider
  $U,W$ as cycles with weight $1$, their intersection product
  can be computed to be
  $$
    U \cdot W = |\Lambda / (\Lambda_U + \Lambda_W)| \cdot (U \cap W).
  $$
\end{lemma}

\begin{proof}
  By definition we have to compute
  $$
    \max\{x_1,y_1\} \cdots \max\{x_r,y_r\} \cdot (U \times W),
  $$
  (where we chose arbitrary coordinates on $\Lambda$).
  Instead of $\max\{x_i,y_i\}$, we can as well subtract the linear
  function
  $y_i$ and use the functions $\max\{x_i-y_i,0\}$. Now we can apply
  \ref{IntersectingLinearFunctions}. In our case, the function
  $H$ is just
  \begin{eqnarray*}
    H : \Lambda \times \Lambda & \rightarrow & \Lambda, \\
        (x, y)                 & \mapsto     & x - y.
  \end{eqnarray*}
  Restricted to $U \times W$, this provides
  $$
    U \cdot W
      = |\Lambda / H(\Lambda_U \times \Lambda_W)| \cdot \pi_*(\ker(H))
      = |\Lambda / (\Lambda_U \mp \Lambda_W)| \cdot (U \cap W).
  $$
\end{proof}

Now, as a combination of \ref{LocalIntersectionOfCycles} 
and \ref{IntersectingLinearSpaces}, we obtain the
following result.

\begin{corollary} \label{TransversalIntersectionOfCycles}
  Let $X, Y$ be two cycles in $V = \RR \otimes \Lambda$
  that intersect transversally. Then
  $X \cdot Y = (X \cap Y, \omega_{X \cap Y})$ with the
  following weight function. Any facet $\tau$ in
  $X \cap Y$ is the intersection of two facets $\sigma, \sigma'$
  in $X$ resp. $Y$. Then the weight of $\tau = \sigma \cap \sigma'$ is
  $$
    \omega_{X \cap Y}(\sigma \cap \sigma') 
      = \omega_X(\sigma) \omega_Y(\sigma') |\Lambda / \Lambda_\sigma + \Lambda_{\sigma'}|.
  $$
\end{corollary}

\subsection{Comparison to the ``fan displacement rule''}
\label{fandisplacementrule}


In \cite{FS94} the authors introduce Minkowski weights to
describe the Chow cohomology groups of a toric variety
combinatorially. Moreover, they compute the cup-product
of these cohomology groups in terms of Minkowski weights.
In this subsection we show explicitly that, when we interpret
Minkowski weights as tropical cycles, this cup-product
coincides with our product of tropical cycles.
Another approach to this topic is given in \cite[section 9]{Katz06}
and \cite{Ka2}. 
\\ 
Let $\Theta$ be a complete fan in a vector space 
$V = \RR \otimes \Lambda$ of dimension $r$
(in \cite{FS94}, the fan is called $\Delta$ and the lattice
is called $N$). Let $\Theta^{(k)}$ denote the set of
$k$-dimensional cones in $\Theta$ (in \cite{FS94}, the
exponent indicates the codimension, i.e. $\Delta^{(k)}$
means $\Theta^{(r-k)}$).

\begin{definition}[cf. \cite{FS94}, section 2]
  A \emph{Minkowski weight $c$ of codimension $k$} is an integer-valued function
  on $\Theta^{(r-k)}$ that satisfies for any $\tau \in \Theta^{(r-k-1)}$
  $$
    \sum_{\substack{\sigma \in \Theta^{(r-k)} \\
                    \tau \subseteq \sigma}}
      c(\sigma) v_{\sigma/\tau}
      \in \Lambda_\tau
  $$
  (in \cite{FS94}, primitive vectors are denoted
  by $n_{\sigma, \tau}$).
\end{definition}

Let $c$ be a Minkowski weight of codimension $k$. 
Of course, if we set $X(c)$ to be
the fan $\bigcup_{0 \leq i \leq r-k} \Theta^{(i)}$ with
weight function $c$, the Minkowski weight condition
precisely coincides with our balancing condition, i.e.
$X(c)$ is a tropical cycle of codimension $k$.
\\
In \cite{FS94} it is shown that Minkowski weights
are in one-to-one correspondence with
the Chow cohomology classes
of the toric variety associated to the fan $\Theta$
and therefore admit a cup-product with the following properties.
Let $c, c'$ be Minkowski weights of codimension $k, k'$. Then
the cup-product $c \cup c'$ is a Minkowski weight of codimension 
$k + k'$ given by
$$
  (c \cup c')(\tau) =
    \sum_{\substack{\sigma \in \Theta^{r-k} \\
                    \sigma' \in \Theta^{r-k'} \\
                    \tau \subseteq \sigma, \sigma'}}
                    m_{\sigma, \sigma'}^{\tau} \cdot 
                    c(\sigma) \cdot c'(\sigma').
$$
Here, the coefficients are not unique but depend on the choice
of a generic vector $v \in V$. If we fix such a vector $v$, then
$$
  m_{\sigma, \sigma'}^{\tau} =
    \begin{cases}
      |\Lambda / \Lambda_\sigma + \Lambda_{\sigma'}| & 
         \text{if } (\sigma + v) \cap \sigma' \neq \emptyset, \\
      0 & \text{otherwise}
    \end{cases}
$$
(cf. \cite[introduction]{FS94}).
With the tools developed in the previous sections, we can show easily (and purely tropically)
that the cup-product of Minkowski weights
coincides with our intersection product of tropical cycles
in $V$. An independent proof of this statement is given in \cite[Theorem 4.4]{Ka2}.

\begin{theorem} \label{FanDisplacementEqualsTropInt}
  Let $c, c'$ be Minkowski weights of codimension $k, k'$. 
  Then the following equation holds.
  $$
    X(c) \cdot X(c') = X(c \cup c')
  $$
\end{theorem}

\begin{proof}
  For each facet $\tau$ in $X(c \cup c')$ we have to
  show
  $$
    \omega_{X(c) \cdot X(c')}(\tau) = (c \cup c')(\tau).
  $$
  First, note that we can compute both sides locally
  on $\Star_{\Theta}(\tau)$, where we of course define
  the ``local'' Minkowski weights by $\bar{c}(\bar{\sigma}) := c(\sigma)$
  and $\bar{c}'(\bar{\sigma}') := c'(\sigma')$.
  For the left hand side this follows from \ref{LocalIntersectionOfCycles}
  and for the right hand side it follows from
  $|\Lambda / \Lambda_\sigma + \Lambda_{\sigma'}|
  = |(\Lambda / \Lambda_\tau) / ((\Lambda_\sigma + \Lambda_{\sigma'})/\Lambda_\tau)|$. \\
  Therefore we can assume $k + k' = r$ and
  $\tau = \{0\}$. In this case, by plugging
  in the definition on the right hand side and choosing
  a generic vector $v \in V$, it remains to show
  $$
    \deg(X(c) \cdot X(c')) 
      = \sum_{\substack{\sigma \in \Theta^{r-k} \\
                    \sigma' \in \Theta^{r-k'} \\
                    (\sigma + v) \cap \sigma' \neq \emptyset}}
                    |\Lambda / \Lambda_\sigma + \Lambda_{\sigma'}| \cdot 
                    c(\sigma) \cdot c'(\sigma').
  $$
  Now, for a generic vector $v \in V$ we can assume that
  $X(c) + v$ and $X(c')$ intersect transversally (in fact, this
  is what the authors of \cite{FS94}
  mean by a generic vector). Note that, in fact, the sum
  on the right hand side runs through all points in the 
  intersection of
  $X(c) + v$ and $X(c')$. Therefore, by \ref{TransversalIntersectionOfCycles}
  it equals $\deg((X(c) + v) \cdot X(c'))$. But as $X(c) + v$ and $X(c)$
  are rationally equivalent, the equation
  $\deg(X(c) \cdot X(c')) = \deg((X(c) + v) \cdot X(c'))$ holds
  and the statement follows. 
\end{proof}

\subsection {Convexity and Positivity}
\label{convexpositive}


A non-zero cycle $X$ is called \emph{positive}, denoted $X > 0$, if
all weights 
are non-negative. By throwing away the facets with weight $0$ (and all polyhedra
contained in only such facets) we can assume all weights to be positive. 
A rational function $\varphi$ on $X$ is called \emph{convex} if it is locally
the restriction of a convex function on $V$.
The pull-back $f^*(\varphi)$ of a convex function is again convex, as the
composition of a convex function and a linear map is again convex.
Moreover, if $Z$ is a subcycle of $X$, then $\varphi|_{|Z|}$ is also
convex on $Z$.
Combining positivity and convexity we get the following result.

\begin{proposition} \label{convex}
  Let $X$ be a positive cycle and let $\varphi$ be a convex function
  on $X$. Then
  \begin{enumerate}
    \item
      $\varphi \cdot X$ is positive and
    \item
      $|\varphi| = |\varphi \cdot X|$.
  \end{enumerate}
\end{proposition}

\begin{proof}
  First of all note that we can assume that $X$ is a one-dimensional fan,
  as all intersection weights can be computed locally modulo the ridge
  (cf. \ref{LocalIntersectionOfFunctions} (c)) and convexity is preserved
  when adding linear functions or when considering the function induced
  on the quotient. Thus we assume that $X = \{\{0\}, \rho_1, \ldots,
  \rho_r\}$ is a one-dimensional fan with positive
  weights $\omega(\rho_i) > 0$. The statements of the lemma translate to
  \begin{enumerate}
    \item
      $\varphi$ convex $\Rightarrow$ $\varphi \cdot X > 0$, 
    \item
      $\varphi$ convex, $\varphi \cdot X = 0$ $\Rightarrow$ $\varphi$ linear.
  \end{enumerate}
  We use the following criteria for linearity and convexity. Let $\varphi$
  be a rational function on $X$ and let us abbreviate the primitive
  vector of the ray $\rho_i$ by $v_i$. Then
  \begin{enumerate}
    \item[i)]
      $\varphi$ is linear if and only if for all 
      $\lambda_1, \ldots, \lambda_r \in \RR$ with
      $\sum_i \lambda_i v_i = 0$ it holds 
      $$\sum_i \lambda_i \varphi(v_i) = 0,$$
    \item[ii)]
      $\varphi$ is convex if and only if for all positive 
      $\lambda_1, \ldots, \lambda_r \geq 0$ with
      $\sum_i \lambda_i v_i = 0$ it holds 
      $$\sum_i \lambda_i \varphi(v_i) \geq 0.$$
  \end{enumerate}
  Now let $\varphi$ be convex. We can apply criterion
  ii) to the coefficients $\omega(\rho_i)$, which are
  positive and satisfy $\sum_i \omega(\rho_i)
  v_i = 0$. This provides 
  $$
    \omega_{\varphi \cdot X}(\{0\}) =
      \sum_i \omega(\rho_i) \varphi(v_i) \geq 0,
  $$
  which proves (a). \\
  For (b), let us assume that $\sum_i \omega(\rho_i) \varphi(v_i) = 0$
  (i.e. $\varphi \cdot X = 0$) but $\varphi$ is not linear. Then by i)
  there exist $\lambda_1, \ldots, \lambda_r$ with
  $\sum_i \lambda_i v_i = 0$ but $\sum_i \lambda_i \varphi(v_i) \neq 0$.
  W.l.o.g.\ we can assume $\sum_i \lambda_i \varphi(v_i) < 0$ (otherwise
  we replace $\lambda_i$ by $-\lambda_i$). For large enough $C \in \RR$ 
  the coefficients $\lambda'_i := \lambda_i + C \omega(\rho_i)$ are all
  positive and still satisfy $\sum_i \lambda'_i v_i = 0$ and
  $\sum_i \lambda'_i \varphi(v_i) < 0$, which contradicts ii). Therefore
  $\varphi$ is linear, which proves $(b)$.
\end{proof}

The following application of this proposition we be useful for us later.

\begin{proposition} \label{ImageOfConvexPullbacks}
  Let $f : X \rightarrow Y$ be a morphism of cycles and let
  us assume that $Y$ is positive. Let furthermore
  $\varphi_1, \ldots, \varphi_l$ denote convex functions
  on $Y$. Then the following equation of sets holds.
  $$
    |f^*(\varphi_1) \cdots f^*(\varphi_l) \cdot X|
      \subseteq f^{-1}(|\varphi_1 \cdots \varphi_l \cdot Y|)
  $$
\end{proposition}

\begin{proof} 
  This can be proven by an easy induction. If $l=1$ we have
  $$
    |f^*(\varphi_1) \cdot X| =
      |f^*(\varphi_1)| \subseteq
      f^{-1}(|\varphi_1|) =
      f^{-1}(|\varphi_1 \cdot Y|),
  $$
  where the equalities follow from
  \ref{convex} (a).
  Now for arbitrary $l$ we can apply
  the case of a single function to $\varphi_l$,
  obtaining 
  $$
    |f^*(\varphi_l) \cdot X| \subseteq f^{-1}(|\varphi_l \cdot Y|).
  $$
  This shows that we can restrict the morphism $f$ to
  $f : f^*(\varphi_l) \cdot X \rightarrow \varphi_l \cdot Y$.
  As $\varphi_l \cdot Y$ is still positive by
  \ref{convex} (b), we
  can apply the induction hypothesis to this restriction,
  which yields the result.
\end{proof}

\subsection {Complete intersections} \label{completeintersection}
\label{complete}


We define the set of $m$-dimensional \emph{complete
intersections} $\Zci_m(X) \subset Z_m(X)$ to be the set of $m$-dimensional
cycles in $X$ obtained as an intersection product $\varphi_1 \cdots
\varphi_l \cdot X$ (where $l = \dim(X) - m$). 
\\
  Let $C, C' \in \Zci_*(X)$ be complete intersections
  given by $C = \varphi_1 \cdots
  \varphi_l \cdot X$ and 
  $C' = \varphi'_1 \cdots
  \varphi'_{l'} \cdot X$. Then we define
  $$
    C \cdot C' := \varphi_1 \cdots 
      \varphi_l \cdot
      \varphi'_1 \cdots
      \varphi'_{l'} \cdot X.
  $$
  Using commutativity of the intersection product of functions,
  this multiplication is independent of the chosen functions,
  commutative and satisfies $|C \cdot C'| = |C| \cap |C'|$.
  Note that, if $X = V = \Lambda \otimes \RR$, it follows from
  \cite[corollary 9.8]{AR07} that this definition
  coincides with the usual intersection
  product of cycles. 
	\\
  Let $C \in \Zci_m(X)$ be given by $C = \varphi_1 \cdots
  \varphi_l \cdot X$ and let $f : Y \rightarrow X$ be a
  tropical morphism. Then we would like to define
  the pull-back of $C$ along $f$ to be
  the complete intersection
  $$
    f^*(C) :=
      f^*(\varphi_1) \cdots f^*(\varphi_l) \cdot Y.
  $$
  However, in general this definition is not independent of
  the chosen functions $\varphi_1, \ldots, \varphi_l$. 
  For us it is enough to consider the case of projections
	where this indeterminacy does not occur.
	
\begin{proposition} \label{PullBackUnderProjection}
  Let $X,Y$ be two cycles and let $\pi : X \times Y \rightarrow X$
  be the projection onto the first factor. Moreover, let $Z$ be 
  a complete intersection of $X \times Y$ and consider the map
  $f = \pi|_Z : Z \rightarrow X$. Now, if $C = \varphi_1 \cdots
  \varphi_l \cdot X$ is a complete intersection in
  $X$, then the pull-back 
  $$
    f^*(C) :=
      f^*(\varphi_1) \cdots f^*(\varphi_l) \cdot Z
  $$
  is well-defined and the equation
  $$
    |f^*(C)| \subseteq f^{-1}(|C|)
  $$
  holds.
\end{proposition}

\begin{proof}
  First, we apply \cite[9.6]{AR07}, which yields
  $$
    \pi^*(\varphi_1) \cdots \pi^*(\varphi_l) 
      \cdot (X \times Y)
      = (\varphi_1 \cdots \varphi_l \cdot X) \times Y
      = C \times Y.
  $$
  Therefore $f^*(\varphi_1) \cdots f^*(\varphi_l) \cdot Z$
  is just the product of the complete intersections
  $C \times Y$ and $Z$, which does not depend on any choices.
  Moreover, its support is contained in $|C \times Y|$ and 
  the equation of sets follows.
\end{proof}

\begin{remark}[Pulling back preserves numerical equivalence]
               \label{PullBackPreservesNumericalEquivalence}
  Let $C, C'$ be 
  complete intersections in $\RR^r$
  and 
  let $f : Y \rightarrow \RR^r$ be a
  tropical morphism. Then, if $C$ and $C'$ are
  numerically equivalent, also $f^*(C)$ and $f^*(C')$
  are numerically equivalent in the following
  sense. If $Z$ is an arbitrary complete intersection
  in $Y$ of complementary dimension, then
  $$
    \deg(f^*(C) \cdot Z) = \deg(f^*(C') \cdot Z)
  $$
  holds.
  This follows from the projection formula.
  $$
    \deg(f^*(C) \cdot Z) =
      \deg(f_*(f^*(C) \cdot Z)) =
      \deg(C \cdot f_*(Z))
  $$
  In particular, if we move around $C$ in $V$, the 
  numerical properties of the pull-backs of the original
  and the translated cycle coincide. 
\end{remark}

\subsection {General position} \label{generalposition}
\label{general}


We now investigate what we can say about the set-theoretic
preimage of a general translation of a cycle under a morphism $f$.
This section is a simple generalization of \cite[Section 3]{MR08}
(where $f$ is (a combination of) evaluation morphisms).

\begin{lemma} \label{generalposition1}
  Let $X$ be a pure-dimensional polyhedral complex and let $f : X \rightarrow \RR^r$ 
  be a morphism of polyhedral complexes (i.e. $f$ is linear
  on every polyhedron of $X$). Furthermore, let 
  $C$ be a polyhedral complex in $\RR^r$ and consider the subcomplex
  $f^{-1}(C)$ of $X$ consisting of
  all polyhedra $\tau \cap f^{-1}(\gamma), \tau \in X, \gamma \in C$.
  Then for a general translation $C' = C + v$ (i.e. $v \in \RR^r$ can
  be chosen from an open dense subset of $\RR^r$) 
      the codimension of each non-empty polyhedron $\tau \cap f^{-1}(\gamma)$
      of $X$ is equal to 
      $$
        \codim_X(\tau \cap f^{-1}(\gamma)) =
          \codim_X(\tau) + \codim_{\RR^r}(\gamma).
      $$
\end{lemma}

\begin{proof}
  For each $\tau$ in $X$ and $\gamma$ in $C$ we consider the 
  map
  $$
    f_\tau : \text{AffSpan}(\tau) \rightarrow \RR^r,
  $$
  induced by $f|_\tau$. Now we are interested
  in $\tau \cap f^{-1}(\gamma') = \tau \cap f_\tau^{-1}(\gamma')$ 
  for general translations $\gamma'$ of $\gamma$. 
  We have to distinguish the cases
  $\text{Im}(f_\tau) + V_\gamma = \RR^r$ and $\text{Im}(f_\tau)
  + V_\gamma \neq \RR^r$. In the latter case, $f_\tau^{-1}(\gamma')$
  is empty for general $\gamma'$. In the former case,
  $f_\tau^{-1}(\gamma')$ is a polyhedron of dimension
  $\dim(\tau) + \dim(\gamma) - r$, and for general $\gamma'$
  it is disjoint from $\tau$ or intersects the interior of $\tau$,
  in which case $\tau \cap f_\tau^{-1}(\gamma')$ has the same
  dimension $\dim(\tau) - \codim_{\RR^r}(\gamma)$, which is 
  the expected dimension. \\
  As there are only finitely many pairs
  $\tau, \gamma$, this holds simultaneously for all pairs for
  general enough translations of $C$.
\end{proof}

This technical statement has the following more applicable consequences.

\begin{proposition}[Preimages of general translations] \label{generalposition2}
  Let $f_k : X \rightarrow \RR^r, k = 1, \ldots, n$ be morphisms of pure-dimensional 
  polyhedral complexes and let 
  $C_k, k = 1, \ldots, n$ be cycles in $\RR^r$. 
  Then for a general translation $C'_k = C_k + v_k, v_k \in \RR^r$ the following
  holds. Either
  $Z := f_1^{-1}(C'_1) \cap \ldots \cap f_n^{-1}(C'_n)$ is empty or
  \begin{enumerate}
    \item 
      the codimension of $Z$
       in $X$ equals the sum 
       $$
         \codim_X(Z) = \sum_{k=1}^n \codim_{\RR^r}(C_k),
       $$
    \item
      $Z$ is pure-dimensional,
    \item 
      if a polyhedron $\alpha$ of $Z$ is
      contained in a polyhedron $\tau$ of $X$, the codimensions
      satisfy \linebreak $\codim_X(\tau) \leq \codim_{Z}(\alpha)$
      (in particular, the interior of a facet of 
      $Z$ is
      contained in the interior of a facet of $X$),
    \item
      if the images $f_k(\alpha)$ of a polyhedron $\alpha$ of 
      $Z$ are contained
      in polyhedra $\gamma_k$ of $C_k$, the codimensions satisfy
      $\sum_{k=1}^n \codim_{C_k}(\gamma_k) \leq 
      \codim_{Z}(\alpha)$.
  \end{enumerate}
\end{proposition}

\begin{proof}
  It is easy to prove the statement in the case $n=1$:
  (a), (b) and (c) are immediate consequences of
  \ref{generalposition1} and (d) follows from
  applying \ref{generalposition1} to the
  $(r-\codim_{Z}(\alpha)-1)$-dimensional skeleton of $C_1$ (if $\gamma_1$ belonged
  to this skeleton, $\alpha$ would be contained in its preimage, which 
  (for general translations) contradicts (a)).
  Now the statement follows if we apply 
  the case of a single morphism
  to $f_1 \times \ldots \times f_n : X \rightarrow (\RR^r)^n$
  and $C := C_1 \times \ldots \times C_n$. 
\end{proof}

\begin{remark} \label{ComparisonOfPullBackAndPreimage}
  Sticking to the notations of the previous statement, let us assume
  that $X$ is a cycle and that the maps $f_k$ are tropical morphisms.
  Moreover, we assume that the maps $f_k$ are projections
  (at least after composing with an isomorphism) and that the complexes 
  $C_k$ are complete intersections. Then $f_1^*(C_1) \cdots f_n^*(C_n)$ is
  also a pure-dimensional complex of the same dimension as 
  $f_1^{-1}(C'_1) \cap \ldots \cap f_n^{-1}(C'_n)$. 
	Moreover, \ref{PullBackUnderProjection} shows that
  $$
    |f_1^*(C_1) \cdots f_n^*(C_n)| \subseteq 
      f_1^{-1}(C'_1) \cap \ldots \cap f_n^{-1}(C'_n)
  $$
  holds. 
	Hence in this case we can think of $f_1^*(C_1) \cdots f_n^*(C_n)$
  as being the polyhedral set
  $f_1^{-1}(C'_1) \cap \ldots \cap f_n^{-1}(C'_n)$ with the
  additional data of weights (some of which might be zero). 
\end{remark}

\section {Intersections on the space of abstract curves} \label {abstractcurves}

 
Let us start with a definition of smooth abstract curves.
As a local model of a curve we will use the following fan.
Let $e_1, \ldots, e_r$ be the standard basis in $\RR^r$ and
set $e_0 := -e_1 - \ldots - e_r$. We define
the one-dimensional fan
$$
  L^r := \{\{0\}, \RR_{\geq} (-e_0), \ldots, \RR_{\geq} (-e_r)\},
$$
with weights $\omega(\RR_{\geq} (-e_i)) = 1$ for all $i$.
This fan is balanced because of $e_0 + \ldots + e_r = 0$.
\begin{center}
  \input{pics/SmoothVertices.pstex_t}
\end{center}
Note that this fan is also irreducible, as $e_0 + \ldots + e_r = 0$
is the only relation that the generating vectors fulfil.

\begin{definition} \label{DefAbstractCurves}
  A \emph{smooth abstract curve $C$} is a one-dimensional connected cycle
  that is locally isomorphic to $L^r$ for suitable $r$, i.e. for each
  vertex $V$ in $C$ we have $\Star_C(V) \cong L^{val(V)}$. The 
  \emph{genus of $C$} is the first Betti number of $|C|$. 
  An \emph{$n$-marked smooth abstract curve $(C,x_1, \ldots, x_n)$}
  is a smooth abstract curve $C$ with $n$ unbounded rays (called
  \emph{leaves}), which are labelled by $x_1, \ldots x_n$. If we
  instead label
  the leaves by elements of some finite set $I$, we will call it an
  $I$-marked curve.
\end{definition}

\begin{remark}
	We will often omit the word ``smooth'' here
	as we will not consider other abstract curves 
	(which are allowed to have different one-dimensional fans
	as local structures).
  Note that
  by definition $C$ is (locally) irreducible. 
  We will always consider
  abstract curves up to isomorphisms.
  \\
  Note that the valence of a vertex $V$ in $C$ completely fixes the local structure
	(which is $L^{val(V)}$). 
	Hence $C$ is in fact completely determined by the underlying metric graph, i.e.\
	the combinatorial graph together with the (lattice length) of the edges. 
	This is the definition in most existing literature,
  in particular in \cite{GKM07}.
	We will later be interested in parametrized curves, i.e.\
	maps $f : C \rightarrow \RR^r$. 
	With the ``metric graph'' definition, the balancing condition has to be incorporated
	in the definition of these maps (see \cite[4.1]{GKM07}).
	With our definition we can just impose
	that $f$ should be a morphism of tropical cycles, giving the same result.
	Note that our definition also requires
  that a global embedding $C \subset \RR^N$ of our curve
  exists (which we then forget as we identify isomorphic curves). 
  This is done to avoid some technicalities involved in glueing 
	abstract tropical cycles.
  However, will see that (at least for rational curves) this
	is not a restriction as any ``metric graph'' curve can be embedded. 
\end{remark}

\begin{remark}[Smoothness criterion] \label{Smoothness}
Let us mention two simple criteria to decide whether a one-dimen\-sio\-nal 
fan with $r+1$ rays is isomorphic to $L^r$ or not (i.e.\ smoothness
criteria).
\\
  Let $X$ be a one dimensional fan in $V = \Lambda \otimes \RR$
  with $r+1$ rays, all with weight $1$
  and generated by the primitive vectors $v_0, \ldots, v_r$. 
	Let $V_X$ be the vector space spanned by $X$.
  Then the following are equivalent.
  \begin{enumerate}
    \item
      $X$ is isomorphic to $L^r$.
    \item
      The equations $v_0 + \ldots + v_r = 0$, $\dim(V_X) = r$
      and $V_X \cap \Lambda = \ZZ v_0 + \ldots + \ZZ v_r$ hold.
    \item
      For arbitrary coefficients 
      $\lambda_0, \ldots, \lambda_r \in \RR$
      we have
      \begin{enumerate}
        \item[i)]
          $$
            \sum_{i=0}^r \lambda_i v_i = 0
              \hspace{2ex} \Leftrightarrow \hspace{2ex}
              \lambda_0 = \ldots = \lambda_r
              \hspace{2ex} \Leftrightarrow \hspace{2ex}
              \lambda_i - \lambda_j = 0
              \text{ for all } i,j,
          $$
        \item[ii)]
          $$
            \sum_{i=0}^r \lambda_i v_i \in \Lambda
              \hspace{2ex} \Leftrightarrow \hspace{2ex}
              \lambda_i - \lambda_j \in \ZZ
              \text{ for all } i,j.
          $$
      \end{enumerate}
  \end{enumerate}
\end{remark}

  The \emph{moduli space of (abstract smooth) 
  $n$-marked rational tropical curves},
  denoted by $\calM_n$, is the fan in $\RR^{\binom{n}{2}} /
  \text{Im}(\Phi_n)$ that parametrizes metric trees 
  with positive lengths on the bounded edges 
  (and infinite lengths on the unbounded edges).
	The explicit construction of this space 
	can be found in	\cite{SS03}, \cite{Mi07} and \cite[section 3]{GKM07}.
  The cones of $\calM_n$
  are in one-to-one correspondence 
  with combinatorial types of $n$-marked trees (with
  $2$-valent vertices removed), and the dimension of a 
  cone equals the number of bounded edges in the respective
  combinatorial type.
  A general point in $\calM_n$ 
  (i.e. an element in the interior of a facet) 
  is a $3$-valent metric tree with $n-3$ bounded edges
  (hence $\dim(\calM_n) = n-3$). When all facets are equipped with
  weight $1$, $\calM_n$ fulfils the balancing
  condition. Hence $\calM_n$ is a tropical fan cycle.
	We denote the leaves by $x_1, \ldots, x_n$.
  If we work with $\calM_{n+1}$, the extra leaf is labelled by $x_0$.
  As $\calM_3$ is just a single point, we assume $n \geq 4$ in
  most cases.
	\\
  The notation $I|J$ denotes a non-trivial
  partition of $[n] = \{1, \ldots, n\}$ (or of $\{0\} \cup [n]$
  if we work with $\calM_{n+1}$) into the two disjoint
  subsets $I$ and $J$. In most cases --- the few
  exceptions will be mentioned ---
  we will consider this partition to be unordered.
  Occasionally, we use $I^c$ to denote the complement of $I$ and
  write $I|I^c$.
  If $|I| \neq 1 \neq |J|$, 
  such a partition describes a ray in $\calM_n$ generated by
  the metric tree $V_{I|J} \in \calM_n$ with only one bounded edge:
  \begin{center}
    \input{pics/IJ-Curves.pstex_t}
  \end{center}
  An edge of a tree is
  uniquely determined by the partition $I|J$ obtained when removing the edge.
  In this sense, we can regard the partitions $I|J$ as ``global'' labels
  of the edges of a tree, where $I|J$ labels the leaf $x_i$
  if $I = \{i\}$ or $J = \{i\}$, and a bounded edge otherwise. 
  A cone $\cell$ of $\calM_n$ is generated by the vectors $V_{I|J}$ 
  for all partitions which correspond to edges
  in the combinatorial type of $\cell$. In particular,
  it is natural to use the lengths of the bounded
  edges as local coordinates of a cone of $\calM_n$ ---
  this identifies each cone $\cell$ of $\calM_n$ with the positive
  orthant of $\RR^{\dim(\cell)}$.

  Let us make some remarks here. 
    We sometimes also think of $V_{I|J}$ as a vector in 
    $\RR^{\binom{n}{2}}$, in which case we also allow $|I| = 1$ 
    or $|J| = 1$ to get simpler formulas. However, as
    \[
		V_{\{k\}|[n] \setminus \{k\}} = 
    \Phi_n(0, \ldots, 0,1,0, \ldots, 0),
		\]
    these vectors vanish modulo $\text{Im}(\Phi_n)$. 
\\
    Note that for the following purposes, the
    underlying lattice of 
    $\RR^{\binom{n}{2}}/\Phi_n(\RR^n)$
    is \emph{not} $\ZZ^{\binom{n}{2}}/\Phi_n(\ZZ^n)$, 
    but is the lattice
    generated by the vectors $V_{I|J}$, denoted by $\Lambda_n$ (see
    \cite[3.3]{GKM07}). This is a technical issue,
    as it does not change the lattices of the
    cones $\Lambda_\cell, \cell \in \calM_n$, 
    but is necessary to make maps such
    as forgetful maps \emph{integer} affine.
\\
    As mentioned above, any metric tree can be realized uniquely by
		a smooth rational curve in the sense of definition \ref{DefAbstractCurves}
		(we actually prove this in proposition
    \ref{AbstractUniversalFamily}).
    Therefore $\calM_n$ really parametrizes what is promised by its
		name. 
\\
    Comparing $\calM_n$ to its classical counterpart, note that
		we will stick to the non-compact part of smooth curves 
		and will not use a compactification. However, the ``recursive structure''
		of the boundary known from the classical moduli space of stable curves
		is already visible in $\calM_n$ (without adding a ``boundary'').
		Namely, let $\cell$ be a cone of $\calM_n$ and let $\Gamma$ denote
    the corresponding combinatorial type of $n$-marked trees. Then
    it is easy to check that the star around $\cell$ satisfies
    \[
      \Star_{\calM_n}(\cell) = 
        \prod_{\substack{\vertex \text{ vertex} \\ 
                         \text{of } \Gamma}}
        \calM_{\val(\vertex)},
    \] 
		i.e.\ can be described as the product of ``smaller'' moduli spaces.

We will now define divisors respectively rational functions
that play the role of ``boundary'' divisors in our moduli space.
More precisely, if we actually would compactify $\calM_n$, these divisors
should be rationally equivalent to the actual boundary divisors.
All these divisors lie in the codimension one skeleton of $\calM_n$,
therefore represent higher-valent curves. 
As $\calM_n$ is simplicial, we can define a rational
function on $\calM_n$ by assigning an integer to each $I|J$: The 
integers are the values of the function at $V_{I|J}$ and on each 
cone we extend the function by linearity.

\begin{definition}
  We define the rational function $\varphi_{I|J}$ by
  $$
    \varphi_{I|J}(V_{I'|J'}) := \left\{ \begin{array}{ll}
                                  1      & \text{if } I=I' \text{ or } I=J', \\		
                                  0      & \text{otherwise}.
                                \end{array} \right.
  $$
  Furthermore, we use the notation
  $$
    \varphi_{k,l} := \varphi_{\{k,l\} | [n] \setminus \{k,l\}}
  $$
  for $k \neq l$.
\end{definition}

The ridges (codimension one cells) of $\calM_n$ correspond to 
combinatorial types of curves with one 4-valent vertex, which 
we will draw like this.
$$
  ^{A}_{D}\!\times^{B}_{C}
$$
Here $A$, $B$, $C$ and $D$ denote the four parts of the 
combinatorial type adjacent to the 4-valent vertex and by
abuse of notations also the sets of leaves belonging 
to this part (as, in most cases, this is the only information 
needed). 
If we want to compute the weight of a ridge 
$^{A}_{D}\!\times^{B}_{C}$ in the divisor of a rational function
on $\calM_n$,
we need to know how $\calM_n$ looks like locally around
$^{A}_{D}\!\times^{B}_{C}$. Obviously,
$\Star_{\calM_n}(^{A}_{D}\!\times^{B}_{C})$ contains three
facets corresponding to the three possibilities of "`resolving"' the $4$-valent
vertex by inserting a new bounded edge. 
\begin {center} 
  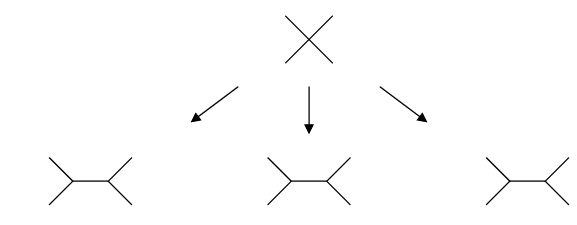 
\end {center}
The (representatives of the)
primitive vectors are $V_{A \cup B|C \cup D}$, $V_{A \cup C|B \cup D}$
and $V_{A \cup D|B \cup C}$. 
For the balancing condition 
around $^{A}_{D}\!\times^{B}_{C}$, it suffices to show the 
equation
$$
  V_{A \cup B|C \cup D} + V_{A \cup C|B \cup D} + V_{A \cup D|B \cup C} =
    V_{A|B \cup C \cup D} + V_{B|A \cup C \cup D} +
    V_{C|A \cup B \cup D} + V_{D|A \cup B \cup C},
$$
as all vectors on the right hand side
lie in the vector space spanned by the ridge 
$^{A}_{D}\!\times^{B}_{C}$, as required. But the equation 
follows from the fact that, on the level of metric trees,
the distance between two marked leaves is identical on both sides.
If both leaves belong to the same set $A,B,C,D$, the distance is
$0$, if not, it is $2$. 

Let us now compute the divisors of the functions $\varphi_{I|J}$.
In the following, a formula involving 
$I$, $J$ and $A$, $B$, $C$, $D$
stands for all permuted formulas as well, e.g.\ $I=A$ means
``$I=A$ or $I=B$ or $J=A$ $\ldots$''.

\begin{lemma} \label{LEMboundarydivisor}
  The \emph{boundary divisor} $\divisor(\varphi_{I|J})$ 
  carries the weight function
  \[
    \omega_{\varphi_{I|J}} (^{A}_{D}\!\times^{B}_{C}) =
      \left\{ \begin{array}{ll}
         1  & \text{if } I = A \cup B, \\		
         -1 & \text{if } I = A, \\
         0  & \text{otherwise}.
      \end{array} \right.
  \]
\end{lemma}

\begin{proof}
  Following from the previous discussion,
  the weight of $^{A}_{D}\!\times^{B}_{C}$ in 
  $\divisor(\varphi_{I|J})$ is by definition
  \begin{eqnarray*}
    \omega_{\varphi_{I|J}} (^{A}_{D}\!\times^{B}_{C}) 
      & = &
      \varphi_{I|J}(V_{A \cup B|C \cup D}) + 
      \varphi_{I|J}(V_{A \cup C|B \cup D}) + 
      \varphi_{I|J}(V_{A \cup D|B \cup C}) \\
    & &
      - \varphi_{I|J}(V_{A|B \cup C \cup D}) 
      - \varphi_{I|J}(V_{B|A \cup C \cup D})
      - \varphi_{I|J}(V_{C|A \cup B \cup D})
      - \varphi_{I|J}(V_{D|A \cup B \cup C}).
  \end{eqnarray*}
  Hence, this weight is $1$ if $I$ is the union
  of two of the sets $A,B,C,D$ and is $-1$ if
  $I$ equals one of the four sets. Otherwise,
  it is $0$. 
\end{proof}

  These divisors were computed before by Matthias Herold 
  (see \cite{H}).

\begin{remark}	
	In terms of the general toric geometry rules,
	the functions $\varphi_{I|J}$ respectively the
  divisors $\divisor(\varphi_{I|J})$ are
  the tropical analogues of the irreducible components
  of the boundary of the classical moduli space of stable curves
  $\Mbar_{0,n}$. 
	Instead of using this fact explicitly, 
	in the following we will show in purely tropical terms
	that the tropical divisors show the same
	intersection-theoretic behaviour as their
	classical counterparts. (We will need this
	anyway when dealing with parametrized curves later on.)
\end{remark}

\begin{lemma}
\label{ZeroIntersectionOf2Phis}
  The equation
  $$
    \varphi_{i,j} \cdot \varphi_{i,k} \cdot \calM_n = 0
  $$
  holds for $n \geq 4$ and pairwise different $i,j,k \in [n]$.
\end{lemma}

\begin{proof}
  An abstract curve $C$ cannot simultaneously have
  bounded edges with partitions 
  $\{i,j\}|\{i,j\}^c$ and $\{i,k\}|\{i,k\}^c$ (as for example
  the first partition forces $i$ and $k$ to be adjacent to 
  the same $3$-valent vertex).
  Let $C$ be a curve in $|\varphi_{i,k}|$. At least after
  resolving a $4$-valent vertex, it contains an edge
  with partition $\{i,k\}|\{i,k\}^c$ and can therefore
  \emph{not} contain an edge with partition
  $\{i,j\}|\{i,j\}^c$.
  But $\varphi_{i,j}$ just measures
  the length of such an edge if present. 
  Thus, 
  $\varphi_{i,j}|_{|\varphi_{i,k}|} \equiv 0$.
\end{proof}

Analogues of Psi-classes on tropical $\calM_n$ have been 
defined by G. Mikhalkin (\cite{Mi07}). Their intersections were studied in
in \cite{KM07}. We use the notion
Psi-divisor instead of Psi-class to emphasize that, in contrast to
the algebro-geometric case, tropically Psi-divisors are \emph{not} defined
up to rational equivalence. (Again, in the toric geometry language,
the tropical Psi-divisors are just the Minkowski weights associated
to the classical ones.)
In order to perform intersections, we need to describe our Psi-divisors by
rational functions.
Let us recall the important
definitions and results of \cite{KM07} here.

\begin{definition}
  We define the \emph{$k$-th Psi-function $\psi_k$} by
  $$
    \psi_k(V_{I|J}) := \frac{|I|(|I|-1)}{(n-1)(n-2)}
  $$
  for all partitions $I|J$ with $|I|, |J| \geq 2$ and $k \in J$.
\end{definition}

\begin{remark}
  Our function $\psi_k$ equals the function 
  $\frac{1}{\binom{n-1}{2}} f_k$ defined in \cite{KM07} (follows from 
  \cite[Lemma 2.6]{KM07}). In particular, $\psi_k$ is a convex
  function (cf. \cite[Remark 2.5]{KM07}). 
  Note that in this paper, $\psi_k$ and
  $\varphi_{I|J}$ denote functions and \emph{not} their corresponding
  divisors. On the other hand, as mentioned in subsection
  \ref{completeintersection}, this is only a matter of notation.
  For intersection-theoretic purposes, the actual choice of a function
  defining the same divisor does not matter.
\end{remark}

\begin{remark}
  Obviously the numbers $\psi_k(V_{I|J})$ 
  are only rational. A generalization of
  intersection theory to rational numbers is straightforward, 
  but also essentially  unnecessary: The weights
  of the divisor of $\psi_k$ turn out to be integers 
  (see the following proposition) and there exist integer rational 
  functions producing the same divisor (see \ref{PsiAsBoundary}). This
  particular function $\psi_k$ was chosen in \cite{KM07} because 
  of its symmetry.	
\end{remark}

\begin{proposition}[see \cite{KM07} 3.5]
  The divisor $\divisor(\psi_k)$ consists of the cones
  corresponding to trees where the 
  marked leaf $k$ is at a $4$-valent vertex, i.e.
  the weight of a facet in
  $\divisor(\psi_k)$ (which is a ridge in $\calM_n$)
  is 
  $$
  \omega_{\psi_k} (^{A}_{D}\!\times^{B}_{C}) =
    \left\{ \begin{array}{ll}
       1  & \text{if } \{k\} = A, \\		
       0  & \text{otherwise}
    \end{array} \right.
  $$
\end{proposition}

\begin{remark}
  As mentioned above, in toric geometry language and using
	the considering the embedding of the classical moduli space 
	of stable curves in the toric variety associated to $\calM_n$,
	it is easy to see that $\divisor(\psi_k)$ indeed represents
	the Minkowski weight associated to the classical $k$-th Psi-class.
	To check this, consider the one-dimensional boundary stratum $S$
	in the classical moduli space corresponding to reducible curves
	with dual graph $^{A}_{D}\!\times^{B}_{C}$. Each of such 
	curves contains exactly one $\PP^1$-component with $4$ special points,
	whereas all other components carry exactly $3$ special points and therefore
	are rigid. Hence $S$ is isomorphic to the moduli space of
	$4$-marked stable curves $\Mbar_{0,4} \cong \PP^1$. Let $\mathcal{L}_k$ be the line
	bundle whose fibre over a point corresponding to a curve $C$ is the cotangent
	space $T^*_{x_k} C$ at the marked point $x_k$. By definition the classical Psi-class
	is just the first Chern class of this line bundle. To compute the associated
	Minkowski weight, we should evaluate this Chern class on the one-dimensional
	boundary stratum described above. If $A = \{k\}$, this means that $x_k$ is
	one of the four special points on the non-rigid component. It follows
	that $\mathcal{L}_k|_S$ is equal to the corresponding Psi line bundle on $\Mbar_{0,4}$
	and we can compute $\deg(\mathcal{L}_k|_S) =1$. (For example, representing
	$\Mbar_{0,4}$ as the pencil of conics through $4$ points in the plane, there is 
	exactly one conic with prescribed tangent line at one of the $4$ points.)
	If $A \neq \{k\}$, the marked point $x_k$ lies on one of the rigid components.
	Therefore $\mathcal{L}_k|_S$ is the trivial line bundle and $\deg(\mathcal{L}_k|_S) =0$.
	This reproduces the weights from our proposition. 
	Again, let us emphasize that we do not really use this derivation as our arguments
	are purely tropical.
\end{remark}

\begin{notation}
  As in the conventional case we will introduce the following
  $\tau$-notation that makes formulas shorter and hides 
  ``unimportant'' data such as the number of marked leaves.
  For any positive integers $a_1, \ldots, a_n$ we define
  $$
    (\tau_{a_1} \cdot \ldots \cdot \tau_{a_n}) :=
      \psi_1^{a_1} \cdot \ldots \cdot \psi_n^{a_n} \cdot \calM_n.
  $$
  Every factor $\tau_{a_k}$ stands for a marked leaf and the index $a_k$
  serves as the exponent with which the corresponding Psi-function appears
  in the intersection product. If $\sum a_k = \dim(\calM_n) = n - 3$, 
  the above cycle is zero-dimensional (in fact, its only point 
  corresponds to the curve without bounded edges where all leaves
  are adjacent to one single vertex)
  and we define
  $$
    \langle\tau_{a_1} \cdot \ldots \cdot \tau_{a_n}\rangle :=
      \deg\big(\psi_1^{a_1} \cdot \ldots \cdot \psi_n^{a_n} 
      \cdot \calM_n\big).
  $$
\end{notation}

The main theorem of \cite{KM07} computes these intersection products
of Psi-divisors.

\begin{theorem}[Intersections of Psi-divisors for abstract curves, see \cite{KM07} 4.1]
                  \label{PsiProductsForAbstractCurves}
  The intersection product $(\tau_{a_1} \cdot \ldots \cdot \tau_{a_n})$ is the 
  subfan of $\calM_n$ consisting of the closure of the cones of
  dimension $n-3-\sum_{i=1}^n a_i$ whose interior curves $C$ 
  have the following property. \\
	Let $k_1, \ldots, k_q \subseteq N$ be the marked 
	leaves adjacent to a vertex $V$ of $C$. Then the valence of
  $V$ is
  $$
    \val(V) = a_{k_1} + \ldots + a_{k_q} + 3.
  $$
  Let us define the multiplicity of this vertex to be
  $\mult(V) := \binom{\val(V)-3}{a_{k_1}, \ldots, a_{k_q}}$.
  Then the weight of such a cone $\sigma$ in $X$ is
  $$
    \omega_X(\sigma) = \prod_V \mult(V),
  $$
  where the product runs through all vertices $V$ of an
  interior curve of $\sigma$.
\end{theorem}

  In this section we reprove the 
  zero-dimensional case of this theorem (see \ref{AbstractInvariants}).
  To do this, we first have to analyse how Psi- and boundary divisors
  intersect and how they behave when pulled back or pushed forward
  along forgetful morphisms.

\begin{lemma}
\label{ZeroIntersectionOfPhiAndPsi}
  It holds
  $$
    \varphi_{i,j} \cdot \psi_i \cdot \calM_n = 0
  $$
  for $n \geq 4$ and $i \neq j \in [n]$.
\end{lemma}

\begin{proof}
  Curves in $|\psi_i|$ cannot contain a bounded edge with partition
  $\{i,j\}|\{i,j\}^c$, as the leaf $i$ does not lie at a $3$-valent vertex. 
  Thus $\varphi_{i,j}$ vanishes on $|\psi_i|$.
\end{proof}

The forgetful map $\calM_{n+1} \rightarrow \calM_n$ 
that forgets the extra leaf $x_0$ is denoted
by $\ft_0$ (cf. \cite[4.1]{GM05} and \cite[3.8]{GKM07}). By \cite[3.9]{GKM07} 
this map is a tropical morphism. 
Therefore we can ask 
how Psi-functions behave when pulled 
back along $\ft_0$.

\begin{lemma}[Pull-back of Psi-functions
] 
                                 \label{PullbackOfPsi}
  Let $n \geq 4$ and let $\ft_0 : \calM_{n+1} \rightarrow \calM_n$ be the morphism 
  that forgets the leaf $x_0$. For $k \in [n]$ it holds
  $$
    \divisor(\psi_k) = \divisor(\ft_0^*\psi_k) + \divisor(\varphi_{0,k}).
  $$
\end{lemma}

\begin{proof}
  This can be proven by explicitly computing the weights of the codimension
  one faces of the three divisors. We distinguish four cases (up to renaming
  $A$, $B$, $C$ and $D$):
  $$
    \begin{array}{c|ccc}
      \omega_f(^{A}_{D}\!\times^{B}_{C}) & f = \psi_k 
        & f = \ft_0^*\psi_k & f = \varphi_{0,k} \\ \hline
      A=\{0,k\}                & 0 & 1 & -1 \\
      A=\{0\}, B=\{k\}         & 1 & 0 & 1 \\
      A=\{0, \ldots\}, B=\{k\} & 1 & 1 & 0 \\
      \text{otherwise}         & 0 & 0 & 0
    \end{array}
  $$
\end{proof}

\begin{corollary}
\label{SelfIntersectionOfPhi}
  Let $n \geq 4$ and let $\ft_0 : \calM_{n+1} \rightarrow \calM_n$ be the morphism 
  that forgets the leaf $x_0$.
  Then for $k \in [n]$ the following formulas hold.
  \begin{enumerate}
    \item
  	  $
  		  \varphi_{0,k}^2 = 
  			  - \ft_0^*(\psi_k) \cdot \varphi_{0,k}
      $ \vspace{1ex}
    \item
      $
        \psi_k^a = \ft_0^*(\psi_k)^a + 
          \ft_0^*(\psi_k)^{a-1} \cdot \varphi_{0,k}
      $ \vspace{1ex}
    \item
      $
        \psi_k^a = \ft_0^*(\psi_k)^a +
          (-1)^{a-1} \varphi_{0,k}^a 
      $
  \end{enumerate}
\end{corollary}

\begin{proof}
  All the formulas follow directly from 
  \ref{ZeroIntersectionOfPhiAndPsi} and \ref{PullbackOfPsi}.
\end{proof}

\begin{lemma} \label{PushforwardOfPhiAndPsi}
  Let $n \geq 4$ and let $\ft_0 : \calM_{n+1} \rightarrow \calM_n$ be the morphism 
  that forgets the leaf $x_0$ and choose $k \in [n]$. Then
  $$
    \ft_{0*}(\divisor(\varphi_{0,k})) = \ft_{0*}(\divisor(\psi_k)) = \calM_n.
  $$
\end{lemma}

\begin{proof}
  We show $\ft_{0*}(\divisor(\varphi_{0,k})) = \calM_n$ by direct
  computation. Let $\sigma'$ be a facet of
  $\calM_n$ corresponding to a $3$-valent combinatorial type.
  Let $V$ be the vertex adjacent to $k$. Then there exists
  precisely one cone $\sigma$ in $\divisor(\varphi_{0,k})$ whose
  image under $\ft_0$ is $\sigma'$, namely the cone obtained by
  attaching the additional leaf $x_0$ to the vertex $V$. Moreover,
  on such a cone, the length of the bounded edges remain unchanged
  under $\ft_0$ and therefore $\ft_0(\Lambda_\sigma) = \Lambda_{\sigma'}$.
  On the other hand,
  cones in $\divisor(\varphi_{0,k})$ with negative weight are not
  mapped injectively, as in this case $x_0$ is 
  adjacent to a $3$-valent vertex
  and stabilization is needed. This shows that 
  $\ft_{0*}(\divisor(\varphi_{0,k})) = \calM_n$. \\
  The equation $\ft_{0*}(\divisor(\psi_k)) = \calM_n$ follows from the 
  same argument or by using \ref{PullbackOfPsi}, the projection formula
  and $\ft_{0*}(\calM_{n+1}) = 0$ (because the dimension is too big).
\end{proof}

It is well-known that for the classical moduli space
$\Mbar_{0,n}$, the forgetful morphism plays the 
role of the universal family (cf.\ \cite[section 1.3]{KV}).
In the tropical setting we can 
prove the following statement

\begin{proposition}[Family property of $\ft_0$ for abstract curves] 
                                 \label{AbstractUniversalFamily}
  Let $p$ be a point in $\calM_n$ and let
  $C_p = \ft_0^{-1}(p)$ be the fibre of $p$
  under the forgetful morphism $\ft_0 : \calM_{n+1} \rightarrow \calM_n$.
  Then the following holds.
  \begin{enumerate}
    \item
      $C_p$ has the canonical structure of a one-dimensional polyhedral complex.
    \item
      The leaves of $C_p$ (as graph itself) are the
      facets where $x_0$ and another leaf $x_i$
      lie at the same $3$-valent vertex (i.e.
      the leaves are given by $L_i := \overline{\{y \in C_p |
      \varphi_{0,i}(y) > 0\}}$). Moreover
      $p \in \calM_n$ represents the $n$-marked 
      metric graph $(C_p,L_1, \ldots, L_n)$.
    \item
      When we equip all its facets with weight $1$, $C_p$ is
      a smooth abstract curve (in the sense
      of \ref{DefAbstractCurves}).
    \item
      Let $\sum_k \mu_k p_k = \varphi_1 \cdot \ldots \cdot
      \varphi_{n-3} \cdot \calM_n$ be a zero-dimensional
      cycle in $\calM_n$ obtained as the intersection
      product of convex functions $\varphi_j$. Then
      $$
        \ft_0^*(\varphi_1) \cdot \ldots \cdot
          \ft_0^*(\varphi_{n-3}) \cdot \calM_{n+1}
          = \sum_k \mu_k C_{p_k}.
      $$
      We write this as $\ft_0^*(\sum_k \mu_k p_k) = 
      \sum_k \mu_k C_{p_k}$.
  \end{enumerate}
\end{proposition}

\begin{proof}
  (a): As polyhedral complex, $C_p$ consists of the polyhedra
  $\ft_0|_\sigma^{-1}(p)$ for each cone $\sigma$ of $\calM_{n+1}$.
  The dimension of these polyhedra can be at most one as
  $\dim(f_0(\sigma)) \geq \dim(\sigma) - 1$ (it depends
  on whether $x_0$ is adjacent to a $3$-valent or higher-valent
  vertex). \\
  (b): Let $\Gamma_p$ denote the $n$-marked 
  metric graph represented by $p$. The bijective map
  $\Gamma_p \rightarrow C_p$
  indicated in the picture identifies the two graphs.
  \begin{center}
    \input{pics/GraphMapsToFiber.pstex_t}
  \end{center}
  (c): Let $V$ be a vertex of $C_p$. It corresponds
  to the metric graph $\Gamma_p$ with the extra leaf $x_0$
  adjacent to one of the vertices. Let us label the other
  edges adjacent to this vertex by $1, \ldots, m$ and let
  us divide the other leaves $[n] = I_1 \dcup \ldots
  \dcup I_m$ according to via which edge one reaches
  $x_i$ from $x_0$. 
  There are $m$ facets in $C_p$ adjacent to $V$ corresponding
  to moving $x_0$ on one of the edges.
  Hereby on has to shorten the edge $I_k|I_k^c$ as much
  as the length of $I_k\cup \{x_0\}|(I_k\cup \{x_0\})^c$
  increases.
  \begin{center}
    \input{pics/Vertex.pstex_t}
  \end{center}
  Thus the primitive integer vector of the corresponding
  facet with respect to $V$ is given by
  $$
    V_k := V_{I_k\cup \{x_0\}}
      - V_{I_k}.
  $$
  Note that this formula as well as the following ones also
  holds in the case that $I_k$ consists only of a single leaf $x_i$ 
  (which means $x_i$ is adjacent to the same vertex as $x_0$),
  as $V_{\{x_i\}} = 0 \in \RR^{\binom{n+1}{2}} / \text{Im}(\Phi_{n+1})$.
  To prove the statement we now use \ref{Smoothness}
  and verify the conditions i) and ii), which can be done by applying some 
  formulas of \cite{KM07}. Let $\mathcal{S}$ be the set of
  two-element subsets of $[n]$ (i.e. not containing $0$). It follows
  from \cite[2.3, 2.4, 2.6]{KM07} that the vectors
  $V_S, S \in \mathcal{S}$
  fulfil i) and ii) (with 
  $V = \RR^{\binom{n+1}{2}} / \text{Im}(\Phi_{n+1})$
  and $\Lambda = \Lambda_n$). Furthermore 
  \cite[2.6]{KM07} gives us a representation
  of our vectors in terms of the vectors $V_S$, namely
  $$
    V_{I_k} 
      = \sum_{\substack{S \in \mathcal{S} \\ S \subseteq I_k}} V_S
  $$
  $$
    V_{I_k\cup \{x_0\}} 
      = \sum_{\substack{S \in \mathcal{S} \\ S \cap I_k = \emptyset}} V_S,  
      = -\Big(\sum_{\substack{S \in \mathcal{S} \\ S \cap I_k \neq \emptyset}} V_S\Big),
  $$
  and therefore
  $$
    V_k = -\Big(\sum_{S \in \mathcal{S}} |S \cap I_k| \cdot V_S\Big).
  $$
  Now let $\lambda_1, \ldots, \lambda_m$ be arbitrary real coefficients. Then
  we obtain the formula
  $$
    \sum_{k = 1}^m \lambda_k V_k 
      = -\Big(\sum_{\substack{\{i,j\} \in \mathcal{S} \\
                              i \in I_k, j \in I_{k'}}}  
              (\lambda_k + \lambda_{k'}) \cdot V_{\{i,j\}}\Big).
  $$
  Now all differences of two coefficients on the left hand side 
  $\lambda_k - \lambda_k'$ can be obtained as differences of two coefficients
  on the right hand side (choose elements $i \in I_k, j \in I_{k'}, l \in
  I_{k''}$; then the coefficients of $V_{\{i,l\}}$ and $V_{\{j,l\}}$ differ
  by $\lambda_k + \lambda_{k''} - \lambda_{k'} - \lambda_{k''} = \lambda_k - \lambda_{k'}$).
  Conversely, a right hand side difference of coefficients equals the sum 
  of two left hand side differences. (The coefficients of $V_{\{i_1,i_2\}}$ and
  $V_{\{j_1,j_2\}}$ differ by
  $(\lambda_{k_1} - \lambda_{l_1}) + (\lambda_{k_2} - \lambda_{l_2})$,
  where $i_1 \in I_{k_1}, i_2 \in I_{k_2}, j_1 \in I_{l_1}, 
  j_2 \in I_{l_2}$.)
  Hence, as conditions \ref{Smoothness} i) and ii) hold for the vectors $V_S$, 
  they also hold for the vectors $V_k$. \\
  (d): First of all, the set-theoretic equation
  $$
    |\ft_0^*(\varphi_1) \cdot \ldots \cdot
      \ft_0^*(\varphi_{n-3}) \cdot \calM_{n+1}|
      \subseteq 
      \ft_0^{-1}(|\varphi_1 \cdot \ldots \cdot \varphi_{n-3} \cdot \calM_n|)  
      = \bigcup_k |C_{p_k}|.
  $$ 
  follows from \ref{ImageOfConvexPullbacks}.
  But the sets $|C_{p_k}|$ are pairwise disjoint (as they are
  fibres of pairwise different points) and belong to irreducible
  cycles (as the curves $C_{p_k}$ are smooth abstract curves).
  Thus any one-dimensional cycle whose support lies in $\bigcup_i |C_{p_k}|$
  is actually a sum $\sum_k \lambda_k C_{p_k}, \lambda_k \in \ZZ$.
  So it remains to check that in our case these coefficients $\lambda_k$
  coincide with $\mu_k$. To do this, we choose an arbitrary leaf $x_i \neq x_0$
  and consider the function $\varphi_{0,i}$ on $C_{p_k}$. On the leaf $L_i$ of
  $C_{p_k}$, where $x_0$ and $x_i$ are adjacent to the same $3$-valent vertex,
  it measures the length of the third edge, elsewhere it is constantly zero.
  Thus $\varphi_{0,i} \cdot C_{p_k} = V_{p_k}$, where $V_{p_k}$ is the vertex
  of $C_{p_k}$ adjacent to $L_i$ (where $x_0$ and $x_i$ lie together at
  a higher-valent vertex). Thus we get
  $$
    \ft_{0*}\big(\varphi_{0,i} \cdot (\sum_k \lambda_k C_{p_k})\big) 
      = \ft_{0*}\big(\sum_k \lambda_k V_{p_k}\big)
      = \sum_k \lambda_k p_k.
  $$
  On the other hand we can use projection formula and 
  \ref{PushforwardOfPhiAndPsi} and compute
  $$
    \ft_{0*}\big(\varphi_{0,i} \cdot 
      \ft_0^*(\varphi_1) \cdot \ldots \cdot
      \ft_0^*(\varphi_{n-3}) \cdot \calM_{n+1}\big)
      = \varphi_1 \cdot \ldots \cdot \varphi_{n-3} \cdot 
      \ft_{0*}(\varphi_{0,i} \cdot \calM_{n+1})
      = \sum_k \mu_k p_k.
  $$
  Comparing the coefficients proves the statement.
\end{proof}

\begin{remark}
  Hence there is a one-to-one correspondence between
  curves according to the ``old'' definition
  (i.e. as metric graphs) and definition \ref{DefAbstractCurves}.
  In particular, $\calM_n$ parametrizes smooth abstract curves
  in our sense. 
\end{remark}

\begin{theorem}[String equation for abstract curves
] \label{StringEquation}
  For zero-dimensional intersection products of Psi-divisors the
  following holds.
  $$
    \langle \tau_0 \prod_{k=1}^n \tau_{a_k} \rangle_d =
      \sum_{i=1}^{n} \langle \tau_{a_i - 1} 
      {\textstyle \prod\limits_{k \neq i}} \tau_{a_k} \rangle_d
  $$
\end{theorem}

\begin{proof}
  The proof is identical to the algebro-geometric one.
  We have to compute degree of the intersection product
  $\prod_{k=1}^n \psi_k^{a_k} \cdot \calM_{n+1}$.
  First we replace each term $\psi_k^{a_k}$ ($k \neq 0$) by
  $\ft_0^*(\psi_k)^{a_k} + \ft_0^*(\psi_k)^{a_k-1} \cdot \varphi_{0,k}$
  using \ref{SelfIntersectionOfPhi} (b) and multiply the product out.
  As $\varphi_{0,k} \cdot \varphi_{0,k'} = 0$ for $k \neq k'$
  (see \ref{ZeroIntersectionOf2Phis}), we only get the
  following $n+1$ terms.
  $$
    \prod_{k=1}^n \ft_0^*(\psi_k)^{a_k} \cdot \calM_{n+1}
      + \sum_{i=1}^n \ft_0^*(\psi_i)^{a_i - 1}
        \cdot {\textstyle \prod\limits_{k \neq i}} \ft_0^*(\psi_k)^{a_k}
        \cdot \varphi_{0,i} \cdot \calM_{n+1}
  $$
  Now we push this cycle forward along $\ft_0$ and use projection formula.
  The first term vanishes for dimension reasons and, as $\varphi_{0,i}$
  pushes forward to $\calM_n$ by \ref{PushforwardOfPhiAndPsi}, the other
  terms provide the desired result.
\end{proof}

\begin{remark}             \label{AbstractInvariants}
  As in the classical case, the string equation suffices
  to compute all intersection numbers of Psi-divisors
  of abstract curves. 
  Namely, if $\sum a_i = n - 3$, the equation
  $$
    \langle \tau_{a_1} \cdot \ldots \cdot \tau_{a_n} \rangle = 
      \frac{(n-3)!}{a_1! \cdot \ldots \cdot a_n!}
  $$
  holds. This was proven in \cite[4.2]{KM07} using the paper's main
  theorem \cite[4.1]{KM07} (cited here in 
  \ref{PsiProductsForAbstractCurves}).
  Note, however, that in order to prove the string equation
  it was not necessary to use \cite[4.1]{KM07}. Another
  independent proof of the above equality is given in
  \cite[Proposition 7.4]{Ka2}.
\end{remark}

\begin{lemma} \label{PushforwardOfPhi2}
  Let $n > 4$ and let $\ft_0 : \calM_{n+1} \rightarrow \calM_n$ be the morphism 
  that forgets the last leaf. Then
  $$
    \ft_{0*}(\divisor(\varphi_{I|J})) = 
      \left\{ \begin{array}{ll}
        \calM_n & \text{if } I = \{0,k\} \text{ or } J = \{0,k\} 
                  \text{ for some } k \in [n], \\		
        0       & \text{otherwise}.
      \end{array} \right.
  $$
\end{lemma}

\begin{proof}
  The first part is shown in \ref{PushforwardOfPhiAndPsi}. So let us prove
  the second part. First, we choose $i\in I$ and $j\in J$, both different
  from $0$. Consider a facet $\sigma'$ in
  $\calM_n$ corresponding to a combinatorial type where $x_i$ and $x_j$
  are adjacent to the same $3$-valent vertex $V$. All ridges in 
  $\calM_{n+1}$ mapping onto $\sigma'$, are obtained by attaching $x_0$ to any of 
  the vertices. If not attached to $V$, the induced partition $A,B,C,D$ cannot
  separate $i$ and $j$. If attached to $V$, the induced partition is
  $\{0\},\{i\},\{j\},D$. It follows from $\{0,i\} \neq I$ and 
  $\{0,j\} \neq J$ that $D$ intersects both $I$ and $J$ and therefore
  none of these types is contained in $\divisor(\varphi_{I|J})$.
  Hence $\sigma'$ is not contained in the push-forward of $\divisor(\varphi_{I|J})$.
  But $\calM_n$ is irreducible, thus $\ft_{0*}(\divisor(\varphi_{I|J})) = 0$.
\end{proof}

\begin{lemma}
\label{PsiAsBoundary}
  For $n \geq 4$ 
  we define
  $$
    (x_1|x_2,x_3) := 
      \sum_{\substack{I|J \\ 1\in I;\, 2,3\in J}} \divisor(\varphi_{I|J}).
  $$
  Then 
  $$
    \divisor(\psi_1) = (x_1|x_2,x_3).
  $$
\end{lemma}
  
\begin{proof}
  We use induction on	the number of leaves $n$. For $n=4$, only
  the partition $\{1,4\}|\{2,3\}$ contributes to the sum. But
  $\divisor(\psi_1)$ as well as $\divisor(\varphi_{1,4|2,3})$ is just
  the single vertex in $\calM_4$ parametrizing the curve $^{1}_{4}\!\times^{2}_{3}$ with
  weight $1$. For the induction step, assume $n \geq 4$ and consider the morphism
  $\ft_0 : \calM_{n+1} \rightarrow \calM_n$
  that forgets the leaf $x_0$ and let $I'|J'$ be a partition
  of $[n]$. Then $\ft_0^*(\varphi_{I'|J'})$ measures the sum
  of the lengths of the edges separating $I'$ and $J'$ if present.
  Hence we obtain 
  $$
    \ft_0^*(\varphi_{I'|J'}) = 
      \varphi_{I' \cup \{0\}|J'} + \varphi_{I'|J' \cup \{0\}}.
  $$
  Using the induction hypothesis, we conclude that $\ft_0^*(\psi_1)$
  equals the sum on the right hand side except for the partition
  $\{0,1\}|\{0,1\}^c$. This missing summand is
  provided by \ref{PullbackOfPsi}.
\end{proof}

\begin{lemma}

\label{PushforwardOfPsi2}
  Let $n \geq 4$ and let $\ft_0 : \calM_{n+1} \rightarrow \calM_n$ be the morphism 
  that forgets the leaf $x_0$. Then
  $$
    \ft_{0*}(\divisor(\psi_0)) = (n-2)\calM_n.
  $$
\end{lemma}

\begin{proof}
	We express $\psi_0$ as $(x_0|x_1,x_2)$ by \ref{PsiAsBoundary}
  and use linearity of the push-forward. Lemma \ref{PushforwardOfPhi2} 
  says that
  we get \emph{one} $\calM_n$ for each $\varphi_{\{0,k\}|\{0,k\}^c}$
  and zero for each other $\varphi_{I|J}$. As $k$ runs through
  $\{3, \ldots, n\}$, the statement follows. 
\end{proof}

\begin{proposition}[Dilaton equation for abstract curves
] 
                                  \label{DilatonEquation}
  Let $\langle \prod_{k=1}^n \tau_{a_k} \rangle$ 
  be a zero-dimensional intersection product. Then
  $$
  	\langle \tau_1 \cdot \prod_{k=1}^n \tau_{a_k} \rangle =
  	  (n-2) \langle \prod_{k=1}^n \tau_{a_k} \rangle.
  $$
\end{proposition}

\begin{proof}
  The proof is identical to the algebro-geometric one,
  using \ref{SelfIntersectionOfPhi}, \ref{ZeroIntersectionOfPhiAndPsi}, 
  \ref{PushforwardOfPhiAndPsi},  
  \ref{PushforwardOfPsi2} and the projection formula. \\
  As degree is preserved, we push forward 
  $(\tau_1 \cdot \prod_{k=1}^n \tau_{a_k})$
  along the forgetful morphism $\ft_0$ forgetting the extra leaf $x_0$
  corresponding to the factor $\tau_1$. To see what happens, we use 
  \ref{SelfIntersectionOfPhi} (b) and replace each term $\psi_k^{a_k}$ 
  by $\ft_0^*(\psi_k)^{a_k} + \ft_0^*(\psi_k)^{a_k-1} \cdot \varphi_{0,k}$. 
  When we multiply the whole product out,
  all summands containing a factor $\varphi_{0,k}$
  vanish when multiplied with $\psi_0$ (see \ref{ZeroIntersectionOfPhiAndPsi}).
  It follows 
  $$
    \psi_0 \cdot \prod_{k=1}^n \psi_k^{a_k} = \psi_0 \cdot \prod_{k=1}^n \ft_0^*(\psi_k)^{a_k} 
  $$
  and the projection formula together with $\ft_{0*}(\divisor(\psi_0)) = (n-2)\calM_n$
  from \ref{PushforwardOfPsi2} gives the desired result. \\
\end{proof}

\section {Intersections on the space of parametrized curves} \label {parametrizedcurves}


In the previous section, we proved that the tropical "`boundary divisors"'
and Psi-divisors satisfy exactly the same intersection-theoretic formulas as
their classical counterparts. 
We did this using purely tropical arguments, but could have used instead
tropicalization methods and toric intersection theory (cf.\ \cite{Ka2}).
The situation changes as we move on to curves together with maps to $\RR^n$. 
In this case, as $\Mbar_{0,n}(\PP^n,d)$ does not admit a nice toric 
embedding and therefore the tropicalization of its intersection ring
is not yet well understood. 
(In fact, we will see that, at least with our definitions,
the intersection theories are not completely identical.) 
So at least from now on, we are somehow forced to take this purely tropical approach.
Let us start with the necessary definitions.

A \emph{(labelled) degree $\Delta$ in $\RR^r$} is a finite set of labels together
with a map $\Delta \rightarrow \ZZ^r\setminus\{0\}$ to
the set of non-zero integer vectors. Furthermore the images
of this map, denoted by $v(x_i), i \in \Delta$ as they will later
play the role of the directions of the leaves $x_i$, sum up to zero, i.e.
$\sum_{i \in \Delta} v(x_i) = 0$. The number of elements in $\Delta$ 
is denoted by $\#\Delta$ (to distinguish it
from the support of a cycle). As an example, we
define the \emph{projective degree $d$ (in dimension $r$)} to be
the set $[(r+1)d]$ with the map
\begin{eqnarray*}
  [(r+1)d] & \rightarrow & \ZZ^r \setminus \{0\}, \\
  1, \ldots, d    & \mapsto     & -e_0, \\
  d+1, \ldots, 2d & \mapsto     & -e_1, \\
  \vdots          &             & \vdots \\
  rd+1, \ldots, (r+1)d & \mapsto & -e_r, 
\end{eqnarray*}
where, as usual, $e_1, \ldots, e_r$ denote the standard basis vectors
and $e_0 := e_1 + \ldots + e_r$.

\begin{definition}     \label{DefParametrizedCurves}
  An \emph{$n$-marked (labelled) parametrized curve of
  degree $\Delta$ in $\RR^r$} is a tuple $(C, h)$, where
  $C$ is an $[n] \dcup \Delta$-marked smooth abstract
  curve and $h : C \rightarrow \RR^r$ is a tropical morphism
  such that for all leaves $x_i$ the ray $h(x_i)
  \subseteq \RR^r$ has direction $v(x_i)$. Here $v(x_i)$ is set to be 
  zero if $i \in [n]$
	(i.e.\ the marked leaves $x_i, i \in [n]$ are contracted to a point). 
	The genus of $(C,h)$ is defined to be the genus
  of $C$.
\end{definition}

\begin{remark}
  The leaves $x_i, i \in [n]$ are called \emph{marked leaves}, as
  they correspond to the marked points of stable maps classically.
  Marked leaves are contracted by $h$. 
  In contrast to that we call the leaves $x_i, i \in \Delta$ 
  \emph{non-contracted leaves}. Our curves are called ``labelled'' as
  also the non-contracted leaves are labelled.\\
  Two parametrized curves $(C,h)$ and $(C',h')$ are called isomorphic
  (and therefore identified in the following) if there exists an
  isomorphism $\Phi : C \rightarrow C'$ identifying the labels and
  satisfying 
  $h = h' \circ \Phi$. \\
  Let us compare our definition to \cite[definition 4.1]{GKM07}.
  Conditions (a) and (b) in that definition make sure that
  $h$ is a tropical morphism in our sense (at least locally;
  but again, considering the family property of $\ft_0, \ev_0$ 
  over $\M{n}{r}{\Delta}$
  we will see that a global integer affine map
  $h$ always exists). Condition (c) is also contained in our definition.
\end{remark}

Let $\M{n}{r}{\Delta}$ be the moduli space of rational $n$-marked labelled
parametrized curves of degree $\Delta$ in $\RR^r$. Its construction as a tropical cycle
can be found in \cite[4.7]{GKM07}. After fixing one of the marked leaves $x_i$
as \emph{anchor leaf} (we avoid ``root leaf'' as, from the botanic point of view,
this does not make much sense), we can identify $\M{n}{r}{\Delta}$ with 
$\calM_{[n] \cup \Delta} \times \RR^r$, where the first factor
parametrizes the abstract curve $C$ and the second
factor contains the coordinates of the image point of the anchor leaf
$x_i$. As our curves are rational, this suffices to fix the morphism $h$.
Indeed, as $\Delta$ determines the directions of all leaves, we can use the 
balancing condition to recursively compute the directions of all bounded
edges as well. Hence $h$ is uniquely determined by the lengths of the
edges and the coordinates of one image point (in our case $h(x_i)$).

So again, cones in $\M{n}{r}{\Delta}$ correspond to combinatorial types
of the underlying abstract curves, but this time the minimal cone is not
zero- but $r$-dimensional because we can move the curve in $\RR^r$. \\
For enumerative purposes,
we would like to identify curves whose only difference is the labelling
of the non-contracted leaves. Let $\calM_{n}(\RR^r,\Delta)$
denote the set of these \emph{unlabelled} curves. Then the number of elements
in a general
fibre of the map $\M{n}{r}{\Delta} \rightarrow 
\calM_{n}(\RR^r,\Delta)$ forgetting the labelling of the
non-contracted leaves equals the number of possibilities 
to label a general unlabelled curve, which is
$$
  \Delta! := \prod_{v \in \ZZ^r \setminus \{0\}} n(v)!,
$$
where $n(v)$ denotes the number of times $v$ occurs as $v(x_i), i \in \Delta$.
Therefore each enumerative invariant
computed on $\M{n}{r}{\Delta}$ must simply be divided by $\Delta!$ 
to get the corresponding one in 
$\calM_{n}(\RR^r,\Delta)$. \\
From now on, 
$I|J$ denotes a (non-empty) partition of $[n] \dcup \Delta$ (or 
$\{0\} \dcup [n] \dcup \Delta$ if we work with $\M{n+1}{r}{\Delta}$). 
Again such partitions can be used as global labels of the edges
of our curves. The direction of the image of the corresponding
edge under $h$ is given by 
\label{directionvector}
$$
  v_{I|J} := \sum_{i \in I} v(x_i) = -(\sum_{j \in J} v(x_j))
$$ 
(as an exception, the ordering of $I$ and $J$ plays a little role here,
namely $v_{I|J} = - v_{J|I}$). We call $I|J$ 
\emph{reducible} if $v_{I|J} = 0$ (i.e. if the
corresponding edge is contracted). This is equivalent to requiring
that the corresponding split sets $\Delta_I = I \cap \Delta$ and 
$\Delta_J = J \cap \Delta$ fulfil the balancing condition, i.e.
are degrees on its own. Also the marked leaves split up into
$[n] = \{i \in I | v(x_i) = 0 \} \dcup \{j \in I | v(x_j) = 0 \}$.
In this sense, the partition corresponds (nearly) to a conventional partition 
$(S', \beta'| S'', \beta'')$ of the marked points $S = S' \dcup S''$
and the degree $\beta = \beta' + \beta''$, occurring for example in the 
classical splitting lemma.
However, note 
that in the tropical setting it is possible to permute 
non-contracted leaves with the same direction vector between $I$
and $J$ without changing the corresponding conventional partition, hence
in general several tropical reducible partitions correspond
to the same conventional partition. The \emph{non-reducible} partitions $I|J$ 
(i.e.\ $v_{I|J} \neq 0$) do not correspond to such a partition.

There exists a forgetful map 
$\ft' : \M{n}{r}{\Delta} \rightarrow \calM_{[n] \cup \Delta}$
forgetting just the position
of a curve in $\RR^r$. 
This forgetful map $\ft' : \M{n}{r}{\Delta} 
\rightarrow \calM_{[n] \cup \Delta}$ is a morphism of
tropical varieties, as after choosing a anchor leaf and identifying $\M{n}{r}{\Delta}$ with
$\calM_{[n] \cup \Delta} \times \RR^r$, $\ft'$ is just the projection
onto the first factor.
We use this to define Psi-functions on $\M{n}{r}{\Delta}$.

\begin{definition}[Psi-functions for parametrized curves]
  For a partition $I|J$ of $[n] \cup \Delta$ we define the
  function $\varphi_{I|J}$ on $\M{n}{r}{\Delta}$ to be 
  $\ft'^*(\varphi_{I|J}^{\text{abstr}})$, where
  $\varphi_{I|J}^{\text{abstr}}$ is the corresponding
  function on $\calM_{[n] \cup \Delta}$. \\
  For $i=1,\ldots,n$ we define \emph{the $k$-th Psi-function on 
  $\M{n}{r}{\Delta}$}
  to be $\psi_k := \ft'^*(\psi_k^{\text{abstr}})$,
  where the $\psi_k^{\text{abstr}}$ is the $k$-th Psi-function 
  on $\calM_{[n] \cup \Delta}$.
\end{definition}

\begin{remark}
  Again, in spite of defining functions we are actually interested in its
  divisors. Note that by \ref{PullBackUnderProjection} the pull-backs
  of the respective divisors do not depend on the particular functions. 
  Note that in \cite[Definition 2.2]{MR08} we used the notation
  $\psi_k$ for the divisor instead of the rational function.
\end{remark}

We can immediately generalize statement \ref{PsiProductsForAbstractCurves}
to parametrized curves (cf. \cite[Lemma 2.4]{MR08}).

\begin{proposition}[Intersections of Psi-divisors for parametrized curves] \label{PsiProductsForParametrizedCurves}
  Let $a_1, \ldots, a_n$ be positive integers and let
  $X = \prod_{k = 1}^n \psi_k^{a_k} \cdot \M{n}{r}{\Delta}$ be
  a product of Psi-divisors. Then $X$ is the subfan 
  of $\M{n}{r}{\Delta}$ consisting of the closure of the cones of
  dimension $n+\#\Delta-3-\sum_{i=1}^n a_i$ whose interior curves $C$ 
  have the following property.\\
	Let $k_1, \ldots, k_q \subseteq [n]$ be the marked 
	leaves adjacent to a vertex $V$ of $C$. Then the valence of
  $V$ is
  $$
    \val(V) = a_{k_1} + \ldots + a_{k_q} + 3.
  $$
  Let us define the multiplicity of this vertex to be
  $\mult(V) := \binom{\val(V)-3}{a_{k_1}, \ldots, a_{k_q}}$.
  Then the weight of such a cone $\sigma$ in $X$ is
  $$
    \omega_X(\sigma) = \prod_V \mult(V),
  $$
  where the product runs through all vertices $V$ of an
  interior curve of $\sigma$.
\end{proposition}

\begin{proof}
  Choose an anchor leaf and identify $\M{n}{r}{\Delta}$ with
$\calM_{[n] \cup \Delta}
  \times \RR^r$. Then $\ft'$ is just the projection on the first
  factor and we can apply \cite[9.6]{AR07}, i.e. instead of 
  intersecting the pull-backs of the $f_k$ on the product, we can just
  intersect the $f_k$ on the first factor and then multiply with
$\RR^2$.
  Thus,
  $$
    X = (\prod_{k = 1}^n (\psi_k^{\text{abstr}})^{a_k} 
    \cdot \calM_{[n] \cup \Delta}) \times \RR^r,
  $$
  where here $\psi_k^{\text{abstr}}$ denotes a Psi-function on
$\calM_{[n] \cup \Delta}$.
  Now, as the weight of $\RR^r$ is one and the combinatorics
  of a curve do not change under $\ft'$, the statements follows
  from \ref{PsiProductsForAbstractCurves}.
\end{proof}

\begin{proposition} \label{ParametrizedVersion}
  Let $\ft_0$ be the map $\M{n+1}{r}{\Delta} \rightarrow 
  \M{n}{r}{\Delta}$ that forgets the extra leaf $x_0$
  and assume $n + \#\Delta \geq 4$ (and $n \geq 1$).
  Furthermore, let $x_i, x_j, x_k$ be pairwise different leaves.
  Then the following equations hold (where
  all the occurring intersection products are computed
  in $\M{n}{r}{\Delta}$ or $\M{n+1}{r}{\Delta}$ respectively):
  \begin{enumerate}
    \item   \label{ZeroIntersectionOf2PhisParametrized}
      $\varphi_{i,j} \cdot \varphi_{i,k} = 0$ \vspace{1ex}
    \item   \label{ZeroIntersectionOfPhiAndPsiParametrized}
      $\varphi_{i,j} \cdot \psi_i = 0$ \vspace{1ex}
    \item   \label{PullbackOfPsiParametrized}
      $\divisor(\psi_k) = \divisor(\ft_0^*\psi_k) + \divisor(\varphi_{0,k})$ \vspace{1ex}
    \item   \label{SelfIntersectionOfPhiParametrizedA}
  	  $\varphi_{0,k}^2 = - \ft_0^*(\psi_k) \cdot \varphi_{0,k}$ \vspace{1ex}
    \item   \label{SelfIntersectionOfPhiParametrizedB}
      $\psi_k^a = \ft_0^*(\psi_k)^a + \ft_0^*(\psi_k)^{a-1} \cdot \varphi_{0,k}$ \vspace{1ex}
    \item   \label{SelfIntersectionOfPhiParametrizedC}
      $\psi_k^a = \ft_0^*(\psi_k)^a + (-1)^{a-1} \varphi_{0,k}^a$ \vspace{1ex}
    \item   \label{PushforwardOfPhiAndPsiParametrized}
      $\ft_{0*}(\divisor(\varphi_{0,k})) = \ft_{0*}(\divisor(\psi_k)) = \M{n}{r}{\Delta}$ \vspace{1ex}
    \item   \label{PushforwardOfPhi2Parametrized}
      $\ft_{0*}(\divisor(\varphi_{I|J})) = 
          \left\{ \begin{array}{ll}
            \M{n}{r}{\Delta} & \text{if } I = \{0,k\} \text{ or } J =\{0,k\} \text{ for some } k \in [n], \\
            0       & \text{otherwise}.
          \end{array} \right.
      $ \vspace{1ex}
    \item   \label{PsiAsBoundaryParametrized}
      $\divisor(\psi_i) = (x_i|x_j,x_k) := 
			  \sum_{\substack{I|J \\ i\in I;\, j,k\in J}} \divisor(\varphi_{I|J})$, \\
      where the sum runs also through \emph{non-reducible}
      partitions. \vspace{1ex}
    \item   \label{PushforwardOfPsi2Parametrized}
      $
        \ft_{0*}(\divisor(\psi_0)) = (n+\#\Delta-2)\M{n}{r}{\Delta},
      $\\
      (which is \emph{different} to the algebro-geometric factor $n-2$
      that equals the abstract case).
  \end{enumerate}
\end{proposition}

\begin{proof}
  As in the proof of \ref{PsiProductsForParametrizedCurves},
  we apply \cite[9.6]{AR07} to the morphism
  $\ft' : \M{n}{r}{\Delta} = \calM_{[n] \cup \Delta} \times \RR^r
  \rightarrow \calM_{[n] \cup \Delta}$ forgetting the position
  in $\RR^r$. This means that instead of computing
  the intersection product on $\M{n}{r}{\Delta}$ we can
  compute them on $\calM_{[n] \cup \Delta}$ and therefore
  use the corresponding statements for abstract curves.
  For statements \ref{PullbackOfPsiParametrized} -- 
  \ref{PushforwardOfPhi2Parametrized} and
  \ref{PushforwardOfPsi2Parametrized}
  we also use $\ft_0 = \ft_0^{\text{abstr}} \times \id$.
\end{proof}

\begin{definition}[Evaluation maps and their pull-backs]
\label{EvaluationMaps}
  The \emph{evaluation map} 
  $\ev_k : \M{n}{r}{\Delta} \rightarrow \RR^r$,
  for $k \in [n]$,
  maps each parametrized curve $(C, h)$ to the position of
  its $k$-th leaf $h(x_k)$ (see \cite[4.2]{GKM07}). 
  If we choose one of the marked leaves, say $x_a$, as anchor leaf, 
  then the evaluation maps are morphisms from 
  $\calM_{[n] \cup \Delta} \times \RR^r$ to $\RR^r$ obeying
  the following mapping rule.
  $$
    (C^\text{abstr},P) \mapsto
      P + \sum_{\substack{I|J \\ a \in I, k \in J}}
          \varphi_{I|J}(C^\text{abstr}) \, v_{I|J}
  $$ 
  In particular, if our
  anchor leaf is chosen to be $x_k$, then $\ev_k$ is just the
projection onto the 
  second factor.
  Let $C \in \Zci_m(\RR^r)$ be given by $C = h_1 \cdot
  \ldots \cdot h_l \cdot X$. Then we can apply
  \ref{PullBackUnderProjection}
  which states that there is a well-defined
  \emph{pull-back of $C$ along $\ev_k$}
  $$
    \ev_k^*(C) :=
      \ev_k^*(h_1) \cdot \ldots \cdot \ev_k^*(h_l).
  $$
\end{definition}

\begin{proposition}[Family property of $\ft_0,\ev_0$ for parametrized curves]
               \label{ParametrizedUniversalFamily}
  Let $p$ be a point in $\M{n}{r}{\Delta}$ and let
  $C_p = \ft_0^{-1}(p)$ be the fibre of $p$
  under the forgetful morphism $\ft_0 : \M{n+1}{r}{\Delta} 
  \rightarrow \M{n}{r}{\Delta}$.
  Then the following holds.
  \begin{enumerate}
    \item
      When we equip all its facets with weight $1$, $C_p$ is
      a rational smooth abstract curve. Its leaves are naturally 
      $[n] \dcup \Delta$-marked by 
      $L_i := \overline{\{y \in C_p |
      \varphi_{0,i}(y) > 0\}}$.
    \item
      The tuple $(C_p, \ev_0|_{|C_p|})$ is an $n$-marked
      parametrized curve of degree $\Delta$. Moreover,
      $p$ represents $(C_p, \ev_0|_{|C_p|})$.
    \item
      Let $\sum_k \mu_k p_k = \varphi_1 \cdot \ldots \cdot
      \varphi_{n + \#\Delta -3} \cdot \M{n}{r}{\Delta}$ be a zero-dimensional
      cycle in $\M{n}{r}{\Delta}$ obtained as the intersection
      product of convex functions $\varphi_j$. Then
      $$
        \ft_0^*(\varphi_1) \cdot \ldots \cdot
          \ft_0^*(\varphi_{n + \#\Delta -3}) \cdot \M{n+1}{r}{\Delta}
          = \sum_k \mu_k C_{p_k}.
      $$
      We write this as $\ft_0^*(\sum_k \mu_k p_k) = 
      \sum_k \mu_k C_{p_k}$.
  \end{enumerate}
\end{proposition}

\begin{proof}
  (a): First of all, let us fix an anchor leaf $x_a, a \in [n]$ in
order
  to identify $\M{n+1}{r}{\Delta} = \calM_{n+\#\Delta+1} \times \RR^r$
and 
  $\M{n}{r}{\Delta} = \calM_{[n] \cup \Delta} \times \RR^r$. 
  We use again $\ft_0 = \ft_0^{\text{abstr}}
  \times \id$, where $\ft_0^{\text{abstr}}$ is the corresponding
forgetful map
  on the abstract spaces. Then the fibre of $p = (p', P)$ equals
  $C_{p'} \times \{P\}$, where $C_{p'}$ is the $[n] \dcup
\Delta$-marked
  rational smooth abstract curve considered in
\ref{AbstractUniversalFamily}
  (a)--(c). \\
  (b): We have to check that the direction of the rays
  $\ev_0(L_i)$ are correct. For curves in $L_i$, the only length
  that varies is that of the third edge adjacent to the 
  same $3$-valent vertex as $x_i$ and $x_0$. Hence we can use the
description
  of $\ev_0$ in \ref{EvaluationMaps} and obtain for all $y \in L_i$
  $$
    \ev_0|_{L_i}(y) = Q + \varphi_{0,i}(y) \cdot
v_{\{0,i\}|\{0,i\}^c},
  $$
  where $Q \in \RR^r$ is some constant vector. But
$v_{\{0,i\}|\{0,i\}^c} =
  v(x_i) + v(x_0) = v(x_i)$ is the expected direction. \\
  To show that $p = (p', P)$ represents $(C_{p}, \ev_0|_{|C_{p}|})$ it
  actually suffices to prove that the anchor leaf $L_a$ of $C_p$ is
mapped
  to the point $P$ under $\ev_0$, which is obviously the case as
  $\ev_0|_{L_a} = \ev_a|_{L_a}$ and $\ev_a$ is just the projection
  on the second factor of $C_{p'} \times \{P\}$.  \\
  (c): We can use literally the same proof as in the abstract case 
  \ref{AbstractUniversalFamily} (d) using \ref{ParametrizedVersion}
  \ref{PushforwardOfPhiAndPsiParametrized}.
\end{proof}

\begin{notation}[Tropical Gromov-Witten invariants]
  Let us now extend our $\tau$-notation to the 
  case of parametrized curves.
  For any positive integers $a_1, \ldots, a_n$ and 
  complete intersection cycles 
  $C_1, \ldots, C_n$ $ \in \Zci_*(\RR^r)$ we define
  $$
    (\tau_{a_1}(C_1) \cdot \ldots \cdot \tau_{a_n}(C_n))_\Delta^{\RR^r} :=
      \psi_1^{a_1} \cdot \ev_1^*(C_1)
      \cdot \ldots \cdot 
      \psi_n^{a_n} \cdot \ev_n^*(C_n) \cdot \M{n}{r}{\Delta}.
  $$
  Once again, each factor $\tau_{a_k}(C_k)$ stands for a marked leaf
  subject
  to $a_k$ Psi-conditions and to the condition that it must meet
  $C_k$. Let
  $c_k$ be the codimension of $C_k$ in $\RR^r$. 
  If $\sum (a_k + c_k) = \dim(\M{n}{r}{\Delta})$ $ = n + \#\Delta + r - 3$, 
  the above cycle is zero-dimensional 
  and we denote its degree by
  $$
    \langle\tau_{a_1}(C_1) \cdot \ldots \cdot \tau_{a_n}(C_n)\rangle_\Delta^{\RR^r} .
  $$
  These numbers are called \emph{tropical descendant 
  Gromov-Witten invariants}. 
	In \cite{MR08} these numbers were studied in the case
	$r = 2, \Delta = d$ and all $C_i$ are
  tropical lines.
\end{notation}

\begin{remark}[Enumerative meaning of tropical Gromov-Witten invariants]
  \label{enumerativerelevance}
  Let $(\tau_{a_1}(C_1) \cdot \ldots \cdot \tau_{a_n}(C_n))$ be
  an intersection product as defined above. If we set 
  $X = \prod_{k = 1}^n \psi_k^{a_k} \cdot \M{n}{r}{\Delta}$ 
  and apply \ref{generalposition2} to the morphisms $\ev_k : X \rightarrow \RR^r$,
  we can conclude the following (as discussed in \ref{ComparisonOfPullBackAndPreimage}):
  After replacing all the cycles $C_k$ by general translations (called
  \emph{general conditions} in the following),
  $Z := \tau_{a_1}(C_1) \cdot \ldots \cdot \tau_{a_n}(C_n))$ 
  is the set
  of curves $C$ such that 
  \begin{itemize}
    \item
      every vertex $V \in C$ with adjacent marked leaves 
      $k_1, \ldots, k_q$ fulfils 
      $$
        \val(V) \geq a_{k_1} + \ldots + a_{k_q} + 3,
      $$
    \item
      for all $k = 1, \ldots, n$ it holds
      $$
        \ev_k(C) \in C_k.
      $$
  \end{itemize}
  Additionally, the facets of $Z$
  (i.e. general curves) are equipped with (possibly zero) weights. \\
  Moreover, assume that all the cycles $C_k$ can be described by convex
  functions $h_1 \cdots h_l \cdot \RR^r$. Then by \ref{convex}, all these
  weights are positive (in particular, $|Z|$ really \emph{is}
  the set of such curves). \\
  Thus, if $Z$ is zero-dimensional, 
  $\deg(Z) = \langle\tau_{a_1}(C_1) \cdot \ldots \cdot \tau_{a_n}(C_n)\rangle$
  is the number of curves satisfying the above properties, counted with 
  a certain integer multiplicity/weight. Now again, if all $C_k$ can be described by
  convex functions, all these multiplicities and in particular
  $\langle\tau_{a_1}(C_1) \cdot \ldots \cdot \tau_{a_n}(C_n)\rangle$
  are positive. 
\end{remark}

  Let $\ft_0 : \M{n+1}{r}{\Delta} \rightarrow 
  \M{n}{r}{\Delta}$ be the morphism 
  that forgets the leaf $x_0$. Then by abuse of
  notation the equation
  $$
    \ft_0^*(\ev_k) = \ev_k
  $$
  holds for all $k \in [n]$.
	This equality directly implies the following extensions
	of the string and dilaton equation to the
	case of parametrized curves.

\begin{theorem}[String equation for parametrized curves
] 
                                \label{StringEquationForMaps}
  Let $(\tau_0(\RR^r) \cdot \prod_{k=1}^n \tau_{a_k}(C_k))_\Delta$
  be a zero-dimensional cycle.
  Then 
  $$
    \langle\tau_0(\RR^r) \cdot \prod_{k=1}^n
\tau_{a_k}(C_k)\rangle_\Delta =
  	  \sum_{k=1}^n \langle\tau_{a_k-1}(C_k) \cdot \prod_{l \neq k}
\tau_{a_l}(C_l)\rangle_\Delta.
  $$  
\end{theorem}

\begin{theorem}[Dilaton equation for parametrized curves
] 
                                \label{DilatonEquationForMaps}
	Let $(\tau_1(\RR^r) \cdot \prod_{k=1}^n \tau_{a_k}(C_k))_\Delta$
  be a zero-dimensional cycle.
  The following equation holds.
  $$
     \langle\tau_1(\RR^r) \cdot \prod_{k=1}^n
\tau_{a_k}(C_k)\rangle_\Delta =
  	  (n+\#\Delta-2) \langle \prod_{k=1}^n
\tau_{a_k}(C_k)\rangle_\Delta.
  $$
\end{theorem}

\begin{proof}[Proofs]
  In both cases, the proofs are completely analogous to the
  abstract case using \ref{ParametrizedVersion} and
	$\ft_0^*(\ev_k) = \ev_k$. 
\end{proof}

\begin{remark}
  Note that the factor appearing in the dilaton equation is
  different from the algebro-geometric one, due to 
  $\ft_{0*}(\psi_0) = (n+\#\Delta-2) \cdot 
  \M{n}{r}{\Delta}$ (cf. \ref{ParametrizedVersion} 
  \ref{PushforwardOfPsi2Parametrized}).
\end{remark}

\begin{lemma}
\label{CommutingEvkEvlInFrontOfPhikl}
  Let $h$ be a rational function. Then
  $$
    \ev_k^*(h) \cdot \varphi_{k,l} \cdot \M{n}{r}{\Delta} =
      \ev_l^*(h) \cdot \varphi_{k,l} \cdot \M{n}{r}{\Delta}
  $$
\end{lemma}

\begin{proof}
  In all curves corresponding to points in $\divisor(\varphi_{k,l})$,
the
  leaves $k$ and $l$ lie at a common vertex. Therefore their
coordinates in
  $\RR^r$ must agree, which means $\ev_k|_{|\divisor(\varphi_{k,l})|}
=
  \ev_l|_{|\divisor(\varphi_{k,l})|}$. The result follows.
\end{proof}

  For a given labelled degree $\Delta$,
  we define $\delta(\Delta)$ to be the associated
  unlabelled degree in the sense of subsection
  \ref{rationalequivalence}: $\delta(\Delta)$
  is the one-dimensional balanced fan
  in $\RR^r$
  consisting of all the rays generated by the 
  direction vectors $v_k, k \in \Delta$ appearing in $\Delta$.
  The weight of such a ray $\RR_\geq v$, where $v$ is primitive, is given by
  $$
    \sum_{\substack{k \in \Delta \\ v_k \in \ZZ_{>0} v}}
      |\ZZ v / \ZZ v_k|.
  $$
  Obviously, if $(C,h) \in \M{n}{r}{\Delta}$ is an arbitrary $n$-marked
  parametrized curve of degree $\Delta$, then by definition
  $\delta(h(C)) = \delta(\Delta)$ holds. \\
  For a given rational function $h$ on $\RR^r$ we define
  $h \cdot \Delta$ to be $\deg(h \cdot \delta(\Delta))$.

\begin{lemma}
\label{PushforwardOfPullbackOfDivisor}
  Let $h$ be a rational function on $\RR^r$
  and define $Y := \ev_0^*(h) \cdot \M{n+1}{r}{\Delta}$. Then
  $$
    \ft_{0*}(Y) = (h \cdot \Delta) \M{n}{r}{\Delta}.
  $$
\end{lemma}

\begin{proof}
  As our moduli space $\M{n}{r}{\Delta}$ is irreducible, 
  we know that $\ft_{0*}(Y) = \alpha \cdot \M{n}{r}{\Delta}$
  for an integer $\alpha$. 
  To compute this number, we set
  $m := n + \#\Delta + r - 3$ and consider the zero-dimensional 
  intersection product 
  $Z = \varphi_1 \cdots \varphi_m \cdot \M{n}{r}{\Delta}$
  of arbitrary convex functions $\varphi_1, \ldots, \varphi_m$
  such that $\deg(Z) \neq 0$
  (e.g. $Z = \psi_1^{m-r} \cdot \ev_1(P)$ for some point 
  $P \in \RR^r$). If we pull back $Z$ along $\ft_0$,
  we know by the projection formula
  $$
    \deg(\ev_0(h) \cdot \ft_0^*(Z)) = \alpha \cdot \deg(Z).
  $$ 
  On the other hand, by the family property of 
  $\ft_0$ we know that $Z$ is the union of the curves
  represented by the points in $Z$ (with according weights)
  and therefore the push-forward $\ev_{0*}(\ft_0^*(Z))$ is rationally
  equivalent to its degree
  $$
    \delta(\ev_{0*}(\ft_0^*(Z))) = \deg(Z) \cdot \delta(\Delta).
  $$
  So, applying the projection formula to $\ev_0$, we obtain
  $$
    \deg(\ev_0(h) \cdot \ft_0^*(Z)) = \deg(Z) \cdot (h \cdot \Delta).
  $$
  But this implies $h \cdot \Delta = \alpha$, which proves the claim.
\end{proof}

We can now prove the following rather general version of the divisor equation.

\begin{theorem}[Divisor equation
]
                              \label{DivisorEquation}
  Let $h$ be a rational function on $\RR^r$
  and let $(\prod_{k=1}^{n} \tau_{a_k}(C_k))_\Delta$ be a 
  one-dimensional cycle.
  Then 
  $$
    \langle \tau_0(h) \cdot \prod_{k=1}^{n} \tau_{a_k}(C_k) \rangle_\Delta
      = (h \cdot \Delta) \langle \prod_{k=1}^{n} \tau_{a_k}(C_k) \rangle_\Delta
      + \sum_{k=1}^{n} \langle \tau_{a_k-1}(h \cdot C_k) 
      \prod_{l \neq k} \tau_{a_l}(C_l)\rangle_\Delta. 
  $$
\end{theorem}

\begin{proof}
  First we use \ref{ParametrizedVersion} \ref{SelfIntersectionOfPhiParametrizedB}
  and \ref{ZeroIntersectionOf2PhisParametrized}: 
  We replace each factor $\psi_k^{a_k}$ by 
  $\ft_0^*(\psi_k)^{a_k} + \ft_0^*(\psi_k)^{a_k-1} \cdot \varphi_{0,k}$
  and multiply out. All terms containing two $\varphi$-factors vanish.
  In terms with only one factor $\varphi_{0,k}$, we replace 
  $\ev_0(h)$ by $\ev_k(h)$ using \ref{CommutingEvkEvlInFrontOfPhikl}.
  Now we push forward along $\ft_0$ and produce the desired equation 
  by applying \ref{PushforwardOfPullbackOfDivisor}
  and $\ft_{0*}(\divisor(\varphi_{0,k})) = \M{n}{r}{\Delta}$.
\end{proof}

Note that the divisor equation can be used to prove the statement
of \cite[Proposition 2.10]{MR08}.

\section{The Splitting Lemma} \label{splittingcurves}


The basic fact used to compute Gromov-Witten type invariants of 
$\overline{M}_{g,n}(X, \beta)$ is the recursive structure
of its boundary: Its irreducible components correspond to reducible curves with a 
certain partition of the combinatoric data and therefore are (nearly) a product
of two ``smaller'' moduli spaces. In this section we will investigate how far
this principle can be carried over to the tropical world.

\subsection{The case of abstract curves}


\begin{definition}
  Let $S$ be a finite set. 
  By $\calM_S$ we denote the moduli space of $|S|$-marked tropical 
  curves $\calM_{|S|}$ where we label the leaves by elements in $S$.
  For each partition $I|J$ of $[n]$ we construct the map
  $\rho_{I|J} : \calM_{I\cup\{x\}} \times \calM_{J\cup\{y\}}
  \rightarrow \varphi_{I|J} \cdot \calM_n$ by the following rule:
  Given two curves $(p_I, p_J) \in 
  \calM_{I\cup\{x\}} \times \calM_{J\cup\{y\}}$, we remove the extra leaves 
  $x$ and $y$ and glue the curves together at the two vertices to 
  which these leaves have been adjacent. In other words, we glue $x$ and 
  $y$ together by creating a bounded edge whose length we define to 
  be $0$. In the coordinates of the space of tree metrics, this map 
  is given by the linear map
  \begin{eqnarray*}
    \rho_{I|J} : \RR^{\binom{I}{2}} \times \RR^{\binom{J}{2}}
      & \rightarrow & \RR^{\binom{n}{2}}, \\
    (p_I, p_J) & \mapsto & p,
  \end{eqnarray*}
  where
  $$
    p_{k,l} :=       
      \left\{ \begin{array}{ll}
        {p_I}_{k,l}             & \text{if } k,l \in I, \\
	
        {p_J}_{k,l}            & \text{if } k,l \in J, \\
        {p_I}_{k,x} + {p_J}_{y,l} & \text{if } k \in I, l \in J. 
      \end{array} \right.
  $$
  Caution: This map does \emph{not} induce a linear map on the 
  corresponding quotients in which our moduli spaces are balanced
  and therefore $\rho_{I|J}$ is \emph{not} a tropical morphism of our
  moduli spaces. 
  Even more, $\rho_{I|J}$ is not even locally linear around ridges
  of our 
  moduli spaces considered as balanced complexes in the quotients. 
  On the other hand, $\rho_{I|J}$ is at least 
  piecewise linear (i.e. it is linear on all cones of 
  $\calM_{I\cup\{x\}} \times \calM_{J\cup\{y\}}$). Its image is a polyhedral
  complex, namely the positive part of 
  $\varphi_{I|J} \cdot \calM_n$ (i.e. it consists of all (faces of)
  facets $^{A}_{D}\!\times^{B}_{C}$ with $A \cup B = I$).
\end{definition}

\begin{definition}[Morphisms of rational polyhedral complexes]
    \label{MorphismPolyhedralComplexes}
  Let $X$ and $Y$ be (rational) polyhedral complexes. Then a \emph{
  morphism of polyhedral complexes} is a map $\rho : |X|
  \rightarrow |Y|$ that satisfies for each polyhedron $\sigma \in X$
  \begin{enumerate}
    \item
      $\rho(\sigma) \in Y$, 
    \item
      $\rho|_\sigma$ is affine linear,
    \item
      $\rho(\Lambda_\sigma) \subseteq \Lambda_{\rho(\sigma)}$.
  \end{enumerate}
  We call $\rho$ an \emph{isomorphism of polyhedral complexes}
  if there exists an inverse morphism. 
  It other words, an isomorphism is a bijection between $|X|$ and $|Y|$ (as well 
  as between $X$ and $Y$) and $\rho(\Lambda_\sigma) = 
  \Lambda_{\rho(\sigma)}$ for all $\sigma \in X$.
\end{definition}

\begin{lemma}[Intersections of Psi-functions with the boundary] \label{HowPhiIJTimesPsiLooksLike}
  The facets of the fan $\varphi_{I|J} \cdot \psi_1^{a_1} \cdot \ldots \cdot 
  \psi_n^{a_n} \cdot \calM_n$ with positive weight are precisely the
  cones $\sigma$ in $\calM_n$ with the following properties. \\
  Consider a curve in the interior of $\sigma$. 
  Let $E(V) \in [n]$ be the set of leaves adjacent to a vertex $V$ and 
  let $P(V)$ be the $val(V)$-fold partition of $[n]$ obtained by 
  removing $V$. Then the following holds.
  \begin{enumerate}
    \item
      There exists one special vertex $V_{\text{spec}}$ whose partition 
      $P(V_{\text{spec}})$ is 
      a subpartition of $I|J$ and whose valence is
      $(\sum_{k \in E(V)} a_k) + 4$.
    \item
      Let $m_I$ be the number of sets in $P(V_{\text{spec}})$
      contained in $I$. Then $m_I + 1 = 
      (\sum_{k \in E(V) \cap I} a_k) + 3$ (together with (a), the
      analogue $m_J + 1 = 
      (\sum_{k \in E(V) \cap J} a_k) + 3$ follows).
      In particular, $m_I, m_J > 1$.
    \item
      The valence of all other vertices $V$ equals
      $(\sum_{k \in E(V)} a_k) + 3$.
  \end{enumerate}
  Furthermore, the facets of $\varphi_{I|J} \cdot \psi_1^{a_1} \cdot \ldots \cdot 
  \psi_n^{a_n} \cdot \calM_n$ with negative weight fulfil the 
  same properties (a) and (c) and the property
  \begin{enumerate}
    \item[(b')]
      Let $m_I$ (resp. $m_J$) be the number of sets in $P(V_{\text{spec}})$
      contained in $I$ (resp. $J$). Then $m_I = 1$ or
      $m_J = 1$, i.e. $I \in P(V_{\text{spec}})$ or $J \in
      P(V_{\text{spec}})$.
  \end{enumerate}
\end{lemma}

\begin{proof}
  We know how $X := \psi_1^{a_1} \cdot \ldots \cdot 
  \psi_n^{a_n} \cdot \calM_n$ looks like by \ref{PsiProductsForAbstractCurves}.
  In the combinatorial type
  of a facet of $X$ the valence of each vertex is $(\sum_{k \in E(V)} a_k) + 3$;
  in the combinatorial type of a ridge, there is one special vertex $V_{\text{spec}}$
  with valence $(\sum_{k \in E(V)} a_k) + 4$. 
  The balancing condition of a ridge is given by the equation
  $$
    \sum_{I'|J'} \omega_{I'|J'} V_{I'|J'}
      = \sum_{\substack{I'|J' \\ I' \in P(V_{\text{spec}})}}
                  \lambda_{I'|J'} V_{I'|J'},
  $$
  where the left hand sum runs through all superpartitions $I'|J'$ of $P(V_{\text{spec}})$ not appearing
  in the right hand sum, $\omega_{I'|J'}$
  denotes the weight of the facet obtained by inserting an edge $I'|J'$ and $\lambda_{I'|J'}$ is
  some (rational) coefficient. Therefore the weight $\omega$ that this ridge obtains when intersecting
  $X$ with $\varphi_{I|J}$ is given by
  $$
    \omega =
    \begin{cases}
      0             & \text{if $I|J$ is \emph{not} a superpartition of $P(V_{\text{spec}})$,} \\
      \lambda_{I|J} & \text{if $I \in P(V_{\text{spec}})$ or $J \in P(V_{\text{spec}})$,} \\
      \omega_{I|J}  & \text{otherwise.}
    \end{cases}
  $$
  This already shows two implications. As all weights
  $\omega_{I'|J'}$ are at least non-negative, a ridge can only obtain a negative weight if
  it fulfils conditions (a), (b') and (c). On the other hand,
  if a ridge of $X$ satisfies properties (a), (b) and (c), then $\omega_{I|J}$
  and hence the ridge obtains a positive weight. It remains to show the converse,
  which can be done by proving that all $\lambda_{I'|J'}$
  are non-negative.
  To see this, we consider the balancing equation in $\RR^{\binom{r}{2}}$ and compare some coordinate
  entries. \\
  Let $K$ be an arbitrary element of $P(V_{\text{spec}})$; we want to show that 
  $\lambda_K := \lambda_{K|K^c}$ is non-negative.
  We choose two more arbitrary elements $L_1, L_2$ in 
  $P(V_{\text{spec}})$ and fix some leaves $k \in K$,
  $l_i \in L_i$. Now the $k,l_i$-entry of the right hand side equals $\lambda_K + \lambda_{L_i}$
  and analogously the $l_1,l_2$-entry equals $\lambda_{L_1} + \lambda_{L_2}$. Therefore, by 
  adding the two $k,l_i$-entries and subtracting the $l_1,l_2$-entry we get $2\lambda_K$. Meanwhile,
  on the left hand side we get
  \begin{eqnarray*}
    2\lambda_K 
      & = & 
      \sum_{\substack{I'|J' \\ k \in I' \\ l_1 \in J'}} \omega_{I'|J'} + 
      \sum_{\substack{I'|J' \\ k \in I' \\ l_2 \in J'}} \omega_{I'|J'} - 
      \sum_{\substack{I'|J' \\ l_1 \in I' \\ l_2 \in J'}} \omega_{I'|J'}  \\
    & = &
      \sum_{I'|J'} \alpha_{I'|J'} \omega_{I'|J'},
  \end{eqnarray*}
  where
  $$
    \alpha_{I'|J'} = \left\{ \begin{array}{ll}
                               2 & \text{if } k \in I', \; l_1,l_2 \in J' \\
                               0 & \text{if } k, l_1 \in I',\; l_2 \in J' \\
                               0 & \text{if } k, l_2 \in I',\; l_1 \in J' \\
                               0 & \text{if } k, l_1,l_2 \in I'.
                             \end{array} \right.
  $$
  But as all the weights $\omega_{I'|J'}$ are non-negative, it follows
  that $\lambda_K$ is non-negative.
\end{proof}

\begin{proposition} \label{IsomPositivePartCurves}
  The map
  $$
    \rho_{I|J} : \big(\prod_{k \in I} \psi_k^{a_k} \cdot \calM_{I\cup\{x\}}\big)
      \times \big(\prod_{k \in J} \psi_k^{a_k} \cdot \calM_{J\cup\{y\}}\big)
      \rightarrow
      (\varphi_{I|J} \cdot \psi_1^{a_1} \cdot \ldots \cdot 
      \psi_n^{a_n} \cdot \calM_n)^+
  $$
  is a well-defined isomorphism of polyhedral complexes.
\end{proposition}

\begin{proof}
  We have to check the conditions of \ref{MorphismPolyhedralComplexes}.
  Using the lengths of the bounded edges as local coordinates on the cones,
  this follows directly from the description of the target complex in \ref{HowPhiIJTimesPsiLooksLike}. 
	The inverse map is given by splitting a
  given curve at its special vertex $V_{\text{spec}}$.
\end{proof}

\subsection{The case of parametrized curves}


\begin{definition}
  Let $I|J$ be a reducible partition and let $\Delta_I, \Delta_J$
  be the corresponding splitting of the tropical degree $\Delta$.
  Let $Z = \max(x_1,y_1) \cdot \ldots \cdot \max(x_r,y_r) \cdot 
  \RR^r \times \RR^r$ denote the diagonal in $\RR^r \times \RR^r$
  and consider the map
  $$
    \ev_x \times \ev_y :
      \M{I\cup\{x\}}{r}{\Delta_I}
      \times \M{J\cup\{y\}}{r}{\Delta_J}
      \rightarrow \RR^r \times \RR^r.
  $$
  We define
  $$
    Z_{I|J} := (\ev_x \times \ev_y)^*(Z)
  $$
  We furthermore define $\pi_{I|J} : Z_{I|J} \rightarrow 
  \M{n}{r}{\Delta}$ by
  \begin{eqnarray*}
    \calM_{I\cup\{x\}} \times \RR^r
      \times \calM_{J\cup\{y\}} \times \RR^r
      & \rightarrow
      & \calM_{[n] \cup \Delta} \times \RR^r \\
    \big((p_I, P), (p_J, Q)\big)
      & \mapsto
      & (\rho(p_I, p_J), P),
  \end{eqnarray*}
  where we choose the same anchor leaf for 
  $\M{I\cup\{x\}}{r}{\Delta_I}$ and 
  $\M{n}{r}{\Delta}$ and $\rho$ is the gluing map for 
	abstract curves from the previous subsection.
\end{definition}

\begin{proposition} \label{IsomPositivePartMaps}
  The map
  $$
    \pi_{I|J} : \psi_1^{a_1} \cdot \ldots \cdot 
      \psi_n^{a_n} \cdot Z_{I|J}
      \rightarrow
      (\varphi_{I|J} \cdot \psi_1^{a_1} \cdot \ldots \cdot 
      \psi_n^{a_n} \cdot \M{n}{r}{\Delta})^+
  $$
  is a well-defined isomorphism of polyhedral complexes.
\end{proposition}

\begin{proof}
  This follows from \ref{IsomPositivePartCurves} and from 
  $\ev_x|_{Z_{I|J}} = \ev_y|_{Z_{I|J}}$ (which follows
  from both
  \ref{ImageOfConvexPullbacks} ($Z$ is described by convex
  functions) as well as from
  \ref{PullBackUnderProjection} 
  ($\ev_x \times \ev_y$ can be considered as a projection)).
\end{proof}

\begin{remark}
  Restricting to curves from $Z_{I|J}$ makes sure that
	the positions of the marked leaves are preserved 
  under $\pi_{I|J}$, i.e.\ (by abuse of notation)
  for all $i \in I$, but also $j \in J$
  we have$\ev_i \circ \pi_{I|J} = \ev_i$ resp.\
  $\ev_j \circ \pi_{I|J} = \ev_j$.
\end{remark}

\begin{lemma}
  Let $E = (\varphi_{I|J} \cdot \tau_{a_1}(C_1) \cdot \ldots \cdot 
  \tau_{a_n}(C_n))_\Delta$ be a 
  zero-dimensional cycle.   Then all points of $E$
  lie in $(\varphi_{I|J} \cdot \psi_1^{a_1} \cdot \ldots \cdot 
  \psi_n^{a_n} \cdot \M{n}{r}{\Delta})^+$. 
\end{lemma}

\begin{proof}
  By \ref{LocalIntersectionOfFunctions} we can compute the weight of a point $p \in E$
  locally around $p$
  in $X := \varphi_{I|J} \cdot \psi_1^{a_1} \cdot \ldots \cdot 
  \psi_n^{a_n} \cdot \M{n}{r}{\Delta}$, namely we can focus on
  $\Star_X(p)$.
  Assume $p \notin (\varphi_{I|J} \cdot \psi_1^{a_1} \cdot \ldots \cdot 
  \psi_n^{a_n} \cdot \M{n}{r}{\Delta})^+$. Then curves corresponding
  to points in $\Star_X(p)$ contain a bounded edge corresponding to
  the
  partition $I|J$ (see \ref{HowPhiIJTimesPsiLooksLike}). But as $I|J$ is chosen
  to be reducible, this edge is a contracted bounded edge whose
  length does not change the positions of the marked leaves in
  $\RR^r$.
  Therefore, if we denote by $\ev = \ev_1 \times \ldots \times
  \ev_n$ the product of all evaluation maps, the image of 
  $\Star_X(p)$ under $\ev$ has smaller dimension which implies
  $\ev_*(\Star_X(p)) = 0$.
  Hence, by projection formula, the weight of $p$ in $E$ must
  be zero.
\end{proof}

The following
statement combines 
\ref{generalposition2}, in particular item
(c), and the preceding result.

\begin{corollary} \label{GeneralEvConditions}
  Let $E = (\varphi_{I|J} \cdot \tau_{a_1}(C_1) \cdot \ldots \cdot 
  \tau_{a_n}(C_n))_\Delta$ be a 
  zero-dimensional cycle. If we substitute the cycles $C_i$
  by general translations, we can assume that all
  points of $E$ lie in the interior of a facet
  of $(\varphi_{I|J} \cdot \psi_1^{a_1}  
  \cdot \ldots \cdot 
  \psi_n^{a_n} \cdot \M{n}{r}{\Delta})^+$. This operation does
  not change
  the degree of $E$ by remark \ref{PullBackPreservesNumericalEquivalence}.
\end{corollary}

As provisional result of this discussion, we can formulate the following:

\begin{proposition} \label{BoundaryVsDiagonal}
  Let $E = (\varphi_{I|J} \cdot \tau_{a_1}(C_1) \cdot \ldots \cdot 
  \tau_{a_n}(C_n))_\Delta$ be a 
  zero-dimensional cycle. Then the equation
  $$
    \langle\varphi_{I|J} \cdot \tau_{a_1}(C_1) \cdot \ldots \cdot 
  \tau_{a_n}(C_n)\rangle_\Delta 
      =
      \langle\tau_{a_1}(C_1) \cdot \ldots \cdot 
      \tau_{a_n}(C_n) \cdot Z_{I|J}\rangle_{\Delta_I, \Delta_J}
  $$
  holds.
\end{proposition}

\begin{proof}
  We denote $X := \psi_1^{a_1} \cdot \ldots \cdot 
  \psi_n^{a_n} \cdot Z_{I|J}$ and
  $Y := \varphi_{I|J} \cdot \psi_1^{a_1} \cdot \ldots
  \cdot \psi_n^{a_n} \cdot \M{n}{r}{\Delta}$ and
  assume that the conditions $C_i$ are general.
  Then \ref{GeneralEvConditions} implies that,
  for each point $p \in E$, we have an
  isomorphism of cycles
  $\pi_{I|J} : \Star_X(\pi_{I|J}^{-1}(p)) 
  \rightarrow \Star_Y(p)$.
  By \ref{LocalIntersectionOfFunctions} this suffices to show that
  the weights of $p$ and $\pi_{I|J}^{-1}(p)$
  in their respective intersection products
  coincide.
\end{proof}

\subsection{Splitting the diagonal}


Up to now, we have seen that intersecting with a ``boundary''
function $\varphi_{I|J}$ leads to intersection products
in two smaller moduli spaces $\M{I\cup\{x\}}{r}{\Delta_I}$
and $\M{J\cup\{y\}}{r}{\Delta_J}$. However, the factor 
$(\ev_x \times \ev_y)^*(Z)$ still connects these
two smaller spaces. 
In order to obtain independent intersection products
on the smaller spaces, we have to split the diagonal 
contribution.
In the algebro-geometric case, this can be easily done as the \emph{class} of the
diagonal $Z$ in e.g. $\PP^r \times \PP^r$ can be written
as the sum of products of classes in the factors
$$
  [Z] = [L^0 \times L^r] + [L^1 \times L^{r-1}]
    + \ldots + [L^r \times L^0],
$$
where $L^i$ denotes an $i$-dimensional 
linear space in $\PP^r$.
But this can \emph{not} copied directly in our setting (see below). 
In some sense, for the first time we meet a disadvantage due to the non-compactness
of our spaces.
Our notion of rational equivalence is ``too strong'' for this 
application, as it is inspired by the idea that two rationally
equivalent objects should be rationally equivalent in \emph{any}
toric compactification. However, we will discuss here how far the conventional
plan can be carried out anyway.

The general plan is the following. Set 
$$
  X_I := (\tau_0(\RR^r) \cdot \prod_{k \in I} \tau_{a_k}(C_k))_{\Delta_I} 
    \hspace{5ex} \text{ in } \M{I\cup\{x\}}{r}{\Delta_I}
$$
and 
$$
  X_J := (\tau_0(\RR^r) \cdot \prod_{k \in J} \tau_{a_k}(C_k))_{\Delta_J} 
    \hspace{5ex} \text{ in } \M{J\cup\{y\}}{r}{\Delta_J}.
$$
We want to compute the degree of 
$$
  (\tau_{a_1}(C_1) \cdot \ldots \cdot 
    \tau_{a_n}(C_n) \cdot Z_{I|J})_{\Delta_I, \Delta_J}
    =(\ev_x \times \ev_y)^*(Z) \cdot (X_I \times X_J),
$$
or, by the projection formula,
$$
  \deg(Z \cdot (\ev_x(X_I) \times \ev_y(X_J))).
$$
Now we would like to replace the diagonal $Z$
by something like 
$$
  S := \sum_\alpha (M_\alpha \times N_\alpha),
$$
where $M_\alpha, N_\alpha$ are cycles in $\RR^r$ such that
$S$ intersects $\ev_x(X_I) \times \ev_y(X_J)$ like $Z$.
But note that $S$ cannot be rationally equivalent to $Z$
(in the sense of \cite{AR08}), as this would imply
that both cycles must have
the same recession fan, i.e. must have the same directions
towards infinity. 
To come out of this, we need more information about how the
push-forwards $\ev_x(X_I)$ and 
$\ev_y(X_J)$ look like; in particular, we would like to know
how their degrees/recession fans can look like. Let us formalize
this first.

Let $\Theta$ be a complete simplicial fan in $\RR^r$
and let $Z_k(\Theta)$ be the group of $k$-dimensional
cycles $X$ whose support lies in the $k$-dimensional
skeleton of $\Theta$, i.e. $|X| \subseteq |\Theta^{(k)}|$.
Fix a basis of $Z_*(\Theta) := \oplus_{k=0}^r Z_k(\Theta)$
denoted by $B_0, \ldots, B_m$ (where we may assume
$B_0 = \{0\}$ and $B_m = \RR^r$).
If the degree $\delta(X)$ of an arbitrary cycle is contained
in $Z_*(\Theta)$, we say $X$ is \emph{$\Theta$-directional}.
For such a cycle there exist integer coefficients $\lambda_e$
such that $X \sim \delta(X) = \sum_{e=1}^m \lambda_e B_e$. \\
For each ray $\rho \in \Theta^{(1)}$ with primitive
vector $v_\rho$ let $\varphi_\rho$ be the 
rational function on $\Theta$ uniquely defined by
$$
  \varphi_\rho(v_{\rho'}) =     \left\{ \begin{array}{ll}
                                  1      & \text{if } \rho'=\rho, \\		
                                  0      & \text{otherwise}.
                                \end{array} \right.
$$

\begin{lemma} \label{AlphaInvertible}
  The linear map
  \begin{eqnarray*}
    Z_*(\Theta)    & \rightarrow & \ZZ^{m+1}, \\
    X              & \mapsto     & (\deg(B_0 \cdot X), 
                            \ldots, \deg(B_m \cdot X)), 
  \end{eqnarray*}
  (where $\deg(.)$ is set to be zero if the dimension of
  the argument is non-zero) is injective.
\end{lemma}

\begin{proof}
  Let $X \in Z_k(\Theta)$ be an element of the kernel, which 
  implies that $\deg(X \cdot Y) = 0$ for all $Y \in Z_{r-k}(\Theta)$.
  Now, in fact the remaining is identical to the 
  proof of \cite[Lemma 6]{AR08}: 
  Assume that $X$ is non-zero and therefore there exists a cone 
  $\sigma \in \Theta^{(k)}$ such that $\omega_X(\sigma) \neq 0$.
  As $\Theta$ is simplicial, this cone is generated by 
  $k$ rays $\rho_1, \ldots, \rho_k$. 
  Let us consider $\varphi_{\rho_k} \cdot X$ and in particular
  the weight of $\tau := \langle \rho_1, \ldots, \rho_{k-1} \rangle$ in 
  this intersection product. As primitive vector $v_{\sigma/\tau}$
  we can use $\frac{1}{|\Lambda_\sigma / (\Lambda_\tau + \Lambda_{\rho_k})|}
  v_{\rho_k}$ (it might not be an integer vector, but modulo
  $V_\tau$, it is a primitive generator of $\sigma$). Analogously,
  we can get any primitive vector around $\tau$ as a multiple of
  an appropriate $v_\rho$. But as $\varphi_{\rho_k}$ is zero 
  on all of these vectors but $v_{\rho_k}$, we get
  $$
    \omega_{\varphi_{\rho_k} \cdot X}(\tau) = 
      \frac{\omega_X(\sigma)}{|\Lambda_\sigma / (\Lambda_\tau + \Lambda_{\rho_k})|}
      \neq 0.
  $$
  Now induction shows 
  $$
    \deg(\varphi_{\rho_1} \cdots \varphi_{\rho_k} \cdot X)
      = \omega_{\varphi_{\rho_1} \cdots \varphi_{\rho_k} \cdot X}(\{0\})
      = \frac{\omega_X(\sigma)}
             {|\Lambda_\sigma / (\Lambda_{\rho_1} + \ldots + \Lambda_{\rho_k})|}
      \neq 0.
  $$
  This means we have found a
  $\Theta$-directional cycle $Y :=
  \varphi_{\rho_1} \cdots \varphi_{\rho_k} \cdot \RR^r
  \in Z_{r-k}(\Theta)$ with $\deg(X \cdot Y) \neq 0$,
  which contradicts the assumption that $X$ is an element of
  the kernel.
\end{proof}

With respect to the basis $B_0, \ldots, B_m$, the map
defined in the previous lemma has the matrix representation $\alpha :=
(\deg(B_e \cdot B_f))_{ef}$.
Obviously $\alpha$ is a symmetric matrix.
The lemma implies that this matrix is invertible over $\QQ$,
and we denote the inverse by $(\beta_{ef})_{ef}$.
The coefficients of this matrix can be used to replace the 
diagonal $Z$ of $\RR^r \times \RR^r$ by a sum of products
of cycles in the two factors (namely $\sum_{e,f} \beta_{ef} (B_e \times B_f)$)
--- at least with respect
to $\Theta$-directional cycles.

\begin{lemma} \label{DiagonalVsBasisProducts}
  Let $X \sim \sum_e \lambda_e B_e, Y \sim \sum_f \mu_e B_e$ be two 
  $\Theta$-directional cycles in $\RR^r$ with complementary dimension.
  Then
  $$
    \deg(Z \cdot (X \times Y)) =
    \deg(X \cdot Y) =
      \sum_{e,f} \deg(X \cdot B_e) \beta_{ef} \deg(Y \cdot B_f).
  $$
\end{lemma}

\begin{proof}
  Denote $\lambda := (\lambda_1, \ldots, \lambda_m), 
  \mu := (\mu_1, \ldots, \mu_m)$. 
  We get 
  \begin{eqnarray*}
    \sum_{e,f} \deg(X \cdot B_e) \beta_{ef} \deg(Y \cdot B_f)
      & = & (\alpha \cdot \lambda)^T \cdot \beta \cdot (\alpha \cdot \mu) \\
      & = & \lambda^T \cdot \alpha^T \cdot \beta \cdot \alpha \cdot \mu \\
      & = & \lambda^T \cdot \alpha \cdot \beta \cdot \alpha \cdot \mu \\
      & = & \lambda^T \cdot \alpha \cdot \mu = \deg(X \cdot Y).
  \end{eqnarray*}
\end{proof}

Using this, our original goal of deriving a tropical splitting lemma
can be formulated as follows.

\begin{theorem}[Splitting Lemma
] \label{SplittingLemma}
  Let $E = (\varphi_{I|J} \cdot \prod_{k=1}^n \tau_{a_k}(C_k))_\Delta^{\RR^r}$ be a 
  zero-dimensional cycle, where $I|J$ is a reducible partition.
  Moreover, let us assume that $\Theta$ is a complete simplicial fan
  such that (with the notations from above) $\ev_x(X_I)$ and $\ev_y(X_J)$
  are $\Theta$-directional. Let $B_0, \ldots, B_m$ be a basis of $Z_*(\Theta)$
  and let $(\beta_{ef})_{ef}$ be the inverse matrix (over $\QQ$) of 
  $(\deg(B_e \cdot B_f))_{ef}$.
  Then the following equation holds.
  $$
    \langle\varphi_{I|J} \cdot {\textstyle \prod\limits_{k=1}^n} \tau_{a_k}(C_k)\rangle_\Delta 
      = 
      \sum_{e,f}
      \langle{\textstyle \prod\limits_{k \in I}} \tau_{a_k}(C_k) \cdot \tau_0(B_e)\rangle_{\Delta_I}
      \; \beta_{ef} \;
      \langle\tau_0(B_f) \cdot {\textstyle \prod\limits_{k \in J}} \tau_{a_k}(C_k)\rangle_{\Delta_J}
  $$
\end{theorem}

\begin{proof}
  Follows from the general plan above and \ref{BoundaryVsDiagonal}.
\end{proof}

\begin{remark}
  Of course, in toric geometry language, 
	the basis $B_0, \ldots, B_m$ 
  corresponds to a basis $\gamma_0, \ldots, \gamma_m$
  of the cohomology groups of $\bX(\Fan)$ (the toric variety 
	associated to $\Theta$).
  As the cup-product and the intersection
  product of cycles are equivalent (cf.\ theorem
  \ref{FanDisplacementEqualsTropInt}), the corresponding matrix
  $(\deg(\gamma_e \cup \gamma_f))_{ef}$ is equal to $\alpha$.
  This implies that the coefficients
  $\beta_{ef}$ appearing in the tropical splitting lemma really are the
  same as in the associated algebro-geometric version.
%
%
\end{remark}

\subsection{The directions of families of curves}


The above splitting lemma
is only useful if, at least for a certain class of invariants,
the fan of directions $\Theta$ is fixed and well-known. 
This is one of the main problems when transferring the algebro-geometric
theory to the tropical set-up. However, in this subsection we will show that
in some cases the problem can be solved.

\begin{remark} \label{DirectionsR=1}
  In the easiest case, namely if $r=1$, the situation is trivial.
  There is one unique complete simplicial fan 
  $\Theta = \{\RR_{\leq 0}, \{0\}, \RR_{\geq 0}\}$ and any subcycle
  is $\Theta$-directional. Also, with $B_0 = \{0\}, B_1 = \RR$, the statement
  of \ref{DiagonalVsBasisProducts} is obvious here.
\end{remark}

Let us now consider curves in the plane, i.e. $r=2$. Let 
$F = (\tau_0(\RR^2) \cdot \prod_{k = 1}^n \tau_{a_k}(C_k))_\Delta^{\RR^2}$
be a one-dimensional family of plane curves (with unrestricted
leaf $x_0$). We define $\Theta(F)$ to be the complete
fan in $\RR^2$ which contains the following rays:
all directions appearing in $\Delta$ and furthermore
all rays in $\delta(C_k)$ if $\dim(C_k) = 1$ and 
$a_k > 0$.

\begin{proposition}                  \label{DirectionsOfPushforwardAlongEv}
  Let $F = (\tau_0(\RR^2) \cdot \prod_{k = 1}^n \tau_{a_k}(C_k))_\Delta^{\RR^2}$
  be a one-dimensional family of plane curves (with unrestricted
  leaf $x_0$). Let us furthermore assume that $a_k \leq 1$ if
  $\dim(C_k) = 2$
  (i.e. if a leaf is not 
  restricted by $\ev$-conditions,
  only one Psi-condition is allowed).
  Then $\ev_{0*}(F)$ is $\Theta(F)$-directional.
\end{proposition}

\begin{proof}
  As before, we replace each factor $\psi_k^{a_k}$ by 
  $\ft_0^*(\psi_k)^{a_k} + \ft_0^*(\psi_k)^{a_k-1} \cdot \varphi_{0,k}$
  and multiply out. Consider the term without $\varphi$-factors. It is the 
  fibre of $(\prod_{k = 1}^n \tau_{a_k}(C_k))_\Delta$ (which is finite)
  under $\ft_0$ (see family property \ref{ParametrizedUniversalFamily}) 
  and moreover the push-forward
  of the fibre along $\ev_0$ is just the sum/union of the images in $\RR^r$
  of the parametrized curves corresponding to the points in 
  $(\prod_{k = 1}^n \tau_{a_k}(C_k))_\Delta$. But these curves have degree $\Delta$,
  thus by definition their images are $\Theta(F)$-directional. \\
  So let us consider the term with the factor $\varphi_{0,k}$. Here,
  $\ev_0$ and $\ev_k$ coincide (see \ref{CommutingEvkEvlInFrontOfPhikl}), 
  so we can in fact compute the push-forward along
  $\ev_k$. As $\ev_k = \ev_k \circ \ft_0$ (by abuse of notation), 
  we can first push-forward along $\ft_0$
  and get the term $(\tau_{a_k - 1}(C_k) \cdot \prod_{l \neq k} \tau_{a_l}(C_l))$. \\
  Now, if $\dim(C_k) = 2$, by our assumptions $a_k - 1 = 0$ -- in which case
  we can use induction to prove the statement -- or this term
  does not appear at all. \\
  On the other hand, if $\dim(C_k) = 0,1$, we can use the fact that the
  push-forward is certainly contained in $C_k$ -- therefore,
  $\dim(C_k) = 0$ is trivial and $\dim(C_k) = 1$ works
  as we added the directions of $C_k$ to $\Theta(F)$ if $a_k > 0$. \\
  This finishes the proof, as all terms with more $\varphi$-factors vanish.
\end{proof}

\begin{remark}
  A weaker version of this lemma
  can be obtained by assuming general conditions 
	and directly studying the behaviour of
	$\ev_0$ on an unbounded ray in $F$ (see \cite[Lemma 3.7]{MR08}).
\end{remark}

\begin{remark}
  Consider the family $F = (\tau_0(\RR^2) \tau_0(P) \tau_2(\RR^2))_1^{\RR^2} = \ev_1^*(P) \cdot \psi_2^2
  \cdot \M{3}{2}{1}$ of curves of projective degree $1$. It consists of the following types of curves.
  \begin{center}
    \input{pics/PushForwardCounterexample.pstex_t}
  \end{center}
  Its push-forward along $\ev_0$ also contains the inverted standard
  directions
  $(1,0)$, $(0,1)$ and $(-1,-1)$. Therefore this family is a counterexample to our statement
  if we drop the condition on the number of Psi-conditions allowed at leaves not restricted by
  incidence conditions.
\end{remark}

\begin{remark} \label{DirectionsForHigherR}
For higher dimensions ($r > 2$), only few cases are explored.
If we restrict to projective degree $d$ and banish all
Psi-conditions, i.e. for a family 
$F = (\tau_0(\RR^r) \cdot \prod_{k = 1}^n \tau_{0}(C_k))_d$
of arbitrary dimension, 
it is proven in \cite{GZ}
that $\ev_{0*}(F)$ is $\Theta$-directional, where $\Theta$
is the complete simplicial fan in $\RR^r$ consisting of all
cones
generated by at most $r$ of the vectors
$-e_0, -e_1, \ldots, -e_r$.
We conjecture that a similar proof
also works for Psi-conditions
at point-conditions. 
\end{remark}

\section {WDVV equations and topological recursion} \label{wdvvequations}


We are now ready to prove the tropical analogues of the
WDVV and topological recursion equations --- under certain
assumptions. With the help of these equations, we show that
certain tropical gravitational descendants coincide with their
classical counterparts. This reduces the computation of the
classical invariants to counting problem for tropical curves
with certain valence and incidence conditions 
(cf.\ remark \ref{enumerativerelevance}).

\subsection{WDVV equations}


Let $x_i,x_j,x_k,x_l$ be pairwise different marked 
leaves and consider the forgetful map $\ft : \M{n}{r}{\Delta}
\rightarrow \calM_{\{i,j,k,l\}}$.

\begin{lemma} \label{M4-coordinate}
  The equation
  $$
    \ft^*(\varphi_{\{i,j\}|\{k,l\}}) = 
      \sum_{\substack{I|J \\ i,j \in I, k,l \in J}} \varphi_{I|J}
  $$
  holds, where the sum on the right side runs through \emph{all} 
  (also non-reducible)
  partitions with $i,j \in I$ and $k,l \in J$.
\end{lemma}

\begin{proof}
  Note that $\ft(V_{I|J}) = V_{I \cap \{i,j,k,l\}|J \cap \{i,j,k,l\}}$.
  Therefore $\varphi(\ft(V_{I|J})) = 1$ if $i,j \in I, k,l \in J$ and zero 
  otherwise.
\end{proof}

Now we face the crucial difference to the conventional setting. The right sum
also runs over non-reducible partitions, which do not 
correspond to something in the algebro-geometric case. Let us add up only those 
$\varphi_{I|J}$ with $I|J$ non-reducible and denote the sum
by $\phi$, i.e.
$$
  \phi_{i,j|k,l} := \sum_{\substack{I|J \text{ non-red.} \\ 
                i,j \in I, k,l \in J}} 
          \varphi_{I|J}
$$
We would like to show that $\phi_{i,j|k,l}$ is bounded, as
then it does not change the degree of a intersection product
and we can derive the same formulas as in the conventional case.
So let us investigate what this function measures. \\
  Let $F = (\prod_{k=1}^n 
  \tau_{a_k}(C_k))_\Delta$ be a 
  one-dimensional family of curves with general conditions.
  Consider a facet $\sigma$ of $F$ representing curves \emph{with} contracted bounded
  edge $E$ (called \emph{reducible curves}). Then we can change the length of $E$
  while keeping all other lengths and our curve will still match the incidence 
  conditions. As our conditions are general,
  the set of curves fulfilling the incidence conditions set-theoretically
  is also one-dimensional. Hence, all curves in $\sigma$ just differ
  by the length of $E$, whereas all other lengths are fixed. But this means that
  $\phi_{i,j|k,l}$ is constant on $\sigma$. \\
  Now, let $\sigma$ be a facet of $F$ representing curves \emph{without} 
  contracted bounded edge $E$ (called \emph{non-reducible curves}). This means,
  for all non-reducible partitions $I|J$, the respective function $\varphi_{I|J}$
  is identically zero on $\sigma$. Therefore, on $\sigma$, $\phi_{i,j|k,l}$ coincides with 
  $\ft^*(\varphi_{\{i,j\}|\{k,l\}})$.

\begin{lemma}
  Let $F = (\prod_{k=1}^n 
  \tau_{a_k}(C_k))_\Delta$ be a 
  one-dimensional family of curves with general conditions.
  Let $\sigma$ be a facet of $F$.
	Then
  $$
    \phi_{i,j|k,l}|_\sigma = \left\{ \begin{array}{ll}
                     \varphi_{\{i,j\}|\{k,l\}} \circ \ft   
                       & \text{if interior curves of $\sigma$ are non-reducible} \\		
                     \text{const} & \text{otherwise}.
              \end{array} \right.  
  $$
\end{lemma}

In other words: Proving that $\phi_{i,j|k,l}$ is bounded on a family one-dimensional family
$F$ is the same as \emph{proving that curves in $F$ with large $\calM_{i,j,k,l}$-coordinate
must contain a contracted bounded edge}. This is the way of speaking in existing literature
(e.g. \cite[proposition 5.1]{GM05}, \cite[proposition 6.1]{KM06}, \cite[section 4]{MR08}).
We will address this problem in its own subsection and first state the desired
results here.

\begin{lemma}
\label{IntersectionWithM4-coordinate}
  Let $F = (\prod_{k=1}^n 
  \tau_{a_k}(C_k))_\Delta$
  be a one-dimensional family of curves. 
  Furthermore assume that $\phi_{i,j|k,l}$ is bounded.
  Then the equation
  $$
    \langle \ft^*(\varphi_{\{i,j\}|\{k,l\}}) \cdot 
      {\textstyle \prod\limits_{k=1}^n} \tau_{a_k}(C_k)\rangle_\Delta
      =
      \sum_{\substack{I|J \text{ reducible} \\ i,j \in I, k,l \in J}}
      \langle \varphi_{I|J} \cdot 
      {\textstyle \prod\limits_{k=1}^n} \tau_{a_k}(C_k)\rangle_\Delta
  $$
  holds.
\end{lemma}

\begin{proof}
  This follows from \ref{M4-coordinate} and the fact
  that the degree of a bounded function intersected
  with a one-dimensional cycle is zero. Therefore, if $\phi_{i,j|k,l}$
  is bounded, the degree of 
  $$
    \langle \phi_{i,j|k,l} \cdot
      \prod_{k=1}^n \tau_{a_k}(C_k)\rangle_\Delta
  $$
  is zero and hence this term can be omitted.
\end{proof}

Finally, we can state the following version of the WDVV equations. A more restrictive
version was proven in \cite[Theorem 8.1]{MR08}. 
Let us emphasize again the difference of the two approaches. 
In \cite{MR08}, similar to previous works such as \cite{GM05}, 
the proof of certain WDVV equations was based on two steps. 
First, under generic conditions it is shown explicitly that the curves
under consideration split into two parts. 
Second, it is shown that the multiplicity of the big curve factors
as a product of the two smaller parts. 
This is done by an involved computation in terms
of a suitable matrix representation of 
$\ev_1 \times \cdots \ev_n \times \ft$ (cf.\ \cite[Lemma 6.6]{MR08}). 
In the present approach, these ad hoc computations
are replaced by intersection-theoretic arguments (e.g.\
the splitting lemma).

As before, we fix a complete simplicial fan $\Theta$ 
and a basis $B_0, \ldots, B_m$  of $Z_*(\Theta)$.
Furthermore, let $(\beta_{ef})_{ef}$ be the inverse matrix (over $\QQ$) of
the matrix $(\deg(B_e \cdot B_f))_{ef}$.

\begin{theorem}[WDVV equations
] \label{TheoremWDVVEquations}
  Let $F = (\prod_{k=1}^n 
  \tau_{a_k}(C_k))_\Delta$
  be a one-dimensional family of curves and fix
  four pairwise different marked leaves $x_i, x_j, x_k, x_l$.
  Moreover, we assume that the following conditions hold.
  \begin{enumerate}
    \item 
      For any reducible partition $I|J$ with $i,j \in I; k,l \in J$
      or $i,k \in I; j,l \in J$ the push-forwards
      $\ev_x(X_I)$ and $\ev_y(X_J)$ are $\Theta$-directional
      (with notations from section \ref{splittingcurves}).
    \item
      The functions $\phi_{i,j|k,l}$ and $\phi_{i,k|j,l}$
      are bounded on $F$.
  \end{enumerate}
  Then the WDVV equation
  $$
    \sum_{\substack{I|J \text{ reducible} \\ i,j \in I, k,l \in J}}
      \sum_{e,f}
      \langle{\textstyle \prod\limits_{k \in I}} \tau_{a_k}(C_k) \cdot \tau_0(B_e)\rangle_{\Delta_I}
      \; \beta_{ef} \;
      \langle\tau_0(B_f) \cdot {\textstyle \prod\limits_{k \in J}} \tau_{a_k}(C_k)\rangle_{\Delta_J}
  $$
  $$  =
    \sum_{\substack{I|J \text{ reducible} \\ i,k \in I, j,l \in J}} 
      \sum_{e,f}
      \langle{\textstyle \prod\limits_{k \in I}} \tau_{a_k}(C_k) \cdot \tau_0(B_e)\rangle_{\Delta_I}
      \; \beta_{ef} \;
      \langle\tau_0(B_f) \cdot {\textstyle \prod\limits_{k \in J}} \tau_{a_k}(C_k)\rangle_{\Delta_J}
  $$
  holds,
  where the sums run through reducible partitions only.
\end{theorem}

\begin{proof}
  The statement follows from 
  \ref{IntersectionWithM4-coordinate}
  and the fact that on $\calM_{\{i,j,k,l\}}$ the functions
  $\varphi_{\{i,j\}|\{k,l\}}$ and $\varphi_{\{i,k\}|\{j,l\}}$
  are rationally equivalent. In fact, they
  only differ by a linear function and therefore have the same
  divisor, namely the single vertex in $\calM_{\{i,j,k,l\}}$.
\end{proof}

\begin{remark} \label{LabelledVsUnlabelledWDVV}
  In the algebro-geometric version of these equations (cf.\ \cite[equation (54) and (55)]{FP}) 
	the big sum(s) usually run like
  $\sum_{\beta_1, \beta_2} \sum_{A,B}$, where $\beta_1, \beta_2$ are cohomology classes
  such that $\beta_1 + \beta_2 = \beta$ and $A \dcup B = [n]$ is a partition of the marks.
  We can proceed accordingly and let our sum run through unlabelled instead of labelled degrees,
  as unlabelled degrees correspond via Minkowski weights to cohomology classes.
  If we collect all reducible partitions $I \dcup J = \Delta \dcup [n]$,
  such that the unlabelled degrees $\delta(\Delta_I),\delta(\Delta_J)$ coincide, 
  we obtain a class of $\frac{\Delta!}{\Delta_I! \cdot \Delta_J!}$ elements.
  On the other hand, as mentioned at the beginning of section 
  \ref{parametrizedcurves}, counting curves with labelled non-contracted leaves
  leads to an overcounting by the factor $\Delta!$, i.e. if $\delta := \delta(\Delta)$
  is an unlabelled degree, we should define
  $$
    \langle\prod_{k=1}^n \tau_{a_k}(C_k)\rangle_\delta := 
        \frac{1}{\Delta!} \langle\prod_{k=1}^n \tau_{a_k}(C_k)\rangle_\Delta.
  $$
  So by switching to ``unlabelled'' invariants, the above factor 
  $\frac{\Delta!}{\Delta_I! \cdot \Delta_J!}$ cancels
  and we obtain
  $$
    \sum_{\substack{\delta_I, \delta_J \\ \delta_I + \delta_J = \delta}}
      \sum_{\substack{A \dcup B = [n] \\ i,j \in A, k,l \in B}}
      \sum_{e,f}
      \langle{\textstyle \prod\limits_{k \in A}} \tau_{a_k}(C_k) \cdot \tau_0(B_e)\rangle_{\delta_I}
      \; \beta_{ef} \;
      \langle\tau_0(B_f) \cdot {\textstyle \prod\limits_{k \in B}} \tau_{a_k}(C_k)\rangle_{\delta_J}
  $$
  $$  =
    \sum_{\substack{\delta_I, \delta_J \\ \delta_I + \delta_J = \delta}}
      \sum_{\substack{A \dcup B = [n] \\ i,k \in A, j,l \in B}}
      \sum_{e,f}
      \langle{\textstyle \prod\limits_{k \in A}} \tau_{a_k}(C_k) \cdot \tau_0(B_e)\rangle_{\delta_I}
      \; \beta_{ef} \;
      \langle\tau_0(B_f) \cdot {\textstyle \prod\limits_{k \in B}} \tau_{a_k}(C_k)\rangle_{\delta_J},
  $$
  which is now combinatorially identical to the algebro-geometric version.
\end{remark}

\subsection{Topological recursion}


In the same flavour as in the previous subsection, we will also formulate a
tropical version of the equations known as ``topological recursion''.

Let $x_i,x_k,x_l$ be pairwise different marked leaves.
We know from \ref{PsiAsBoundary} that
we can express the Psi-divisor $\psi_i$ in terms of ``boundary''
divisors, namely
$$
  \divisor(\psi_i) = \sum_{\substack{I|J \\ i \in I, k,l \in J}} \divisor(\varphi_{I|J}).
$$
Now again we give a name to the term that has no algebro-geometric counterpart
$$
  \phi_{i|k,l} = \sum_{\substack{I|J \text{ non-red.} \\ i \in I; k,l \in J}} \varphi_{I|J}.
$$
As in the previous subsection, we can describe this function as follows.

\begin{lemma}
  Let $F = (\prod_{k=1}^n 
  \tau_{a_k}(C_k))_\Delta$ be a 
  one-dimensional family of curves with general conditions.
  Let $\sigma$ be a facet of $F$.
	Then
  $$
    \phi_{i|k,l}|_\sigma = \left\{ \begin{array}{ll}
                     \sum \text{length of edge that separates $i$ from $k,l$}
                       & \text{if interior curves of $\sigma$ are non-reducible,} \\		
                     \text{constant} & \text{otherwise}.
              \end{array} \right.  
  $$
\end{lemma}

Again, we fix a complete simplicial fan $\Theta$ 
and a basis $B_0, \ldots, B_m$  of $Z_*(\Theta)$.
Furthermore, let $(\beta_{ef})_{ef}$ be the inverse matrix (over $\QQ$) of
the matrix $(\deg(B_e \cdot B_f))_{ef}$.

\begin{theorem}[Topological recursion
] \label{TopologicalRecursion}
  Let $F = (\prod_{k=1}^n 
  \tau_{a_k}(C_k))_\Delta$
  be a one-dimensional family of curves 
  and fix
  three pairwise different marked leaves $x_i, x_k, x_l$.
  Moreover, we assume that the following conditions hold.
  \begin{enumerate}
    \item 
      For any reducible partition $I|J$ with $i \in I; k,l \in J$
      the push-forwards
      $\ev_x(X_I)$ and $\ev_y(X_J)$ are $\Theta$-directional
      (with notations from section \ref{splittingcurves}).
    \item
      The function $\phi_{i|k,l}$ is bounded on $F$.
  \end{enumerate}
  Then the 
  topological recursion
  $$
    \langle \psi_i \cdot {\textstyle \prod\limits_{k=1}^n} \tau_{a_k}(C_k) \rangle_\Delta
  =
    \sum_{\substack{I|J \text{ reducible} \\ i \in I, k,l \in J}}
      \sum_{e,f}
      \langle{\textstyle \prod\limits_{k \in I}} \tau_{a_k}(C_k) \cdot \tau_0(B_e)\rangle_{\Delta_I}
      \; \beta_{ef} \;
      \langle\tau_0(B_f) \cdot {\textstyle \prod\limits_{k \in J}} \tau_{a_k}(C_k)\rangle_{\Delta_J}
  $$
  holds,
  where the sum runs through reducible partitions only.
\end{theorem}

\begin{remark} \label{LabelledVsUnlabelledTopRec}
  In the same way as in \ref{LabelledVsUnlabelledWDVV} we obtain
  the ``unlabelled'' version
  $$
    \langle \psi_i \cdot {\textstyle \prod\limits_{k=1}^n} \tau_{a_k}(C_k) \rangle_\delta
  =
    \sum_{\substack{\delta_I, \delta_J \\ \delta_I + \delta_J = \delta}}
      \sum_{\substack{A \dcup B = [n] \\ i \in A, k,l \in B}}
      \sum_{e,f}
      \langle{\textstyle \prod\limits_{k \in A}} \tau_{a_k}(C_k) \cdot \tau_0(B_e)\rangle_{\delta_I}
      \; \beta_{ef} \;
      \langle\tau_0(B_f) \cdot {\textstyle \prod\limits_{k \in B}} \tau_{a_k}(C_k)\rangle_{\delta_J},
  $$  
  which coincides combinatorially with the algebro-geometric version of this equation.
\end{remark}

\subsection{Contracted bounded edges} \label{contractededges}


As a preparation for the more difficult case of plane curves,
we first assume $r=1$. 

\begin{proposition} \label{PhiBoundedForR=1}
  Let $P_1, \ldots, P_n$ be points in general position in $\RR^1$ and
  let $F = (\prod_{k=1}^n \tau_{a_k}(P_k))_d^{\RR^1}$ be
  a one-dimensional family in $\M{n}{1}{d}$. Then for any choice
  of marked leaves $x_i,x_j,x_k,x_l$, the functions
  $\phi_{i,j|k,l}$ and $\phi_{i|k,l}$ are bounded on $F$.
\end{proposition}

\begin{proof}
  For general conditions, $F$ set-theoretically coincides with the
  set of curves satisfying the given 
  incidence and valence conditions. 
  Consider a general curve $C \in F$. Then $C$ is also a general
  curve in the Psi-product $X := \prod_{k=1}^n \psi_k^{a_k}$. As
  we cut down $X$ by $n$ point conditions and $\dim(F) =1$, 
  the dimension of $X$ must be $n+1$, hence $C$ contains $n$ bounded edges.
  This implies that $C$, as it is a rational curve, has $n+1$ vertices.
  Therefore there exists a vertex $V$ not adjacent to
  a marked leaf $x_k, k \in [n]$. Now one of the three
  adjacent
  edges might be a contracted bounded edge. Then the deformation of $C$ in $F$ is
  given by changing the length of this edge, but this does not affect
  $\phi_{i,j|k,l}$ or $\phi_{i|k,l}$ by definition. Otherwise, if
  all of the adjacent edges are non-contracted, the deformation
  of $C$ in $F$ is given by moving $V$ (and changing the lengths
  accordingly).
  \begin{center}
    \input{pics/MovingVertex.pstex_t}
  \end{center}
  Note that the edge $v$ cannot be unbounded as its direction ``vector''
  is not  primitive. Therefore, if this deformation is supposed to
  be unbounded, $v_1,v_2$ must be unbounded. But in this case
  only the length of $v$ grows infinitely. But as $v$ does not separate
  any marked leaves, this does not change $\phi_{i,j|k,l}$ and $\phi_{i|k,l}$. 
\end{proof}

  Now let us consider the case of plane curves, i.e. $r=2$. The whole
  subsection should be compared with \cite[Section 4]{MR08}, where 
  we dealt with the special case $\Delta = d$.
  We fix the following notation. Let
  $F = (\prod_{k=1}^n 
  \tau_{a_k}(C_k))_\Delta^{\RR^2}$
  be a one-dimensional family of plane curves with general conditions
  and and let $L \dcup M \dcup N = [n]$ be the partition
  of the labels such that
  $$
    \codim(C_k) = \begin{cases}
                  0 & \text{if } k \in L, \\
                  1 & \text{if } k \in M, \\
                  2 & \text{if } k \in N. \\
                \end{cases}
  $$
  First we study how the deformation of a general curve $C$ in 
  $F$ can look like.
  
\begin{lemma}[Variation of \cite{MR08} 4.4] \label{deformations}
  Let us assume 
  \begin{enumerate}
    \item[i)]
      $a_k = 0$ for all
      $k \in L \cup M$, i.e. Psi-conditions are
      only allowed together with point conditions.
  \end{enumerate}
  Then the following holds. \\
  Let $\sigma$ be a facet of $F$ and let $C \in \sigma$
  be a general curve. Then
  the deformation of $C$ inside $\sigma$ is
  described by one of the following cases.
  \begin{enumerate}
    \item[\textbf{(I)}]
      $C$ contains a \emph{contracted bounded edge}. Then the deformation
      inside $\sigma$ is given by changing the length of this
      edge arbitrarily.
    \item[\textbf{(II)}]
      $C$ has a three-valent \emph{degenerated vertex} $V$ 
       of one of the following
      three types.
      \begin{enumerate}
	      \item
	        One of the adjacent edges is a marked leaf $i\in  L$.
        \item
          One of the adjacent edges is a marked leaf $j \in M$
          and the linear spans of the corresponding line $C_j$ at
          $\ev_j(C)$ and of the other two edges adjacent to $V$
          coincide (i.e. the curves $C$ and the $C_j$ do \emph{not}
          intersect transversally at $\ev_j(C)$).
        \item
          All edges adjacent to $V$ are non-contracted, but their
          span near $V$ is still only one-dimen\-sional; w.l.o.g.\ we denote
          the edge alone on one side of $V$ by $v$ and the two edges on
          the other side by $v_1, v_2$.
      \end{enumerate}  
      \begin{center}
        \input{pics/DegeneratedVertex.pstex_t}
      \end{center}
      In all these cases the deformation inside $\sigma$ is given by moving $V$.
    \item[\textbf{(III)}]
      $C$ contains a \emph{movable string $S$}, i.e. a two-valent subgraph
      of $C$ homeomorphic to $\RR$ such that all edges are non-contracted and
      all vertices of $S$ are three-valent in $C$ and not degenerated in the
      sense of case (II).
      Then the deformation of $C$ is given by moving $S$ while all vertices
      not contained in $S$ remain fixed (in particular, only edges in or
      adjacent to $S$ change their lengths).
  \end{enumerate}
\end{lemma}

\begin{proof}
  Again, for general conditions, $F$ set-theoretically coincides with the set of
  curves satisfying the given incidence and valence conditions. Thus finding
  the deformation of $C$ inside $\sigma$ is the same as finding a way of changing
  the position and the length of the bounded edges of $C$ such that 
  the resulting curve still meets the incidence conditions $C_k$. \\
  It is obvious that in the cases (I) and (II)
  changing the length of the contracted bounded edge respectively 
  moving the degenerated vertex $V$ leads to such deformations. \\
  In case (III) the non-degeneracy of the vertices makes sure that both
  ends of $S$ consist of non-contracted ends and that a small movement
  of one of these ends leads to a well-defined movement of the whole string
  (a more detailed description can be found in the proof of
  \cite[4.4]{MR08}). \\
  Finally, this list of cases is really complete, as $C$ 
  always contains a string whose vertices are three-valent in $C$ 
  and whose ends are either non-contracted leaves or marked leaves in $L$.
  This follows from the same calculation
  as in \cite[4.3]{MR08}, with the only difference that we have to
  replace the number $3d$ by $\#\Delta$.
\end{proof}

We have now seen how a general curve $C \in F$ can be deformed. In a second step,
we will now focus on unbounded deformations.

\begin{definition} \label{REMstronglyunimodular}
  A complete fan $\Theta$ in $\RR^2$ is called 
  \emph{del Pezzo} 
	if the associated toric surface is a smooth del Pezzo surface.
	Here is a complete list, up to the action of $\text{SL}(2,\ZZ)$.
	\begin{center}
    \begin{picture}(0,0)%
\includegraphics{pics/Unim_f.pstex}%
\end{picture}%
\setlength{\unitlength}{3947sp}%
\begingroup\makeatletter\ifx\SetFigFont\undefined%
\gdef\SetFigFont#1#2#3#4#5{%
  \reset@font\fontsize{#1}{#2pt}%
  \fontfamily{#3}\fontseries{#4}\fontshape{#5}%
  \selectfont}%
\fi\endgroup%
\begin{picture}(1077,936)(-14,-2098)
\put(  1,-1261){\makebox(0,0)[lb]{\smash{{\SetFigFont{8}{9.6}{\familydefault}{\mddefault}{\updefault}$\Theta_{\PP^2}$}}}}
\end{picture}%
 \qquad\qquad
    \begin{picture}(0,0)%
\includegraphics{pics/Unim_c1.pstex}%
\end{picture}%
\setlength{\unitlength}{3947sp}%
\begingroup\makeatletter\ifx\SetFigFont\undefined%
\gdef\SetFigFont#1#2#3#4#5{%
  \reset@font\fontsize{#1}{#2pt}%
  \fontfamily{#3}\fontseries{#4}\fontshape{#5}%
  \selectfont}%
\fi\endgroup%
\begin{picture}(1227,936)(-14,-2098)
\put(  1,-1261){\makebox(0,0)[lb]{\smash{{\SetFigFont{8}{9.6}{\familydefault}{\mddefault}{\updefault}$\Theta_{\PP^1 \times \PP^1}$}}}}
\end{picture}%
 \\[5ex]
    \begin{picture}(0,0)%
\includegraphics{pics/Unim_c2.pstex}%
\end{picture}%
\setlength{\unitlength}{3947sp}%
\begingroup\makeatletter\ifx\SetFigFont\undefined%
\gdef\SetFigFont#1#2#3#4#5{%
  \reset@font\fontsize{#1}{#2pt}%
  \fontfamily{#3}\fontseries{#4}\fontshape{#5}%
  \selectfont}%
\fi\endgroup%
\begin{picture}(1227,936)(-14,-2098)
\put(  1,-1261){\makebox(0,0)[lb]{\smash{{\SetFigFont{8}{9.6}{\familydefault}{\mddefault}{\updefault}$\Theta_{\FF_1}$}}}}
\end{picture}%
 \qquad\qquad
    \begin{picture}(0,0)%
\includegraphics{pics/Unim_d.pstex}%
\end{picture}%
\setlength{\unitlength}{3947sp}%
\begingroup\makeatletter\ifx\SetFigFont\undefined%
\gdef\SetFigFont#1#2#3#4#5{%
  \reset@font\fontsize{#1}{#2pt}%
  \fontfamily{#3}\fontseries{#4}\fontshape{#5}%
  \selectfont}%
\fi\endgroup%
\begin{picture}(1227,936)(-14,-2098)
\put(  1,-1261){\makebox(0,0)[lb]{\smash{{\SetFigFont{8}{9.6}{\familydefault}{\mddefault}{\updefault}$\Theta_{\Bl_2(\PP^2)}$}}}}
\end{picture}%
 \qquad\qquad
    \begin{picture}(0,0)%
\includegraphics{pics/Unim_e.pstex}%
\end{picture}%
\setlength{\unitlength}{3947sp}%
\begingroup\makeatletter\ifx\SetFigFont\undefined%
\gdef\SetFigFont#1#2#3#4#5{%
  \reset@font\fontsize{#1}{#2pt}%
  \fontfamily{#3}\fontseries{#4}\fontshape{#5}%
  \selectfont}%
\fi\endgroup%
\begin{picture}(1227,936)(-14,-2098)
\put(  1,-1261){\makebox(0,0)[lb]{\smash{{\SetFigFont{8}{9.6}{\familydefault}{\mddefault}{\updefault}$\Theta_{\Bl_3(\PP^2)}$}}}}
\end{picture}%
 
  \end{center}
  It is easy to see that an alternative characterization of these
	fans is as follows: \emph{Any} two independent
  primitive vectors generating rays of $\Theta$ form a basis of $\ZZ^2$.
  A degree $\Delta$ in $\RR^2$ is called \emph{del Pezzo} if $\Theta(\Delta)$
  is del Pezzo and if all direction vectors appearing in $\Delta$ are
  primitive. This ensures that for every
  pair of independent vectors $v_1, v_2$ appearing in $\Delta$, the dual triangle to the
  fan spanned by $v_1$, $v_2$ and $-(v_1 + v_2)$ does not contain lattice points apart
	from its vertices.
  %
\end{definition}
%

\begin{lemma}[Variation of \cite{MR08} 4.4] \label{UnboundedDeformations}
  We assume 
  \begin{enumerate}
    \item[i)]
      $a_k = 0$ for all $k \in L \cup M$,
    \item[ii)]
       $\Delta$ is del Pezzo.
  \end{enumerate}
  Then the following holds. \\
  Let $\sigma$ be a \emph{unbounded} facet of $F$ and let $C \in \sigma$
  be a general curve. Then
  the deformation of $C$ in $\sigma$ is
  described by one of the following cases.
  \begin{enumerate}
    \item[\textbf{(I)}]
      $C$ contains a contracted bounded edge whose length 
      can be changed arbitrarily.
    \item[\textbf{(II)}]
      $C$ has a three-valent degenerated vertex $V$ of one the three types
      described above. Furthermore, in the cases (a) and (b) 
      (of \ref{deformations} (II))
      one of the edges $v_1, v_2$ is bounded, the other one unbounded,
      whereas in case (c) the edge $v$ is bounded and $v_1, v_2$ are
      unbounded.
    \item[\textbf{(III)}]
      $C$ contains a movable string $S$ with two non-contracted
      leaves $v_1, v_2$ and \emph{only one adjacent bounded edge $w$}. The
      deformation of $C$ is given by increasing the length of $w$.
      \begin{center}
        \input{pics/caseIII.pstex_t}
      \end{center}
      Furthermore, if $x_k, k \in M$ is a marked leaf adjacent to $S$, then
      $h(x_k)$ is a general point in an unbounded facet of $C_k$ whose outgoing
      direction vector $v$ lies in the interior of the cone spanned by 
      $v_1, v_2$.
  \end{enumerate}
\end{lemma}

\begin{proof}
  Nothing happens in the cases (I), (II) (a) and (b). 
  In case (II) (c), the edge
  $v$ cannot be unbounded as $v = -v_1 - v_2$ is not primitive.
  Therefore the two edges on the other side of $V$ must be unbounded. 
  
  In case (III), the proof of the first statement is fully contained
  in the last part of the proof of \cite[4.4]{MR08}.
  We assume that we have a string $S$ with two unique non-contracted
	ends and all of its vertices are three-valent and not degenerated
	in the sense of case (II). The deformation only
	moves the string $S$; the adjacent edges are shortened or elongated and
	the other parts of the curve remain fixed. We want to show that $S$ has only 
	one adjacent bounded edge. 
	
	If
  there are bounded edges adjacent to $S$ to both sides of $S$ as in
  picture (a) below then the movement of the string is bounded. (This is true because if we move the string to either side, we can only move until the length of one of the adjacent bounded edges shrinks to $0$.) 
So we only have to consider the case when all adjacent bounded edges of $ S $
  are on the same side of $S $, say on the right side as in picture (b) below. Label the edges of $S$ (respectively, their direction vectors) by $ v_1,\dots,v_k $ and the adjacent bounded
  edges of the curve by $ w_1,\dots,w_{k-1} $ as in the picture. As above the
  movement of the string to the right is bounded. If
  one of the directions $ w_{i+1} $ is obtained from $ w_i $ by a left turn (as
  it is the case for $ i=1 $ in the picture) then the edges $ w_i $ and $
  w_{i+1} $ meet on the left of $S$. This restricts the movement of the string to the left, too, since the corresponding
  edge $ v_{i+1} $ then shrinks to length $0$. 
\begin{center}
\input{pics/sides.pstex_t}
\end{center}
  So we can assume that for all $i$ the direction $ w_{i+1} $ is
  either the same as $ w_i $ or obtained from $ w_i $ by a right turn as in
  picture (c). The balancing condition then shows that for all $i$ both the
  directions $ v_{i+1} $ and $ -w_{i+1} $ lie in the angle between $ v_i $ and
  $ -w_i $ (shaded in the picture above). Therefore, all directions $ v_i
  $ and $ -w_i $ lie within the angle between $ v_1 $ and $ -w_1 $. In
  particular, the image of the string $ S $ cannot have any self-intersections in $
  \RR^2 $. We can therefore pass to the (local) dual picture (d)  where the edges dual to $ w_i $ correspond to a
  concave side of the polygon whose other two edges are the ones dual to $ v_1
  $ and $ v_k $.

  But from our assumption that $\Delta$ is del Pezzo we know
  that the triangle dual to $ v_1 $ and $ v_k $ does not contain more integer 
  points than its vertices. We conclude that the concave side of the polygon
  in picture (d) actually must coincide with the 
  triangle dual to $ v_1 $ and $ v_k $ and 
  therefore the string consists of the two
  ends $ v_1 $ and $ v_2 $ that are connected to the rest of the curve
  by exactly one bounded edge $ w_1=w $ (as shown in picture
  (e)). 
  
  The second
  statement concerning adjacent marked leaves $x_k, k \in M$ is obvious
  as the deformation is supposed to be unbounded. 
\end{proof}

\begin{theorem} \label{PhiBounded}
  Let $x_i, x_j, x_k, x_l$ be pairwise different marked leaves
  and let us assume 
  \begin{enumerate}
    \item[i)]
      $a_k = 0$ for all $k \in L \cup M$,
    \item[ii)]
       $\Delta$ is del Pezzo,
    \item[iii)]
      if $i,j \in M$ (resp. $k,l \in M$), then for any pair of
      independent direction vectors $v_1,v_2$ appearing in $\Delta$, 
      the interior of the cone
      spanned by $v_1,v_2$ does not intersect both degrees
      $\delta(C_i)$
      and $\delta(C_j)$ (resp. $\delta(C_k)$
      and $\delta(C_l)$).
  \end{enumerate}
  Then $\phi_{i,j|k,l}$ is bounded. \\
  If we additionally require
  \begin{enumerate}
    \item[iv)]
      $i \in N$,
  \end{enumerate}
  then also $\phi_{i|k,l}$ is bounded. 
\end{theorem}

\begin{proof}
  As conditions i) and ii) hold, we can apply \ref{UnboundedDeformations}, which
  describes the unbounded facets of $F$. We have to show that $\phi_{i,j|k,l}$
  (resp. $\phi_{i|k,l}$) is bounded on these facets. 
  In case (I), the only changing length is that of a contracted edge and therefore
  not measured by both $\phi_{i,j|k,l}$ and $\phi_{i|k,l}$. 
  In case (II), the edge whose length is growing infinitely cannot separate more
  then one marked leaf $x_k, k \in L \cup M$ from the others. 
  Therefore this length cannot contribute
  to $\phi_{i,j|k,l}$ and --- by condition iv) --- to $\phi_{i|k,l}$. 
  Finally, condition iii) (and also condition iv)) is made such that 
  $\phi_{i,j|k,l}$ and $\phi_{i|k,l}$ are also bounded in case (III).
\end{proof}

\begin{remark} \label{CounterexamplesForContractedEdge}
  The conditions i) -- iv) appearing in the above statements are not only sufficient
  but, in most cases, also necessary for the statements to hold.
  \begin{enumerate}
    \item[iv)]
      If condition iv) in \ref{PhiBounded} is not satisfied,
      we can get the following things. 
      \begin{itemize}
        \item
          If $i \in L$, then
          the degenerated vertex of type (a) leads to an
          unbounded $\phi_{i|k,l}$.
        \item
          If $i \in M$ and $\rho$ is a ray in $C_i$ whose direction
          vector $v_\rho$ also appears in $\Delta$, then
          in general we will find curves in $F$ with a degenerated 
          vertex of type (b), whose unbounded movement will 
          make $\phi_{i|k,l}$ unbounded.
        \item
          If $i \in M$ and $\rho$ is a ray in $C_i$ whose direction
          vector $v_\rho$ lies between two direction vectors $v_1, v_2$
          appearing in $\Delta$, this will in general lead to curves
          in $F$ with unbounded deformations of case (III) such that the outward
          directions are $v_1,v_2$ and such that $x_i$ is adjacent to the moved
          string. So again, $\phi_{i|k,l}$ is in general unbounded.
      \end{itemize}  
    \item[iii)]
      If condition iii) is not satisfied,
      we will in general get unbounded deformations of the following
      type.
      \begin{center}
        \input{pics/CounterexCaseIII.pstex_t}
      \end{center}
      In this case we have $i,j \in M$ and the interior of the cone
      spanned by $v_1,v_2$ contains direction vectors of
      both $C_i$ and $C_j$. As in general $x_k,x_l$ will lie on the
      other side of the growing edge $w$, $\phi_{i,j|k,l}$ will be unbounded.
    \item[ii)]
       If we drop condition ii), i.e. if we allow non-del Pezzo degrees
       $\Delta$, two things can happen.
       If we allow non-primitive direction vectors, then 
       we get deformations of type (II) (c) with unbounded
       edge $v$. Therefore the lengths of $v_1$ and $v_2$, which can
       in general separate arbitrary marked leaves, grow infinitely.
       If $\Theta(\Delta)$ is not supposed to be del Pezzo, 
       then the description of unbounded deformations of case (III)
       in \ref{UnboundedDeformations} becomes incorrect, as there will appear
       more complicated strings with more adjacent bounded edges
       than just one. The example of $\FF_2$ is analysed in detail in
       \cite{Fra} and \cite[e.g. 2.10]{FM}.
    \item[i)]
      If we drop condition i), i.e. if we allow Psi-conditions also
      at marked leaves which are not fixed by points, we end up
      with more complicated kinds of deformations
      of general curves in $F$. The following picture shows
      an example of an unbounded deformation in a one-dimensional
      family of plane curves of projective degree $2$.
      \begin{center}
        \input{pics/CounterexPsiAtLines.pstex_t}
      \end{center}
      Here, $C$ has to meet all the four tropical lines $C_1, \ldots, C_4$ 
      with one Psi-condition. Note that the indicated deformation of $C$
      is indeed unbounded and that the length of the $(1,-1)$-edge $e$ grows
      infinitely. This example can be extended in the following way.
      One can glue
      arbitrary (fixed) curves to the non-contracted leaves of $C$ in direction
      $(1,1)$, obtaining more families admitting such a deformation. In particular,
      the edge $e$ can separate arbitrary kinds of points,
      showing that in general $\phi_{i,j|k,l}$ and $\phi_{i|k,l}$ can be unbounded
      for any choice of $i,j,k,l$.
  \end{enumerate}
\end{remark}

In higher dimensions ($r \geq 3$), up to now only the following case is studied.

\begin{theorem}[\cite {Zim} 4.86] \label{PhiBoundedForHigherR}
  Let $F = (\prod_{k=1}^n \tau_0(V_k))_d^{\RR^r}$ be a one-dimensional
  family of curves of projective degree $d$ in $\RR^r$ which
  do \emph{not} satisfy Psi-conditions, but incidence conditions
  given by conventional
  linear spaces $V_k \subseteq \RR^r$. Then for any choice of
  $\{i,j,k,l\} \in [n]$ the function $\phi_{i,j|k,l}$
  is bounded on $F$.
\end{theorem}

\subsection{Comparison to the algebro-geometric invariants} 
                                      \label{compofnumbers}

In the special case of an empty degree, denoted by $\Delta = 0$,
the situation is analogous to the algebro-geometric one.

\begin{proposition} \label{DegreeZero}
  Let $Z = (\prod_{k=1}^n \tau_{a_k}(C_k))_0$ be a zero-dimensional
  intersection product in $\M{n}{r}{0}$. Then $\deg(Z)$ is non-zero
  if and only if $\sum_{k=1}^n \codim(C_k) = r$ (or equivalently
  $\sum_{k=1}^n a_k = n-3$). In this case, 
  $$
    \deg(Z) = \binom{n-3}{a_1, \ldots, a_n} 
              \deg(C_1 \cdots C_k)
  $$
  holds.
\end{proposition}

\begin{proof}
  By definition $\M{n}{r}{0}$ is isomorphic to
  $\calM_n \times \RR^r$. Moreover, as $\Delta = 0$, all
  evaluation maps $\ev_i$ coincide with the projection
  onto the second factor, which we therefore denote by $\ev$.
  Now let $X := \prod_{k=1}^n \psi_k^{a_k} =
  (\prod_{k=1}^n (\psi_k^{\text{abstr}})^{a_k}) \times \RR^r$
  be the intersection of all Psi-divisors. Then the projection
  formula applied to $\ev$ yields
  $$
    \deg(Z) = \deg(C_1 \cdots C_n \cdot \ev_*(X)).
  $$
  But $\ev_*(X)$ is non-zero if and only if 
  $\sum_{k=1}^n a_k = n-3$. If so, by \ref{AbstractInvariants}
  we know  $\ev_*(X) = \binom{n-3}{a_1, \ldots, a_n} \cdot \RR^r$,
  which proves the statement.
\end{proof}

\begin{remark} \label{TropicalClassicalPositivity}
  The goal of the following theorem is to show that
  certain tropical and classical gravitational descendants
  coincide. The idea is to show that --- under the restrictions
  which we accumulated in the preceding sections --- both 
  sets of numbers satisfy the same WDVV and topological 
  recursion equations, which are sufficient to determine
  the numbers from some initial values. 
  However, there is one further problem concerning
  this plan, which we already mentioned
  in remark \ref{LabelledVsUnlabelledWDVV}.
  The classical WDVV and topological recursion
  equations run through splittings
  of the given cohomology class $\beta$ into
  sums $\beta = \beta_1 + \beta_2$. As
  $\Mbar_{0,n}(\bX, \beta)$ is empty if $\beta$
  is not effective, we can restrict to effective
  classes $\beta, \beta_1, \beta_2$.
  
  Now, for $\PP^2$ and $\PP^1 \times \PP^1$, 
  effectivity is equivalent to the fact that
  the associated one-dimensional tropical fans
  are positive (as $\PP^2$ and $\PP^1 \times  \PP^1$
  do not contain curves with negative self-intersection).
  So a splitting $\beta = \beta_1 + \beta_2$ of
  effective cohomology classes corresponds bijectively
  to a sum of unlabelled tropical degrees
  $\delta = \delta_1 + \delta_2$, and therefore
  the tropical and classical equations are
  really equivalent in this case.
  
  However, for the blow ups of $\PP^2$ in up to
  three torus-fixed points (i.e.\ for
  $\FF_2$, $\Bl_2(\PP^2)$ and $\Bl_3(\PP^2)$,
  cf.\ definition \ref{REMstronglyunimodular}), the
  same argument fails as the exceptional
  divisors induce tropical fans with
  negative weights. The following picture
  shows the example of the
  tropical fan associated to the exceptional
  divisor $V(\ray)$ of $\FF_1$.
  \begin {center} \input {pics/EffectivePositive.pstex_t} \end {center}
  In these cases, i.e.\ when classical curves
	can split into reducible curves with a rigid
	component in the toric boundary, we cannot
	expect that our purely non-compact approach
	will yield the same results. 
	It should be possible to deal with this
	by (partially) compactifying our spaces and/or
	adding suitable correction terms
	(as in \cite{FM}). This needs to be addressed 
  in further work. 
	For now, we just restrict ourselves to $\PP^2$
  and $\PP^1 \times \PP^1$. 
\end{remark}

Now we are finally ready to compare the tropical invariants for plane
tropical curves to the algebro-geometric ones, for some cases,
using the equations proven in 
the previous subsections. The theorem is an extension of \cite[Theorem 8.4]{MR08},
which proves the statement for the case $\PP^2$.

\begin{theorem} \label{CompOfNumbersPlane}
  Let 
  \begin{itemize}
    \item
      $\Fan$ be $\Fan_{\PP^2}$ or $\Fan_{\PP^1 \times \PP^1}$, 
			and set $\bX:=\bX(\Fan)$ (i.e.\ $\bX = \PP^2$ or $\bX = \PP^1 \times \PP^1$),
    \item
      $C_1, \ldots, C_n$ be $\Fan$-directional tropical cycles,
      and let $\gamma_1, \ldots, \gamma_n \in A^*(\bX)$ be the 
      associated cohomology classes of $\bX$, 
    \item
      $\Delta$ be a labelled degree with primitive 
      direction vectors whose unlabelled degree $\delta(\Delta)$
      is $\Fan$-directional, 
      and let $\beta \in A^{r-1}(\bX)$ be the corresponding
      cohomology class,
    \item
      $a_1, \ldots, a_n$ be non-negative integers such that $a_k = 0$ 
      if $\dim(C_k) > 0$.
  \end{itemize}
  Then the tropical and algebro-geometric gravitational descendants
  are equal, i.e.\
  \[
    \frac{1}{\Delta!} \langle \tau_{a_1}(C_1) \cdots \tau_{a_n}(C_n) \rangle_\Delta^{\RR^2} =
      \langle \tau_{a_1}(\gamma_1) \cdots \tau_{a_n}(\gamma_n) \rangle_\beta^\bX.
  \]
\end{theorem}

\begin{proof}
  First we choose a basis $B_0, \ldots, B_m$ of $Z_*(\Fan)$. 
  This also determines a basis $\eta_0, \ldots, \eta_m$ of $A^*(\bX)$, 
  and we know from
  the comparison to the fan displacement rule
  (cf.\ theorem \ref{FanDisplacementEqualsTropInt}) that
  \[
    \deg(B_e \cdot B_f) = \deg(\eta_e \cdot \eta_f)
  \]
  holds. 
  This implies that, if we use WDVV equations or 
  topological recursion with
  respect to these bases,
  then the diagonal coefficients $\beta_{ef}$ appearing 
  in the tropical and in 
  the algebro-geometric setting coincide. 
  Thus, using the results of the previous sections we know that
  the numbers 
  $\frac{1}{\Delta!} \langle \tau_{a_1}(C_1) \cdots $
  $\tau_{a_n}(C_n) \rangle_\Delta
  = \langle \tau_{a_1}(C_1) 
  \cdots $ $\tau_{a_n}(C_n) \rangle_{\delta(\Delta)}$
  and $\langle \tau_{a_1}(\gamma_1) 
  \cdots$ $\tau_{a_n}(\gamma_n) \rangle_\beta^\bX$
  satisfy a certain set of identical equations, namely 
  the WDVV and topological recursion equations 
  (where on the tropical side
  we have to be slightly more careful about $i,j,k,l$ satisfying condition
  iii) and iv) of theorem \ref{PhiBounded}) 
  as well as the string and divisor equation. 
  Therefore we can finish the proof by showing that
  the numbers can be computed recursively, using these equations,
  from some initial numbers and proving that these initial
  numbers coincide. 
  
  We separate the labels of the marked leaves into the sets 
  $L \dcup M \dcup N = [n]$ according
  to the dimension of $C_k$ as 
  in subsection \ref{contractededges}. First we 
  use topological recursion to reduce the number of Psi-conditions: We pick a marked
  leaf $x_i$ with $a_i > 0$ (and therefore $i \in N$) and an arbitrary pair of
  marked leaves $x_k,x_l$ satisfying condition iii) of \ref{PhiBounded}. If such $x_k,x_l$ do
  not exist, we can add them using the divisor equation backwards with appropriate
  rational functions $h_k,h_l$. Namely, if $\bX = \PP^1 \times \PP^1$ we can use
  $h_k = h_l = \max\{0,x,y,x+y\}$, otherwise we can use 
  $h_k = h_l = \max\{0,x,y\}$. Note also that this choice ensures that 
  $h_k \cdot \Delta = h_l \cdot \Delta$ is non-zero for every possible degree, so
  we do not divide by zero. 
  
  After eliminating all Psi-conditions in this way,
  we can assume $a_k = 0$ for all $k \in [n]$, i.e. we are back in the case of
  usual (primary) Gromov-Witten invariants. After applying the
  string and divisor equation
  we can assume that $L = M = \emptyset$ and it remains to
  compute invariants of the form $\langle\prod_{k=1}^n \tau_0(P_k)\rangle_\Delta$
  for points $P_1, \ldots, P_n \in \RR^2$. Comparing dimension shows $\#\Delta = n + 1$.
  Let us first consider the general case $n \geq 3$. Here we consider
  the one-dimensional family 
  $F = (\tau_0(C_i) \tau_0(C_j) \prod_{k=1}^{n-1} \tau_0(P_k))_\Delta$ with
  arbitrary $\Fan$-directional curves $C_i, C_j$ such that $C_i \cdot C_j$
  is non-zero and such that condition iii) of \ref{PhiBounded} is satisfied
  (e.g.\ we can choose the divisors of the functions chosen above). 
  We let $x_i,x_j$ be the first two marked leaves as indicated, and choose 
  $k,l \in [n-1]$ arbitrarily. In the corresponding WDVV equation only
  one extremal partition $I|J$ with $\Delta_I = 0, \Delta_J = \Delta$ does not
  vanish. This follows from lemma \ref{DegreeZero} and the fact
  that the three sums
  $\codim(P_k) + \codim(P_l), \codim(C_i) + \codim(P_k), 
  \codim(C_j) + \codim(P_l)$ are greater than $2$.
  Moreover, the only remaining extremal partition $I = \{i,j\}, J= \Delta \dcup [n-1]$
  provides the term
  \[
    \langle \tau_0(C_i) \tau_0(C_j) \tau_0(\RR^2) \rangle_0 \cdot 
      \langle \tau_0(P) \prod_{k=1}^{n-1} \tau_0(P_k)\rangle_\Delta
    = \deg(C_i \cdot C_j) \cdot \langle\prod_{k=1}^n \tau_0(P_k)\rangle_\Delta.
  \] 
  Hence, we can reduce the computation of 
  $\langle\prod_{k=1}^n \tau_0(P_k)\rangle_\Delta$
  to invariants of smaller degree. We can repeat this until we arrive at the initial invariants
  with $n=1$ or $n=2$. In these cases $\#\Delta = 2$ or $\#\Delta = 3$ and therefore the only
  possible degrees (up to identification via linear isomorphisms of $\ZZ^2$) are
  $\Delta = \{-e_1, e_1\}$ and $\Delta = \{-e_1, -e_2, e_1 + e_2\}$. In both cases, it is easy
  to show by direct computation that $\langle \tau_0(P_1) \rangle_\Delta = 1$ and 
  $\langle \tau_0(P_1) \tau_0(P_2) \rangle_\Delta = 1$ hold (given a point in $\RR^2$,
	there is exactly one horizontal line through it; given two points, there is exactly 
	one tropical line connecting them). 
  But now, as discussed above, the same recursion for the classical
  numbers proves the claim.
\end{proof}

\begin{remark}[Multiplicities of tropical curves] \label{multiplicity}
  The above theorem reduces the computation of the classical 
  gravitational descendants to the count of certain
  tropical curves $C$ with multiplicities $\mult(C)$
  (cf.\ remark \ref{enumerativerelevance}).
  In the above case of plane curves, an easy formula
  for this multiplicity exists (cf.\ \cite[lemma 9.3]{MR08}).
  Namely, if we assume general position, the multiplicity
  of a curve in the count is obtained as the product
  \[
    \mult(C) = \prod_V \mult(V),
  \]
  where the product runs through all vertices to which no
  marked leaf is adjacent and $\mult(V)$ of these
  necessarily $3$-valent vertices is the well-known
  vertex multiplicity introduced by Mikhalkin (cf.\ 
  \cite[definition 2.16]{Mi03}).
  This is correct for labelled curves $C$, but we can
  as well count unlabelled curves $\tilde{C}$ (as the
  incidence and valence conditions do not depend on
  the labelling) by setting
  \[
    \mult(\tilde{C}) =
      \frac{1}{\#\mathrm{Aut}(\tilde{C})} \mult(C).
  \]
  Here $\#\mathrm{Aut}(\tilde{C})$ denotes the number
  of automorphisms of $\tilde{C}$.
  
  Moreover, as well as for the usual Gromov-Witten 
  invariants considered in \cite{Mi03},
  there exists a so-called lattice path algorithm
  to compute these counts easily (cf.\ \cite[section 9]{MR08}).
\end{remark}

\begin{remark} \label{CompOfNumbersR=1}
  Similarly we can deal with the case $r=1$, i.e. we can prove
  $$
    \frac{1}{d!^2} \langle \tau_0(\RR^1)^l \prod_{k=1}^n \tau_{a_k}(P_k) \rangle_d^{\RR^1} =
      \langle \tau_0([\PP^1])^l \prod_{k=1}^n \tau_{a_k}([pt]) \rangle_d^{\PP^1},
  $$
  where the left hand side is a tropical, the right hand side a conventional
  invariant and $[pt]$ denotes the class of a point $pt \in \PP^1$. In fact,
  after applying the string equation, we are left with the case where $l=0$.
  Now we use \ref{PhiBoundedForR=1} and topological recursion to reduce
  the number of Psi-conditions (where, if $n < 3$, we first add more
  marked leaves using the divisor equation). Finally, when $a_k = 0$ for
  all $k \in [n]$, it follows $d=1$ and we can compute directly
  $\langle \tau_0(P) \rangle_1^{\RR^1} = 1$. \\
	This fits with the previously known result for rational Hurwitz numbers
	$H_d^0 := \langle \tau_{1}([pt])^{2d-2} \rangle_d^{\PP^1}$
	(cf. \cite[lemma 9.7]{CJM08}).
\end{remark}

\begin{remark}
  The discussion in \ref{CounterexamplesForContractedEdge} 
  and the factor $n+\#\Delta-2$ appearing in the tropical dilaton
  equation \ref{DilatonEquationForMaps}, instead of $n-2$ in the algebro-geometric
  version, show that for more difficult degrees $\Delta$ (if $r=2$) and for Psi-conditions
  at marked leaves $x_k$ with $\dim(C_k) > 0$, the corresponding tropical and 
  conventional invariants are in general different. 
  For example, if we add a marked leaf that
  has to satisfy only a Psi-condition, the different factors in the dilaton equations 
  immediately lead to different invariants.
\end{remark}

\begin{remark}
  Of course, the machinery developed here is ready to use in higher dimensions as well.
  For example, by remark \ref{DirectionsForHigherR} and theorem 
  \ref{PhiBoundedForHigherR}, the
  same approach can be used to show that tropical and classical
  Gromov-Witten invariants (without Psi-classes) of $\PP^r$,
  $r$ arbitrary, coincide.
\end{remark}


\begin {thebibliography}{MMMM}


\bibitem [AR07]{AR07}
  \arxivjournal{Lars Allermann, Johannes Rau}
        {First steps in tropical intersection theory}
				{Math.\ Z., Volume 264 (2010), no.\ 3, 633--670}
				{0709.3705}
				
\bibitem [AR08]{AR08}
  \arxiv{Lars Allermann, Johannes Rau}
        {Tropical rational equivalence on $\RR^r$}
        {0811.2860}

\bibitem [CJM08]{CJM08}
  \arxivjournal{Renzo Cavalieri, Paul Johnson, Hannah Markwig}
        {Tropical Hurwitz Numbers}
				{J.\ Alg.\ Comb., Volume 32 (2010), no.\ 2}
        {0804.0579}

\bibitem[Fra]{Fra}
  Marina Franz, \emph{The tropical Kontsevich formula
  for toric surfaces}, diploma thesis (2008).
  
\bibitem[FM]{FM}
  \arxivjournal{Marina Franz, Hannah Markwig}
        {Tropical enumerative invariants of $\FF_0$ and $\FF_2$}
				{Adv.\ Geom., Volume  11 (2011), no.\ 1, 49--72}
        {0808.3452}

\bibitem[FP]{FP} 
  \arxivjournal{William Fulton, Rahul Pandharipande}
        {Notes on stable maps and quantum cohomology} 
        {Proc.\ Symp.\ Pure Math.\, Volume 62 (1997), part 2, 45--96}
        {alg-geom/9608011}

\bibitem[FS94]{FS94}
  \arxivjournal{William Fulton, Bernd Sturmfels} 
  	    {Intersection Theory on Toric Varieties}
        {Topology, Volume 36 (1997), no.\ 2, 335--353}
        {alg-geom/9403002}
        

\bibitem [GKM07]{GKM07}
  \arxivjournal{Andreas Gathmann, Michael Kerber, Hannah Markwig}
        {Tropical fans and the moduli spaces of tropical curves}
        {Compos.\ Math., Volume 145 (2009), no.\ 1, 173--195}
        {0708.2268}

\bibitem [GM05]{GM05}
  \arxivjournal{Andreas Gathmann, Hannah Markwig}
        {Kontsevich's formula and the WDVV equations in tropical geometry}
        {Adv.\ Math., Volume 217 (2008), 537--560}
        {math/0509628}
        

\bibitem [GZ]{GZ}
	Andreas Gathmann, Eva-Maria Zimmermann, \emph {The WDVV
  equations in tropical geometry}, in preparation.

\bibitem[H]{H}
  Matthias Herold, \emph{Intersection theory of the tropical
  moduli spaces of curves}, diploma thesis (2007).
  
  
\bibitem [Ka06]{Katz06}
  \arxivjournal{Eric Katz}
        {A Tropical Toolkit}
        {Expo.\ Math., Volume 27 (2009), no.\ 1, 1--36}
        {math/0610878}  

\bibitem [Ka09]{Ka2}
  \arxiv{Eric Katz}
        {Tropical Intersection Theory from Toric Varieties}
        {0907.2488}  

	
  
\bibitem [KM06]{KM06}
  \arxivjournal{Michael Kerber, Hannah Markwig}
        {Counting tropical elliptic plane curves with fixed $j$-invariant}
        {Comment.\ Math.\ Helv., Volume 84 (2009), no.\ 2, 387--427}
        {math/0608472}
        
\bibitem [KM07]{KM07}
  \arxivjournal{Michael Kerber, Hannah Markwig}
        {Intersecting Psi-classes on tropical $M_{0,n}$}
        {Int.\ Math.\ Res.\ Not., Volume 2009, no.\ 2, 221--240}
        {0709.3953}

\bibitem[KV07]{KV}
  Joachim Kock, Israel Vainsencher,
  \emph{An Invitation to Quantum Cohomology},
  Progress in Mathematics, Volume 249 (2007), Birkhäuser Boston.



\bibitem[MR08]{MR08}
  \arxivjournal{Hannah Markwig, Johannes Rau}
        {Tropical descendant Gromov-Witten invariants}
        {Manuscripta Mathematica, Volume 129 (2009), no.\ 3, 293--335}
        {0809.1102}
        
\bibitem[Mi03]{Mi03}
  \arxivjournal{Grigory Mikhalkin}
        {Enumerative tropical geometry in $\RR^2$}
        {J.\ Amer.\ Math.\ Soc., Volume 18 (2005), 313--377}
        {math/0312530}

\bibitem[Mi06]{Mi06}
  \arxivjournal{Grigory Mikhalkin}
        {Tropical geometry and its applications}
        {Int.\ Cong.\ Math., Volume II (2006), 827--852}
        {math/0601041}

\bibitem[Mi07]{Mi07}
  \arxiv{Grigory Mikhalkin}
        {Moduli spaces of rational tropical curves}
        {0704.0839}
  

\bibitem[SS03]{SS03}
  \arxivjournal{David Speyer, Bernd Sturmfels}
        {Tropical Grassmannians}
				{Adv.\ Geom., Volume 4 (2004), no.\ 3, 389--411}
				{math/0304218}


\bibitem[Zim]{Zim}
  Eva-Maria Zimmermann, \emph{Generalizations of the tropical
  Kontsevich formula to higher dimensions},
  diploma thesis (2007).
  
\end {thebibliography}
\end {document}